\def\version{12.9.2023}
\def\users{us}  %
\def\users{final-layout}   % when activated, ``our'' debugging is suppressed
\documentclass[11pt]{article}
\usepackage{upgreek}
\usepackage{xcolor}
\usepackage{bm,amsmath,amsthm,hyperref,amsfonts,amssymb,color}
\usepackage{mathrsfs} % THIS PACKAGE ALLOWS FOR ``scr'' FONTS
\usepackage{cite}

\usepackage{ifthen}
\ifthenelse{\equal{\users}{final-layout}}{}{
\usepackage{fancyhdr}
\pagestyle{fancy}
\headheight=28pt\headwidth=17cm
\definecolor{gray}{gray}{0.5}
\rhead{\color{gray}Thermo-visco-elastodynamics of solids\\
T.Roub\'\i\v cek}
\chead{}
\lhead{Version\,\version, file:\,\jobname.tex
\\
compiled:
\number\day.\number\month.\number\year\ at
\the\hour:\ifnum\minute<10 0\fi\the\minute\ h\ \ \ \ \ 
}
}

\newcount\hour \newcount\minute
\hour=\time
\divide \hour by 60
\minute=\time
\loop \ifnum \minute > 59 \advance \minute by -60 \repeat

\usepackage[normalem]{ulem}

\usepackage[normalem]{ulem}
\usepackage{ifthen}
\usepackage{color}

\ifthenelse{\equal{\users}{final-layout}}{

	\newcommand{\COMMENT}[1]{}
	\newcommand{\COMMENTGT}[1]{}
	\newcommand{\TODO}[1]{}
	\newcommand{\INTERNAL}[1]{}
	\newcommand{\QUESTION}[1]{}
	\newcommand{\DELETE}[1]{}

	\newcommand{\REM}[1]{\marginpar{\bfseries\tiny{\color{blue}}}}
    \newcommand{\MARGINOTE}[1]{}
}
{
	
	\newcommand{\COMMENT}[1]{{\color{red}\uuline{#1}\color{black}}}
	\newcommand{\COMMENTGT}[1]{{\hfill\large\color{red}***{#1}***\color{black}\hfill}\\}
	\newcommand{\TODO}[1]{{\color{red}\uuline{#1}\color{black}}}
	\newcommand{\INTERNAL}[1]{\footnote{#1}}
	\newcommand{\QUESTION}[1]{{\color{brown}\uuline{#1}\color{black}}}
	\newcommand{\DELETE}[1]{{\color{red}\sout{#1}\color{black}}}

	\newcommand{\REM}[1]{\marginpar{\bfseries\tiny{\color{blue}#1}}}
\newcommand{\MARGINOTE}[1]{\marginpar{\color{red}\tiny\texttt{#1}}}
}

\newcommand\DT[1]{\mathchoice
                 {{\buildrel{\hspace*{.1em}\text{\LARGE.}}\over{#1}}}
                 {{\buildrel{\hspace*{.1em}\text{\LARGE.}}\over{#1}}}
                 {{\buildrel{\hspace*{.1em}\text{\Large.}}\over{#1}}}
                 {{\buildrel{\hspace*{.1em}\text{\large.}}\over{#1}}}}
\newcommand\pdt[1]{\frac{\partial{#1}}{\partial t}} %Partial Derivative w.r.t. t
\newcommand{\lineunder}[2]{\LU{\begin{array}[t]{c}\underbrace{#1}\vspace*{.5em}\end{array}}{\mbox{\footnotesize\rm #2}}}
\newcommand{\LU}[2]{\begin{array}[t]{c}#1\vspace*{-1em}\\_{#2}\end{array}}
\newcommand{\linesunder}[3]{\LSU{\begin{array}[t]{c}\underbrace{#1}\vspace*{.5em}\end{array}}{\mbox{\footnotesize\rm #2}}{\mbox{\footnotesize\rm #3}}}
\newcommand{\LSU}[3]{\begin{array}[t]{c}#1\vspace*{-1em}\\_{#2}\vspace*{-.5em}\\_{#3}\end{array}}
\newcommand{\morelinesunder}[4]{\LSUU{\begin{array}[t]{c}\underbrace{#1}\vspace*{.5em}\end{array}}{\mbox{\footnotesize\rm #2}}{\mbox{\footnotesize\rm #3}}{\mbox{\footnotesize\rm #4}}}
\newcommand{\LSUU}[4]{\begin{array}[t]{c}#1\vspace*{-1em}\\_{#2}\vspace*{-.5em}\\_{#3}\vspace*{-.5em}\\_{#4}\end{array}}
\newcommand{\Item}[2]{\parbox[t]{.055\textwidth}{#1}\hfill%
      \parbox[t]{.945\textwidth}{#2}\vspace*{.8mm}} 
\newcommand{\divS}{\mathrm{div}_{\scriptscriptstyle\textrm{\hspace*{-.1em}S}}^{}}
\newcommand{\nablaS}{\nabla_{\scriptscriptstyle\textrm{\hspace*{-.3em}S}}^{}}
\newcommand{\NablaS}{\Nabla_{\scriptscriptstyle\textrm{\hspace*{-.3em}S}}^{}}
\def\Vdots{\!\mbox{\setlength{\unitlength}{1em}
\begin{picture}(0,0)
\put(-.07,0){.}
\put(-.07,.3){.}
\put(-.07,.6){.}
\end{picture}\hspace*{.2em}}}
\usepackage{mathrsfs}   % loading \mathscr
\usepackage{eucal}      % changing \mathcal into EulerCaligraphic  

  \def\bbI{{\mathbb I}}

\def\FG{\boldsymbol}
   
 \def\ee{{\FG e}} \def\ff{{\FG f}} 
  
\def\jj{{\FG j}}   
\def\mm{{\FG m}} \def\nn{{\FG n}}  
  \def\rr{{\FG r}} 
   
\def\vv{{\FG v}}  \def\xx{{\FG x}} 
\def\yy{{\FG y}} \def\zz{{\FG z}} 
   
\def\DD{{\FG D}} %\def\EE{{\FG E}}
\def\FF{{\FG F}}

\def\MM{{\FG M}}   
   
\def\SS{{\FG S}} \def\TT{{\FG T}} %\def\UU{{\FG U}} 
  \def\XX{{\FG X}} 
 \def\ZZ{{\FG Z}}

\newcommand{\R}{\mathbb R}
\newcommand{\N}{\mathbb N}
\newcommand{\Nabla}{{\nabla}}
\newcommand{\Fe}{\FF}
\newcommand{\FFeps}{\FF_{\!\EPS}}
\newcommand{\FFepsk}{\FF_{\!\EPS k}}
\newcommand{\EE}{{\bm e}}
\newcommand{\pl}{\partial}
\newcommand{\eq}[1]{(\ref{#1})}
\renewcommand{\d}{\mathrm d}  
\newcommand{\barOmega}{\hspace*{.2em}{\overline{\hspace*{-.2em}\varOmega}}}

\newtheorem{theorem}{Theorem}[section]
\newtheorem{lemma}[theorem]{Lemma}
\newtheorem{definition}[theorem]{Definition}
\newtheorem{example}[theorem]{Example}

\newtheorem{remark}[theorem]{Remark}

\numberwithin{equation}{section}

\begin{document}
\begin{sloppypar}

\allowdisplaybreaks

\noindent{\LARGE\bf Thermodynamics
of viscoelastic solids,\\[.3em]its Eulerian formulation, and\\[.3em]existence of weak solutions}\ \footnote{Support from M\v SMT \v CR (Ministry of Education of the
Czech Republic) project CZ.02.1.01/0.0/0.0/15-003/0000493,
CSF grant no.\,23-06220S, and from the institutional support
RVO:61388998 (\v CR)  is acknowledged.}

\bigskip
\bigskip

\noindent{\large Tom\'a\v s Roub\'\i\v cek}

\bigskip

\noindent{Mathematical Institute, Charles University,\\
Sokolovsk\'a 83, CZ-186~75~Praha~8, Czech Republic,
\\
and\\
Institute of Thermomechanics of the Czech Academy of Sciences,\\ 
Dolej\v skova 5, CZ-182 00 Praha 8, Czech Republic}\\
email: {tomas.roubicek@mff.cuni.cz}

%\thanks

%\date{}

\begin{abstract}
The thermodynamic model of visco-elastic deformable solids
at finite strains is formulated in a fully Eulerian way in rates.
Also, effects of thermal expansion or buoyancy due to evolving mass
density in a gravity field are covered.
The Kelvin-Voigt rheology with a higher-order viscosity
(exploiting the concept of multipolar materials) is used,
allowing for physically relevant frame-indifferent stored energies and
for local invertibility of deformation. The model complies with
energy conservation and Clausius-Duhem entropy inequality. 
Existence and a certain regularity of weak solutions
are proved by a Faedo-Galerkin semi-discretization and
a suitable regularization. Subtle physical limitations of the model
are illustrated on thermally expanding neo-Hookean materials or
materials with phase transitions.

  \medskip

{\noindent{\bf Mathematics Subject Classification}.\!\!
%\subjclass{\noindent{\bf Mathematics Subject Classification}.
35Q49,\! % Transport equations
35Q74,\! % PDEs in connection with mechanics of deformable solids
35Q79,\! % PDEs in connection with classical thermodynamics and heat transfer
%35Q86, % PDEs in connection with geophysics
65M60,\! % PDE, IVP, ... Galerkin methods, 
%74A15, % Thermodynamics
74A30,\!\! % Nonsimple materials
% 74A45, % Theories of fracture and damage
%74C15, %Large-strain, rate-independent theories(including nonlinear plasticity)
74Dxx,\! %Materials of strain-rate type.
%74F10, % Fluid-solid interactions (...aero- hydro-elasticity, porosity, etc.)
74J30,\! % Nonlinear waves
%74L05, % Geophysical solid mechanics
%74R20, % Anelastic fracture and damage
%76S05. % Flows in porous media; filtration; seepage
80A20.\! % Heat and mass transfer, heat flow
%86A17. % Global dynamics, earthquake problems
}

\medskip

%\keywords
{\noindent{\bf Keywords.} 
 Elastodynamics, Kelvin-Voigt viscoelasticity,
thermal coupling, large strains, multipolar continua,
semi-Galerkin discretization, weak solutions.
}
\end{abstract}

%\maketitle

%\baselineskip=16pt

\def\TRACTION{\bm{f}}
\def\GRAVITY{\bm{g}}
\def\rhoR{\varrho_\text{\sc r}^{}}
\def\RRR{\text{\sc r}}
\def\M{m}
\def\MM{M}
\def\COUPLING{\gamma}
\def\LAM{\lambda}
\def\W{w}
\def\OMEGA{\omega}
\def\EPS{\varepsilon}

\def\ALPH{1}
\def\ONEALPH{2}
\def\DIS{\DD}
\def\TWO{2}
\def\TWOprime{2}
\def\EXP{\mu}
\def\wh{\widehat}

\section{Introduction}
%        ~~~~~~~~~~~~

Even in the isothermal cases, 
the visco-elastodynamics at finite strains in its basic simple-material
variant and the {\it Kelvin-Voigt viscoelastic rheology}
has been articulated in \cite{Ball02SOPE,Ball10PPNE} as a
difficult and essentially open problem as far as the existence of weak
solutions concerns. This holds, a-fortiori, in anisothermal situations

Let us remark that the adjective ``finite'' is often used equivalently
to ``large'' strains, although sometimes it wants to emphasize some
bounds on strains, as used also here since the deformation gradient
ranges bounded sets.

In general, there are two basic approaches: the Lagrangian and the Eulerian
ones. The former formulates the equations in a certain fixed
``reference'' configuration which, in some cases, may have a
good meaning as a configuration at which the body was manufactured
and to which the deformed body is always related in some sense.
This approach allows easily for deformation of the shape of the body
and mass density can be fixed but, on the other hand, the description of
real deforming configuration in the fixed reference
configuration needs careful pull-back and push-forward operators
to be involved in the equations to comply with frame indifference,
which was not always well respected in literature (e.g.\
\cite{CarDol04TSDN,LRXY21GWPC,KLST20EUVK}),
as pointed out in \cite{Antm98PUVS}. Also an interaction with outer spatial
fields (such as gravitation or electromagnetic) is quite complicated.
This approach is believed as proper for solids and, in the isothermal
situations, has been used e.g.\ in \cite{Demo00WSCN,Tved08QEVR} with some
restrictions to the stored energy (admitting self-interpenetration). Due to
inevitably nonlinear nonmonotone character of the problem, even the
mere existence of weak solutions is very problematic and various
concepts of generalized solutions have been therefore devised
(cf.\ also discussion of difficulties in the quasistatic
situations \cite{MiOrSe14ANVM}) or concepts of nonsimple materials with higher
gradients of deformation in the conservative part of the mechanical system
have been exploited \cite{LRXY21GWPC,MieRou20TKVR,RouSte19FTCS,RouTom20DEBL},
cf.\ also \cite[Chap.9]{KruRou19MMCM}.
For completeness, let us still report inviscid purely elastodynamic
studies under strong qualification (certain convexity) of the stored energy
\cite{Dafe86QHSI,DafHru86EMQH,HuKaMa77WPQL}
or using a certain very weak (measure-valued) concept of solutions
\cite{CarRie04YMAE,DeStTz01VAST,DeStTz12WSUD,Proh08CFEB,Rieg03YMSN}.
This applies also to anisothermal extension leading to thermodynamically
consistent models in the Lagrangian setting in the Kelvin-Voigt rheology in
\cite{ChrTza18REHP,KruRou19MMCM,MieRou20TKVR}, inviscid
in \cite{BLGG21FOHF,ChGaTz20DVSI}, and in Maxwellian rheology in \cite{RouSte19FTCS}.

The latter mentioned {\it approach},  {\it Eulerian}, is (with some
exception as e.g.\ \cite{PruTum21TFHG}) standardly believed
to be well fitted with fluids. It is particularly suitable in situations
when there is no natural reference configuration or where a reference
configuration becomes less and less relevant during long-time evolution,
which may however apply also to solids. A formulation of equations in the
current deforming configuration needs rather
velocity/strain than displacement to be involved in the momentum
equation. The advantage is an easier possibility to involve
interaction with outer spatial fields (here gravitational) and
avoiding the pull-back and push-forward manipulation. On the other
hand, there is a necessity to involve convective derivative and transport
equations and also evolving the shape of the body is troublesome.
In isothermal situations, such a model was formulated and analyzed as
incompressible in \cite{LeLiZh08GSIV,LiuWal01EDFC} and as compressible in
\cite{HuMas16GSRH,QiaZha10GWPC}. The mentioned higher gradients that
would allow for reasonable analysis are now to be involved rather in the
dissipative than the conservative part, so that their influence
manifests only in fast evolutions. Besides, the higher gradients
may lead to easier propagation of elastic waves with less dispersion
and less attenuation in some frequency ranges, compared with the
usual 1st-order gradients, as pointed out in \cite[Sect.3.1]{Roub23GTLV}.
In the isothermal situations, it was used in a quasistatic case in \cite{Roub??VELS} and
in the dynamical case in \cite{RouSte19FTCS} when considering
the stored energy in the actual configuration, which then gives
an energy pressure in the stress tensor. In anisothermal situations,
such free-energy pressure would be directly added to stress tensor in an
non-integrable way and likely would cause technical difficulties. 

In comparison with \cite{Roub??VELS,RouSte19FTCS},
the novelty of this paper is to apply the Eulerian approach
to solids in anisothermal situations, using the free energy
in a reference configuration (as in \cite{PruTum21TFHG} without
any analysis), which does not see the energy-type pressure in
the stress tensor and which is also better fitted with
usually available experimental data typically obtained in the
undeformed (reference) configuration. The analysis combines
$L^1$-theory for the heat equation adapted to the convective
time derivatives and the techniques from compressible fluid dynamics
adapted for solids. The main attributes of the devised model are:\\
\Item{{\bf---}}{Concept of {\it hyperelastic materials} (whose conservative-stress
response comes from a potential being thermodynamically a free energy)
combined with the {\it Kelvin-Voigt viscoelastic rheology}.}
\Item{{\bf---}}{The rate formulation in terms of velocity and deformation
gradient is used while the deformation itself does not explicitly occur.}
\Item{{\bf---}}{Mechanical consistency in the sense that
{\it frame indifference} of the free energy (which is in particular
{\it nonconvex} in terms of deformation gradient) is admitted, as well
as its {\it singularity} under infinite compression in relation with
{\it local non-interpenetration}.}
\Item{{\bf---}}{Thermodynamic consistency of the thermally coupled system
in the sense that the {\it total energy is conserved} in a closed system,
the {\it Clausius-Duhem entropy inequality}
holds, and temperature stays non-negative.}
\Item{{\bf---}}{The nonconservative part of the stress in the Kelvin-Voigt model
containing a higher-order component reflecting the concept of nonsimple
(here 2nd-grade) {\it multipolar media} is exploited.}
\Item{{\bf---}}{The model allows for rigorous mathematical analysis as far as
existence and certain regularity of energy-conserving weak solutions
concerns.}

\vspace*{.3em}

\noindent
As far as the non-negativity of temperature, below we will be able
to prove only that at least some solutions enjoy this attribute, although
there is an intuitive belief that all possible solutions will make it and
a hope that more advanced analytical techniques would rigorously prove it.

The main notation used in this paper is summarized in the following table:

\begin{center}
\fbox{
\begin{minipage}[t]{14.5em}\small

$\vv$ velocity (in m/s),

$\varrho$ mass density (in kg/m$^3$),

$\rhoR$ referential mass density, 

$\Fe$ deformation gradient,

$\theta$ temperature (in K)

$\TT$ Cauchy stress (symmetric, in Pa),

$\mathfrak{H}$ hyper-stress (in Pa\,m),

$\jj$ heat flux (in W/m$^2$),

$\ff$ traction load,

\vspace*{.3em}

$\det(\cdot)$ determinant of a matrix,

${\rm Cof}(\cdot)$ cofactor matrix,

$\nu_\flat>0$ a boundary viscosity

$\R_{\rm sym}^{d\times d}=\{A\in\R^{d\times d};\ A^\top=A\}$.

\end{minipage}\hspace*{1em}
%\hfill 
\begin{minipage}[t]{23em}

$\psi=\psi(\FF,\theta)$ referential free energy (in J/m$^3$=Pa)

$\varphi=\varphi(\Fe)$ referential stored energy (in J/m$^3$=Pa),

$\COUPLING=\COUPLING(\Fe,\theta)$ referential heat part of free energy,

$\ee(\vv)=\frac12\Nabla\vv^\top\!+\frac12\Nabla\vv$ small strain rate (in s$^{-1}$),

$\DIS=\DIS(\FF,\theta;\ee(\vv))$ dissipative part of Cauchy stress,

$\W$ heat part of internal energy (enthalpy, in J/m$^3$),

$(^{_{_{\bullet}}})\!\DT{^{}}=\pdt{}{^{_{_{\bullet}}}}+(\vv{\cdot}\Nabla)^{_{_{\bullet}}}$
convective time derivative,

$\cdot$ or $:$ scalar products of vectors or matrices, 

$\Vdots\ \ $ scalar products of 3rd-order tensors,

$\kappa=\kappa(\Fe,\theta)$ thermal conductivity (in W/m$^{-2}$K$^{-1}$),

$c=c(\Fe,\theta)$ heat capacity (in Pa/K),

$\nu>0$ a bulk hyper-viscosity coefficient,

$\GRAVITY$ external bulk load (gravity acceleration in m/s$^{2}$).

\smallskip \end{minipage}
}\end{center}

\vspace{-.5em}

\begin{center}
{\small\sl Table\,1.\ }
{\small
\,Summary of the basic notation used through the paper.
}
\end{center}

The paper is organized as follows: The formulation of the model
in actual Eulerian configuration and its
energetics and thermodynamics is presented in Section~\ref{sec-model}.
Then, in Section~\ref{sec-anal}, the rigorous analysis by a suitable
regularization and a (semi) Faedo-Galerkin approximation is performed,
combined with the theory of transport by regular velocity field which is
briefly presented in the Appendix~\ref{sec-app3}. Some applicable
examples are then stated in Section~\ref{sec-examples}.

\section{The thermodynamic model and its energetics}\label{sec-model}
%        ~~~~~~~~~~~~~~~~~~~~~~~~~~~~~~~~~~~~~~~~~~~~

It is important to distinguish carefully the referential and
the actual time-evolving coordinates. We aim to formulate
the model eventually in actual configurations, i.e.\ the Eulerian
formulation, reflecting also the reality in many (or even most)
situations that the concept of a reference
configuration is only an artificial construction and, even
if relevant in some situations, becomes successively more and
more irrelevant during evolution at truly finite strains and large
displacements. A very typical example is geophysical models for
long time scales. On the other hand, some experimental material data 
are related to some reference configuration -- typically, it
concerns mass density and stored or free energies per referential
volume (in J/m$^3$=Pa) as considered here or per mass (in J/kg).

\subsection{Finite-strain kinematics and mass and momentum transport}
%           ~~~~~~~~~~~~~~~~~~~~~~~~~~~~~~~~~~~~~~~~~~~~~~~~~~~~~~~~~~
We briefly remind the fundamental concepts and formulas which can
mostly be found in the monographs, as e.g.\ \cite[Part~XI]{GuFrAn10MTC},
\cite[Sect.~7.2]{Mart19PCM}, \cite[Sect.~22.1]{Silh97MTCM}, or \cite{True69RT}.

In finite-strain continuum mechanics, the basic geometrical concept is the
time-evolving deformation $\yy:\varOmega\to\R^d$ as a mapping from a reference
configuration of the body $\varOmega\subset\R^d$ into a physical space $\R^d$.
The ``Lagrangian'' space variable in the reference configuration will be
denoted as $\XX\in\varOmega$ while in the ``Eulerian'' physical-space
variable by $\xx\in\R^d$. The basic geometrical object is the deformation
gradient $\FF=\Nabla_{\!\XX}^{}\yy$.

We will be interested in deformations $\xx=\yy(t,\XX)$ evolving in time,
which are sometimes called ``motions''. The difference $\xx{-}\XX$
represents the {\it displacement}. The important quantity is the
Eulerian velocity $\vv=\DT\yy=\pdt{}\yy+(\vv{\cdot}\Nabla)\yy$.
Here and throughout the whole article, we use the dot-notation
$(\cdot)\!\DT{^{}}$ for the {\it convective time derivative}
applied to scalars or, component-wise, to vectors or tensors.

Then, the velocity gradient
$\Nabla\vv=\nabla_{\!\XX}^{}\vv\nabla_{\!\xx}^{}\XX=\DT\FF\FF^{-1}$,
where we used the chain-rule calculus  and
$\FF^{-1}=(\nabla_{\!\XX}^{}\xx)^{-1}=\nabla_{\!\xx}^{}\XX$. 
This gives the {\it transport equation-and-evolution  for the
deformation gradient} as
\begin{align}%\nonumber\\[-2.7em]
\DT\FF=(\nabla\vv)\FF\,.
  \label{ultimate}\end{align}
From this, we also obtain the evolution-and-transport equation for the
determinant $\det\FF$ as
\begin{align}%\nonumber
\label{DT-det}\DT{\overline{\det\FF}}&=
{\rm Cof}\FF{:}\DT\FF=(\det\FF)\FF^{-\top}{:}\DT\FF
%\\[-.3em]&
=(\det\FF)
\bbI{:}\DT\FF\FF^{-1}=(\det\FF)\bbI{:}\Nabla\vv=(\det\FF){\rm div}\,\vv\,,
\end{align}
where $\bbI$ denotes the identity matrix, as well as the
evolution-and-transport equation for $1/\!\det\FF$ as
\begin{align}%\nonumber\\[-2.9em]
\DT{\overline{\!\!\!\bigg(\frac1{\det\FF}\bigg)\!\!\!}}\ 
=-\frac{{\rm div}\,\vv}{\det\FF}\,.
\label{DT-det-1}\end{align}

The understanding of \eqref{ultimate} is a bit delicate because it mixes
the Eulerian $\xx$ and the Lagrangian $\XX$; note that
$\nabla\vv=\nabla_\xx\vv(\xx)$ while standardly
$\FF=\nabla_{\!\XX}^{}\yy=\FF(\XX)$. In fact, we consider
$\FF{\circ}\bm\xi$ where $\bm\xi:\xx\mapsto\yy^{-1}(t,\XX)$
is the so-called {\it return} (sometimes called also a {\it reference})
{\it mapping}. Thus, $\FF$ depends on $\xx$ and \eqref{ultimate} and
is an equality which holds for a.a.\ $\xx$.
The same holds for \eqref{DT-det} and \eqref{DT-det-1}. The reference
mapping $\bm\xi$, which is well defined through its transport equation
\begin{align}\nonumber\\[-2.2em]
\DT{\bm\xi}=\bm0\,,
\label{transport-xi}\end{align}
actually does not explicitly occur in the formulation of the problem
(except Remark~\ref{rem-inhomogenous}). Here, we will benefit from the
boundary condition $\vv{\cdot}\nn=0$ below, which causes that the actual
domain $\varOmega$ does not evolve in time. The same convention concerns
temperature $\theta$ and thus also $\TT$, $\eta$, and $\DIS$ in
\eq{stress-entropy} and \eq{Cauchy-dissip} below, which will make the
problem indeed fully Eulerian. Cf.\ the continuum-mechanics textbooks as e.g.\
\cite{GuFrAn10MTC,Mart19PCM}.

Most physical variables can be classified either as {\it extensive
variables} (as mass density, or free or internal energies, or entropy)
or {\it intensive variables} (as temperature or velocity). These two
are transported differently in compressible media. 

The mass density (in kg/m$^3$) is an extensive variable and its transport
(expressing that the conservation of mass) writes as the {\it continuity
equation} $\pdt{}\varrho+{\rm div}(\varrho\vv)=0$,
or, equivalently, the {\it mass transport equation}
\begin{align}\nonumber\\[-2.7em]
\DT\varrho=-\varrho\,{\rm div}\,\vv\,.
\label{cont-eq+}\end{align}
This equation also ensures the useful equality, namely 
transport of the momentum $\varrho\vv$ as another extensive
(vector-valued) variable:
\begin{align}\label{inertial}
\frac{\partial}{\partial t} (\varrho \bm{v}) +
\text{\rm div}(\varrho \bm{v}{\otimes} \bm{v})
=\varrho\pdt\vv+\pdt\varrho\vv+{\rm div}(\varrho\vv)\vv
+\varrho(\vv{\cdot}\nabla)\vv=\varrho\DT\vv\,.
\end{align}
One can determine the density $\varrho$ instead of 
the transport equation for mass density \eq{cont-eq+}
from the algebraic relation
\begin{align}\nonumber\\[-2.7em]
\varrho=\frac{\rhoR}{\det\FF}\,,
\label{density-algebraically}\end{align}
where $\rhoR$ is the mass density in the reference
configuration. Indeed, relying on \eq{DT-det}, one has the calculus
\begin{align}%\nonumber\\[-3.2em]
\frac{\DT\varrho}{\varrho}=
\Bigg(\rhoR\,\,\DT{\overline{\!\!\bigg(\frac1{\det\FF}\bigg)\!\!}}\,\,
+\frac{\DT\rhoR}{\det\FF}\Bigg)
\frac{\det\FF}{\rhoR}=-\frac{\DT{\overline{\det\FF}}}{\det\FF}
=-{\rm div}\,\vv
\label{towards-cont-eq}\end{align}
because $\DT\rhoR=0$. Later, we will consider the initial
conditions $\FF_0$ for \eq{ultimate} and $\varrho_0$ for \eq{cont-eq+}. 
These two should be related by $\varrho(0)=\rhoR/\!\det\FF_0$.

Alternatively to \eq{cont-eq+}, we can also consider an
evolution-and-transport equation for the ``mass sparsity'' as the inverse
mass density $1/\varrho$:
\begin{align}\nonumber\\[-2.7em]
\DT{\overline{\!\!\bigg(\frac1\varrho\bigg)\!\!}}\
=\frac{{\rm div}\,\vv}{\varrho}\,.
\label{cont-eq-inverse}\end{align}

\subsection{Hyperelasticity and visco-elasticity, and its thermodynamics}
%           ~~~~~~~~~~~~~~~~~~~~~~~~~~~~~~~~~~~~~~~~~~~~~~~~~~~~~~~~~~~~~

The main ingredients of the model are the (volumetric) {\it free energy}
$\psi$ depending on deformation gradient $\FF$ and temperature $\theta$,
 and the temperature-dependent {\it dissipative stress} $\DIS$
not necessarily possessing any underlying potential. The free energy
$\psi=\psi(\FF,\theta)$ is considered per the {\it referential volume}, while
the free energy per actual deformed volume is $\psi(\FF,\theta)/\!\det\FF$.
Considering the free energy by reference volume is more standard in
continuum physics \cite{GuFrAn10MTC,Mart19PCM} than the free energy per
actual evolving volume. This corresponds to experimentally available data
and seems particularly more suitable for analysis of thermally coupled
systems. This last benefit is related to the fact that the
referential free energy does not give an energy pressure contribution
to the Cauchy stress (cf.\ the last term in \eq{referential-stress} below
or\ \cite[Rem.\,2]{Roub??VELS}) and allows for more easy decoupling estimation
strategy decoupling the mechanical part and the thermal part of the coupled
system.

While $\psi$ is referential, the free energy per actual volume is then
$\psi(\Fe,\theta)/\!\det\Fe$. From it, we can read the {\it conservative
part of the actual Cauchy stress} $\TT$ and the {\it actual entropy}
$\eta$ as:
\begin{align}
\TT=\frac{\psi_\Fe'(\Fe,\theta)\Fe^\top\!\!}{\det\Fe}\ \ \ \ \text{ and }\ \ \ \ 
\eta=-\,\frac{\!\psi_\theta'(\Fe,\theta)}{\det\Fe}\,.
\label{stress-entropy}\end{align}

In the already anticipated viscoelastic {\it Kelvin-Voigt rheological model},
we will use also a dissipative contribution to the Cauchy stress,
which will make the system parabolic. 
In addition to the usual first-order stress $\DD$, we consider a dissipative
contribution to the Cauchy stress involving also a higher-order 2nd-grade
{\it hyper-stress} $\mathcal{H}$, so that the overall dissipative stress is:
\begin{align}\nonumber
\DD-{\rm div}\,\mathcal{H}
\ \ \ &\text{ with }\ \DD=\DD(\FF,\theta;\ee(\vv)) \ \text{ and }\ 
\mathcal{H}=\mathcal{H}(\Nabla\ee(\vv))
\\
&\text{ for some }
\DD(\FF,\theta;\cdot):\R_{\rm sym}^{d\times d}\to\R_{\rm sym}^{d\times d}\ \text{ and }\
\mathcal{H}({\bm E})=\nu|{\bm E}|^{p-2}{\bm E}\,.
\label{Cauchy-dissip}\end{align}
More discussion about the hyper-stress contribution will be in
Remarks~\ref{rem-grad} below.

The {\it momentum equilibrium} equation then balances the divergence of
the total Cauchy stress with the inertial and gravity force:
\begin{align}
\varrho\DT\vv-{\rm div}\big(\TT{+}\DD
{-}{\rm div}\,\mathcal{H})
=\varrho\GRAVITY
\label{Euler-thermodynam1-}\end{align}
with $\TT$ from \eq{stress-entropy} and $\DD$ and $\mathcal{H}$ from
\eq{Cauchy-dissip}. 

The second ingredient in \eq{stress-entropy} is subjected to 
the {\it entropy equation}:
\begin{align}
  \pdt\eta+{\rm div}\big(\vv\,\eta\big)
  =\frac{\xi-{\rm div}\,{\bm j}}\theta\ \ \ \ \text{ with }\ 
  \jj=-\kappa(\Fe,\theta)\nabla\theta \,
\label{entropy-eq}\end{align}
with $\xi$ the heat production rate and ${\bm j}$ the heat flux.
The latter equality in \eq{entropy-eq} is the {\it Fourier law}
determining phenomenologically the heat flux $\jj$ proportional to the
negative gradient of temperature $\theta$ through the thermal conductivity 
coefficient $\kappa=\kappa(\Fe,\theta)$. Assuming $\xi\ge0$ and
$\kappa\ge0$ and integrating \eq{entropy-eq} over the domain $\varOmega$
while imposing the non-penetrability of the boundary in the sense that
the normal velocity $\vv{\cdot}\nn$ vanishes across the boundary
$\varGamma$ of $\varOmega$, we obtain the {\it Clausius-Duhem inequality}:
\begin{align}
\frac{\d}{\d t}\int_\varOmega\eta\,\d\xx
=\int_\varOmega\!\!\!\!\!\!\!\!\!\lineunder{\frac\xi\theta+\kappa\frac{|\nabla\theta|^2}{\theta^2}}{entropy production rate}\!\!\!\!\!\!\!\!\!\!\!\d\xx
+\int_\varGamma\!\!\!\!\lineunder{\Big(\kappa\frac{\nabla\theta}\theta-\eta\vv\Big)}{entropy flux}\!\!\!\!\!\!\cdot\nn\,\d S
\ge\int_\varGamma\!\kappa\frac{\nabla\theta{\cdot}\nn}\theta\,\d S\,.
\label{entropy-ineq}\end{align}
If the system is thermally isolated in the sense that the normal heat flux
$\jj{\cdot}\nn$ vanishes across the boundary $\varGamma$, we recover the
{\it 2nd law of thermodynamics}, i.e.\ the total entropy in isolated
systems is nondecreasing in time.

Substituting $\eta$ from \eq{stress-entropy} into \eq{entropy-eq} written
in the form $\theta\DT\eta=\xi-{\rm div}\,\jj-\theta\eta{\rm div}\,\vv$
and using the calculus $(\psi_{\theta}'(\FF,\theta)/\!\det\FF)_{\!\FF}'{:}\DT\FF
+\psi_\theta'(\FF,\theta)\,({\rm div}\,\vv)/\!\det\FF
=\psi_{\FF\theta}''(\FF,\theta){:}\DT\FF/\!\det\FF$, 
we obtain the {\it heat-transfer equation} 
\begin{align}\nonumber
&c(\FF,\theta)\DT\theta
=\DIS(\FF,\theta;\EE(\vv)){:}\EE(\vv)+\nu|\nabla\EE(\vv)|^p\!
+\theta\frac{\psi_{\FF\theta}''(\FF,\theta)}{\det\FF}{:}\DT\FF
-{\rm div}\,{\bm j}
\\[-.6em]&\hspace{16em}\text{ with the heat capacity }\
c(\FF,\theta)=-\theta\,\frac{\!\psi_{\theta\theta}''(\FF,\theta)}{\det\FF}\,.
\label{heat-eq+}\end{align}

The referential {\it internal energy} is given by the {\it Gibbs relation}
$\psi-\theta\psi_\theta'$.
In terms of the {\it thermal part of the actual internal energy}
$\W(\FF,\theta):=(\psi(\FF,\theta){-}\theta\psi_\theta'(\FF,\theta)
{-}\psi(\FF,0))/\det\FF$, the {\it heat equation} can be written in the
so-called {\it enthalpy formulation}: 
\begin{align}\nonumber
\pdt\W+{\rm div}\big(\vv\W{+}{\bm j}\big)=
\DIS(\FF,\theta;\EE(\vv)){:}\EE(\vv)+\nu|\nabla\EE(\vv)|^p
+\frac{\psi'_{\Fe}(\Fe,\theta){-}\psi'_{\Fe}(\Fe,0)}{\det\Fe}{:}\DT\Fe\
\\\ \text{ with }\ \W=
\frac{\psi(\Fe,\theta)-\theta\psi_\theta'(\Fe,\theta)-\psi(\FF,0)}{\det\Fe}\,.
\label{Euler-thermodynam3-}\end{align}
Here, alternatively, we could use also $\DT\Fe=(\nabla\vv)\Fe$ for writing 
$\psi'_{\Fe}(\Fe,\theta){:}\DT\Fe
=\psi'_{\Fe}(\Fe,\theta)\Fe^\top\!{:}(\nabla\vv)$
$=\psi'_{\Fe}(\Fe,\theta)\Fe^\top\!{:}\ee(\vv)$; here the assumed
frame indifference of $\psi(\cdot,\theta)$ assumed in
\eq{frame-indifference} below, leading to symmetry of
$\psi'_{\Fe}(\Fe,\theta)\Fe^\top$, was employed, cf.\ \eq{Euler-thermodynam3}
below. In comparison with \eq{heat-eq+} which needs differentiability of
$\psi\theta'$, the form \eq{Euler-thermodynam3-} needs 
differentiability of $\psi$ only.
 Using again the algebra $F^{-1}={\rm Cof}\,F^\top\!/\!\det F$ and the
 calculus $\det'(F)={\rm Cof}\,F$ and \eq{ultimate}, we realize that
\begin{align}\nonumber
 & \pdt\W+{\rm div}\big(\vv\,\W\big)=
 \DT\W+\W\,{\rm div}\,\vv
 =\ \DT{\overline{\!\!\Big(\frac{\COUPLING(\FF,\theta)-\theta\COUPLING_\theta'(\FF,\theta)}{\det\FF}\Big)\!\!}}
\,+\frac{\COUPLING(\FF,\theta)-\theta\COUPLING_\theta'(\FF,\theta)\!}{\det\FF}\,{\rm div}\,\vv
\\[.2em]&\nonumber\ =\Big(\Big[\frac\COUPLING\det\Big]_\FF'(\Fe,\theta)
-\theta\Big[\frac{\COUPLING_\theta'}\det\Big]_\FF'(\Fe,\theta)\Big){:}\DT\FF
-\theta\frac{\COUPLING_{\theta\theta}''(\Fe,\theta)}{\det\FF}\DT\theta
 +\frac{\COUPLING(\FF,\theta)-\theta\COUPLING_\theta'(\FF,\theta)\!}{\det\FF}\,{\rm div}\,\vv
\\&\nonumber\
=\Big(\frac{\COUPLING_{\!\FF}'(\FF,\theta)}{\det\FF}
-\theta\frac{\COUPLING_{\FF\theta}''(\FF,\theta)}{\det\FF}\Big){:}\DT\FF
+c(\Fe,\theta)\DT\theta
%\\[-.4em]&\nonumber\hspace*{6.2em}
+\hspace{-1em}\lineunder{\big(\COUPLING(\FF,\theta)-\theta\COUPLING_\theta'(\FF,\theta)\big)\Big(\frac1{\det\FF}\Big)'{:}\DT\FF
+\frac{\COUPLING(\Fe,\theta){-}\theta\COUPLING_\theta'(\Fe,\theta)\!}{\det\Fe}\,{\rm div}\,\vv}{$\ \ =0$}
\\[-1em]&\
=c(\Fe,\theta)\DT\theta+
\frac{\COUPLING_{\!\FF}'(\FF,\theta)}{\det\FF}{:}\DT\FF
-\frac{\theta\COUPLING_{\FF\theta}''(\FF,\theta)}{\det\FF}{:}\DT\FF\,,
 \end{align}
 where we abbreviated $\COUPLING(\FF,\theta)=\psi(\FF,\theta)-\psi(\FF,0)$.
Thus we can see that \eq{Euler-thermodynam3-} is indeed equivalent
with \eq{heat-eq+}.

Let us formulate the thermo-visco-elastodynamic system for
$(\varrho,\vv,\FF,\theta)$, composing the equations \eq{ultimate},
\eq{cont-eq+}, \eq{Euler-thermodynam1-}, and \eq{Euler-thermodynam3-}:
\begin{subequations}\label{Euler-thermodynam}
\begin{align}\label{Euler-thermodynam0}
&\pdt\varrho=-\,{\rm div}(\varrho\vv)\,,
\\\nonumber
&\pdt{}(\varrho\vv)={\rm div}\Big(\TT(\Fe,\theta){+}
\DIS(\FF,\theta;\EE(\vv))
{-}{\rm div}\mathcal{H}(\nabla\EE(\vv))
-\varrho\vv{\otimes}\vv\Big)
  +\varrho\GRAVITY\,
 \\[-.3em]
    &\hspace*{5em}
\ \text{ with }\ \TT(\Fe,\theta)=\frac{\ \psi_{\Fe}'(\Fe,\theta)\Fe^\top\!\!}{\det\Fe}
\ \ \ \text{and }\ \,\mathcal{H}(\nabla\EE(\vv))=\nu|\nabla\EE(\vv)|^{p-2}\nabla\EE(\vv)\,,
    \label{Euler-thermodynam1}\\[-.4em]
&\pdt\Fe=(\Nabla\vv)\Fe-(\vv{\cdot}\nabla)\Fe\,,
\label{Euler-thermodynam2}
\\&\nonumber
\pdt{\W}
=\DIS(\FF,\theta;\EE(\vv)){:}\EE(\vv)+\nu|\nabla\EE(\vv)|^p
+\frac{\COUPLING'_{\Fe}(\Fe,\theta)\Fe^\top\!\!}{\det\Fe}
{:}\ee(\vv)+{\rm div}\big(\kappa(\FF,\theta)\nabla\theta-\W\vv\big)
\\&\hspace{5em}
\text{ with }\
\W=\OMEGA(\Fe,\theta){:=}\frac{\COUPLING(\Fe,\theta){-}\theta\COUPLING_\theta'(\Fe,\theta)\!}
{\det\Fe}
\ \text{ where }\ \COUPLING(\FF,\theta)=\psi(\FF,\theta){-}\psi(\FF,0).
\label{Euler-thermodynam3}
\end{align}\end{subequations}

Denoting by $\nn$ the unit outward normal to the (fixed) boundary $\varGamma$
of the domain $\varOmega$, we complete this system by suitable boundary
conditions: 
\begin{subequations}\label{Euler-thermodynam-BC}
\begin{align}\label{Euler-thermodynam-BC-1}
&\vv{\cdot}\nn=0,\ \ \ \ 
\big[(\TT{+}\DD{-}{\rm div}\mathcal{H})\nn{-}\divS\big(\mathcal{H}
\nn\big)\big]_\text{\sc t}^{}\!+\nu_\flat\vv=\ff\,,\ \ 
\\&\label{Euler-thermodynam-BC-2}
\Nabla\ee(\vv){:}(\nn{\otimes}\nn)={\bm0}\,,\ \ \ \text{ and }\ \ \ 
\kappa(\FF,\theta)\nabla\theta{\cdot}\nn=h(\theta)+\frac{\nu_\flat}2|\vv|^2
\end{align}\end{subequations}
with $\nu_\flat>0$ a boundary viscosity coefficient and with
$[\,\cdot\,]_\text{\sc t}^{}$ a tangential part of a vector. Here $\divS={\rm tr}(\nablaS)$ denotes the $(d{-}1)$-dimensional
surface divergence with ${\rm tr}(\cdot)$ being the trace of a
$(d{-}1){\times}(d{-}1)$-matrix and
$\nablaS v=\nabla v-\frac{\partial v}{\partial\nn}\nn$ 
being the surface gradient of $v$. 
The first condition (i.e.\ normal velocity zero) expresses 
nonpenetrability of the boundary was used already for \eq{entropy-ineq}
and is most frequently adopted in literature for Eulerian formulation.
This simplifying assumption fixes the shape of $\varOmega$ in its
referential configuration allows also for considering the fixed boundary even
for such time-evolving Eulerian description. The latter condition in
\eq{Euler-thermodynam-BC-1} involving a boundary viscosity comes from
the Navier boundary condition largely used in fluid dynamics and is here
connected with the technique used below, which is based on the total energy
balance as the departing point and which, unfortunately, does not allow to cope
with $\nu_\flat=0$ and simultaneously $\ff\ne0$. This boundary viscosity
naturally may contribute to the heat production on the boundary as well as
to the outflow of the heat energy to the outer space. For notational
simplicity, we consider that it is just equally distributed, one part
remaining on the boundary of $\varOmega$ and the other part 
leaving outside, which is related to the coefficient $1/2$ in 
the latter condition in \eq{Euler-thermodynam-BC-2}.

\begin{remark}[{\sl Gradient theories in rates}]\label{rem-grad}\upshape
So-called gradient theories in continuum-mechanical viscoelastic models
are nowadays very standard, referred as {\it nonsimple materials}, determining
some internal length scales, allowing for modelling various dispersion
and attenuation of propagation of elastic waves, and often facilitating
mathematical analysis, cf.\ \cite{Roub23GTLV}. They can be applied to the
conservative stress through the free energy or to the dissipation stress.
Here, in contrast to the Lagrangian solid-mechanical models as in
\cite{KruRou19MMCM,MieRou20TKVR,RouSte19FTCS,RouTom20DEBL},
we have used the latter option in \eq{Cauchy-dissip} which is better
fitted to the rate formulation and which can make the velocity field enough
regular, as vitally needed for the transport of $\varrho$ and $\FF$ in
the Eulerian models. The higher gradient hyper-stress as used below in
\eq{Cauchy-dissip}
follows the theory by E.~Fried and M.~Gurtin \cite{FriGur06TBBC}, as already
anticipated in the general nonlinear context of {\it multipolar fluids} by
J.~Ne\v cas at al.\ \cite{Neca94TMF,NeNoSi91GSCI,NecRuz92GSIV}
or solids \cite{Ruzi92MPTM,Silh92MVMS}, inspired by 
R.A.\,Toupin \cite{Toup62EMCS} and R.D.\,Mindlin \cite{Mind64MSLE}.
Let us emphasize that, without such higher-order gradients, the (possibly)
irregular velocity fields may make the treatment of the transport problem
extremely nontrivial due to the possible onset of singularities, whose
occurrence in solids may be debatable, cf.\ \cite{AlCrMa19LRCE}.
\end{remark}

\begin{remark}[{\sl Spatially inhomogeneous media}]\upshape
\label{rem-inhomogenous}
Since the free energy $\psi$ is considered in the referential domain,
we can naturally generalize the model for an initially inhomogeneous
medium making $\psi$ and also $\DIS$ $\XX$-dependent, i.e.\
$\psi=\psi(\XX,\FF,\theta)$ and $\DIS=\DIS(\XX,\FF,\theta;\EE)$. In the
Eulerian formulation, we should complete the system \eq{Euler-thermodynam} by
the transport equation \eq{transport-xi} for the return mapping $\bm{\xi}$
and compose $\psi{\circ}\bm{\xi}$ defined as
$[\psi{\circ}\bm{\xi}](\xx,\FF,\theta)=\psi(\bm{\xi}(\xx),\FF,\theta)$
and analogously for $\DIS$. It yields also inhomogeneous heat capacity.
Analogous generalization concerns also an inhomogeneous heat transfer
coefficient $\kappa=\kappa(\XX,\FF,\theta)$ or the external heat
flux $h=h(\XX,\theta)$. Another generalization may concern an anisotropic
materials with $\kappa$ being $\R_{\rm sym}^{d\times d}$-valued.
\end{remark}

\subsection{Energetics behind the system \eq{Euler-thermodynam}--\eq{Euler-thermodynam-BC}}\label{sec-app2}
%        ~~~~~~~~~~~~~~~~~~~~~~~~~~~~~~~~~~~~~~~~~~~~~~~~

The mechanical energy-dissipation balance of the visco-elastodynamic
model (\ref{Euler-thermodynam}a--c) can be seen when testing the momentum
equation \eq{Euler-thermodynam1} by $\vv$ while using the continuity equation
\eq{Euler-thermodynam0} and the evolution-and-transport equation
\eq{Euler-thermodynam2} for $\FF$.
We will select out the temperature-independent stored energy $\varphi$
by denoting the mere stored energy $\varphi:=\psi(\cdot,0)$ and
the resting temperature-dependent part $\COUPLING(\cdot,\theta):=
\psi(\cdot,\theta)-\varphi(\cdot)$ as used in \eq{Euler-thermodynam3}.
Thus, we have the split:
\begin{align}
\psi(\Fe,\theta)=\varphi(\Fe)+\COUPLING(\Fe,\theta)\ \ \ \text{ with }\ \ 
\COUPLING(\Fe,0)=0\,.
\label{ansatz}\end{align}
As already said, $\psi$ together with $\varphi$ and $\COUPLING$ are
considered per the referential volume.

Using the algebra $F^{-1}={\rm Cof}\,F^\top\!/\!\det F$ and the calculus
$\det'(F)={\rm Cof}\,F$, we can write the conservative part of the Cauchy
stress as
\begin{align}\nonumber
\frac{\varphi'(\FF)}{\det\FF}\FF^\top
&=\ \frac{\varphi'(\FF)-\varphi(\FF)\FF^{-\top}\!\!\!\!}{\det\FF}\FF^\top\!+
\frac{\varphi(\FF)}{\det\FF}\bbI
\\&
=\bigg(\frac{\varphi'(\FF)}{\det\FF}
-\frac{\varphi(\FF){\rm Cof}\FF}{(\det\FF)^2}\bigg)\FF^\top+
\frac{\varphi(\FF)}{\det\FF}\bbI
=\Big[\frac{\varphi(\FF)}{\det\FF}\Big]'\FF^\top\!+\frac{\varphi(\FF)}{\det\FF}\bbI\,.
\label{referential-stress}\end{align}
Let us recall that $[\varphi/\!\det](\FF)$ in \eq{referential-stress}
is the stored energy per actual (not referential)
volume. Using the calculus \eq{referential-stress}, we obtain
\begin{align}\nonumber
\int_\varOmega{\rm div}\,\TT{\cdot}\vv\,\d\xx
&=\!\int_\varGamma(\TT\nn){\cdot}\vv\,\d S-\!\int_\varOmega\!\TT{:}\ee(\vv)
\\&\nonumber=\!\int_\varGamma(\TT\nn){\cdot}\vv\,\d S
-\!\int_\varOmega\!\Big(\frac{\varphi'(\FF)\!}{\det\FF}
+\frac{\COUPLING_\FF'(\FF,\theta)}{\det\FF}\Big)\FF^\top{:}\ee(\vv)\,\d\xx
\\&\nonumber
=\!\int_\varGamma(\TT\nn){\cdot}\vv\,\d S
-\!\int_\varOmega\!\Big(\Big[\frac{\varphi(\FF)}{\det\FF}\Big]'\FF^\top
+\frac{\varphi(\FF)}{\det\FF}\bbI
+\frac{\COUPLING_\FF'(\FF,\theta)}{\det\FF}\FF^\top\Big){:}\ee(\vv)\,\d\xx
\\&\nonumber
=\!\int_\varGamma(\TT\nn){\cdot}\vv\,\d S
-\!\int_\varOmega\!\Big[\frac{\varphi(\FF)}{\det\FF}\Big]'
{:}(\Nabla\vv)\FF
+\frac{\varphi(\FF)}{\det\FF\!}\,{\rm div}\,\vv
+\frac{\COUPLING_\FF'(\FF,\theta)\FF^\top\!\!\!}{\det\FF}{:}\ee(\vv)\,\d\xx
\\[-.0em]&=\!\int_\varGamma(\TT\nn){\cdot}\vv\,\d S
-\frac{\d}{\d t}\int_\varOmega\frac{\varphi(\FF)}{\det\FF}\,\d\xx-
\int_\varOmega\!\frac{\COUPLING_\FF'(\FF,\theta)\FF^\top\!\!\!}{\det\FF}{:}\ee(\vv)
\,\d\xx\,.
\label{Euler-large-thermo}\end{align}
Here, we used the matrix algebra
$A{:}(BC)=(B^\top A){:}C=(AC^\top){:}B$ for any square matrices $A$, $B$, and $C$
and also we used \eq{Euler-thermodynam2} together with the Green formula
and the nonpenetrability boundary condition for
\begin{align}\nonumber
\int_\varOmega\!\Big[\frac{\varphi(\FF)}{\det\FF}\Big]'
{:}(\Nabla\vv)\FF
+\frac{\varphi(\FF)}{\det\FF\!}\,{\rm div}\,\vv\,\d\xx
&=\int_\varOmega\!\Big[\frac{\varphi(\FF)}{\det\FF}\Big]'
{:}\Big(\pdt\FF+(\vv{\cdot}\nabla)\FF\Big)+\frac{\varphi(\FF)}{\det\FF\!}\,{\rm div}\,\vv\,\d\xx
\\&\nonumber
=\frac{\d}{\d t}\int_\varOmega\frac{\varphi(\FF)}{\det\FF}\,\d\xx
+\!\int_\varOmega\!\nabla\Big(\frac{\varphi(\FF)}{\det\FF}\Big){\cdot}\vv+\frac{\varphi(\FF)}{\det\FF\!}\,{\rm div}\,\vv\,\d\xx
\\&
=\frac{\d}{\d t}\int_\varOmega\frac{\varphi(\FF)}{\det\FF}\,\d\xx
+\!\int_\varGamma\frac{\varphi(\FF)}{\det\FF}(\hspace*{-.7em}\lineunder{\vv{\cdot}\nn}{$=0$}\hspace*{-.7em})\,\d S\,.
\end{align}

The further contribution from the dissipative part of the Cauchy stress
uses Green's formula over $\varOmega$ twice and the surface Green formula
over $\varGamma$. We abbreviate the hyper-stress
$\mathcal{H}(\nabla\ee(\vv))=\nu|\nabla\EE(\vv)|^{p-2}\nabla\EE(\vv)$. Then
\begin{align}\nonumber
&\int_\varOmega{\rm div}
\big(\DIS(\FF,\theta;\ee(\vv))-{\rm div}\mathcal{H}(\Nabla\ee(\vv))\big){\cdot}\vv\,\d\xx
\\&\nonumber=\int_\varGamma\vv{\cdot}
\big(\DIS(\FF,\theta;\ee(\vv)){-}{\rm div}\mathcal{H}(\Nabla\ee(\vv))\big)\nn\,\d S
-\int_\varOmega\!\big(\DIS(\FF,\theta;\ee(\vv))-{\rm div}\mathcal{H}(\Nabla\ee(\vv))\big){:}\Nabla\vv\,\d\xx
\\&\nonumber=\int_\varGamma\vv{\cdot}
\big(\DIS(\FF,\theta;\ee(\vv)){-}{\rm div}\mathcal{H}(\Nabla\ee(\vv))\big)\nn
-\nn{\cdot}\mathcal{H}(\Nabla\ee(\vv)){:}\Nabla\vv\,\d S
\\[-.4em]&\nonumber\hspace{15em}
-\int_\varOmega\!
\DIS(\FF,\theta;\ee(\vv)){:}\ee(\vv)+\mathcal{H}(\Nabla\ee(\vv))\Vdots\Nabla^2\vv\,\d\xx
\\&\nonumber=\int_\varGamma\!
\mathcal{H}(\Nabla\ee(\vv)){\Vdots}(\pl_\nn\vv{\otimes}\nn{\otimes}\nn)
+\nn{\cdot}\mathcal{H}(\Nabla\ee(\vv)){:}\NablaS\vv+\vv{\cdot}
\big(\DIS(\FF,\theta;\ee(\vv)){-}{\rm div}\mathcal{H}(\Nabla\ee(\vv))\big)\nn\,\d S
\\[-.4em]&\nonumber\hspace{15em}
-\int_\varOmega\DIS(\FF,\theta;\ee(\vv)){:}\ee(\vv)
+\mathcal{H}(\Nabla\ee(\vv))\Vdots\Nabla^2\vv\,\d\xx
\\&\nonumber=\int_\varGamma\bigg(\mathcal{H}(\Nabla\ee(\vv))
{\Vdots}(\pl_\nn\vv{\otimes}\nn{\otimes}\nn)-
\Big(
\big(\DIS(\FF,\theta;\ee(\vv)){-}{\rm div}\mathcal{H}(\Nabla\ee(\vv))\big)\nn
\\[-.4em]&\hspace{8em}
+\divS(\nn{\cdot}\mathcal{H}(\Nabla\ee(\vv)))\Big)
{\cdot}\vv\bigg)\,\d S
-\int_\varOmega\DIS(\FF,\theta;\ee(\vv)){:}\ee(\vv)
+\nu|\nabla \ee(\vv)|^p\,\d\xx\,,
\label{Euler-test-momentum++}\end{align}
where we used the decomposition of $\Nabla\vv$ into its normal part
$\pl_\nn\vv$ and the tangential part, i.e.\ written componentwise
$\nabla\vv_i=(\nn{\cdot}\nabla\vv_i)\nn+\nablaS\vv_i$.

Furthermore, the inertial force \eq{inertial} in \eq{Euler-thermodynam1}
tested by $\vv$ gives the rate of kinetic energy $\varrho|\vv|^2/2$
when using again the continuity equation \eq{cont-eq+} tested by $|\vv|^2/2$
for the identity
\begin{align}
\pdt{}\Big(\frac\varrho2|\vv|^2\Big)=\varrho\vv{\cdot}\pdt\vv+\pdt\varrho\frac{|\vv|^2}2&
=\varrho\vv{\cdot}\pdt\vv-{\rm div}(\varrho\vv)\frac{|\vv|^2}2
\,.\label{rate-of-kinetic}\end{align}
Integrating over $\varOmega$ and using the Green formula with the boundary
condition $\vv{\cdot}\nn=0$ and relying on \eq{inertial}, we obtain
\begin{align}\nonumber
&\int_\varOmega\Big(\pdt{}(\varrho\vv)+{\rm div}(\varrho\vv{\otimes}\vv)\Big){\cdot}\vv\,\d\xx=\int_\varOmega\varrho\DT\vv{\cdot}\vv\,\d\xx
\\&\qquad\qquad=\int_\varOmega
\varrho\vv{\cdot}\pdt\vv
+\varrho\vv{\cdot}(\vv{\cdot}\nabla)\vv\,\d\xx
=\frac{\d}{\d t}\int_\varOmega\frac\varrho2|\vv|^2\,\d\xx
+\!\int_\Gamma\varrho|\vv|^2\hspace{-.7em}\lineunder{\vv{\cdot}\nn}{$=0$}\hspace{-.7em}\d S\,.
\label{calculus-convective-in-F}
\end{align}

Merging the boundary integrals in \eq{Euler-large-thermo}
and in \eq{Euler-test-momentum++} and using the boundary condition
$((\TT{+}\DD{-}{\rm div}\mathcal{H})
\nn-\divS(\mathcal{H}\nn))_\text{\sc t}^{}+\nu_\flat\vv=\ff$, we thus obtain
(for this moment formally) the {\it mechanical energy dissipation balance}
\begin{align}\nonumber
  &\hspace*{0em}\frac{\d}{\d t}
  \int_\varOmega\!\!\!\!
  \linesunder{\frac\varrho2|\vv|^2}{kinetic}{energy}\!\!\!\!+\!\!\!
  \linesunder{\frac{\varphi(\Fe)}{\det\Fe}}{stored}{energy}\!\!\d\xx
+\!\int_\varOmega\!\!\!\!\lineunder{\DIS(\Fe,\theta;\EE(\vv)){:}\EE(\vv)
+\nu|\Nabla\EE(\vv)|^p_{_{_{}}}\!}{bulk dissipation rate}\!\!\d\xx
+\!\int_\varGamma\!\!\!\!\!\!\!\!\!\!\linesunder{\nu_\flat|\vv|^2_{_{_{}}}\!}{boundary}{dissipation rate}\!\!\!\!\!\!\!\!\d S
\\[-.6em]&\hspace{15em}
=\int_\varOmega\!\!\!\!\!\!\!\!\!\linesunder{\varrho\,\GRAVITY{\cdot}\vv_{_{_{_{}}}}\!\!}{power of}{gravity field}\!\!\!\!\!\!\!-\!\!\!\!
\linesunder{\frac{\COUPLING_{\Fe}'(\Fe,\theta)\Fe^\top\!\!\!}{\det\Fe}{:}\ee(\vv)
\!}{\ \ \ \ power of}{\ \ \ adiabatic effects}\!\!\d\xx+\!\int_\varGamma\!\!\!\!\!\!\linesunder{\!\!\TRACTION{\cdot}\vv_{_{_{_{}}}}\!\!}{power of}{traction}\!\!\!\!\!\d S\,.
\label{thermodynamic-Euler-mech-engr}
\end{align}

When we add \eq{Euler-thermodynam3} tested by 1, the adiabatic and the
dissipative heat sources cancel with those in
\eq{thermodynamic-Euler-mech-engr}. Thus we obtain (at least formally)
the {\it total energy balance }
\begin{align}
  &\!\!\!\!\!\frac{\d}{\d t}
  \int_\varOmega\!\!\!\!
  \linesunder{\frac\varrho2|\vv|^2}{kinetic}{energy}\!\!\!\!+\!\!\!\!\!
  \linesunder{\frac{\varphi(\Fe)}{\det\Fe}}{stored}{energy}\!\!\!\!
  +\!\!\!\!\!\linesunder{\OMEGA(\FF,\theta)}{heat}{energy}\!\!\d\xx
+\!\int_\varGamma\!\!\!\!\!\!\!\!\!\!\!\!\morelinesunder{\frac{\nu_\flat}2|\vv|^2_{_{_{}}}\!}{boundary}{heat production}{outflow}\!\!\!\!\!\!\!\!\!\!\d S
=\int_\varOmega\!\!\!\!\!\!\!\!\linesunder{\varrho\GRAVITY{\cdot}\vv}{power of}{gravity field}\!\!\!\!\!\!\!\!\!
\d\xx
+\!\int_\varGamma\!\!\!\!\linesunder{\!\!\TRACTION{\cdot}\vv\!\!}{power of}{traction\ }\!\!\!\!\!\!+\!\!\!\!\!\linesunder{h(\theta)}{heat}{flux}\!\!\!\d S\,,
%\nonumber\\[-1.8em]
\label{thermodynamic-Euler-engr}\end{align}
which expresses the 1st law of thermodynamics. Also, the 2nd law of
thermodynamics is satisfied, viz the {\it Clausius-Duhem inequality}
\eq{entropy-ineq} above. 
Another aspect important both thermodynamically and also for mathematical
analysis is the non-negativity of temperature, related to the 3rd law of
thermodynamics. This will be demonstrated later when we exploit some
information about the quality of the velocity field extracted from
\eq{thermodynamic-Euler-mech-engr}, cf.\ \eq{Euler-thermo-test-nonnegative} below.

\section{The analysis -- weak solutions of \eq{Euler-thermodynam}}
%        ~~~~~~~~~~~~~~~~~~~~~~~~~~~~~~~~~~~~~~~~~~~~~~~~~~~~~~~
\label{sec-anal}

We will provide a proof of existence and certain regularity of weak
solutions. To this aim, the concept of multipolar viscosity is essential
but, anyhow, still quite nontrivial and carefully ordered arguments
will be needed. The peculiarities are that 
the inertial term in the Eulerian setting involves varying mass density 
requiring sophisticated techniques from compressible fluid dynamics, 
the momentum equation is very geometrically nonlinear, and the heat equation 
has an $L^1$-structure with an $\FF$-dependent heat capacity and
with the convective time derivative, besides ever-troubling adiabatic
effects due to the necessarily general (nonlinear) coupling of
mechanical and thermal effect in the deforming configuration.
Moreover, it is well known that the so-called \emph{Lavrentiev phenomenon}
may occur in nonlinear static elasticity where the
stored energy has singularity at $\det\FF=0$ in the sense
that minimizing energy on $W^{1,\infty}(\varOmega;\R^d)$ may yield a strictly
bigger infimum than minimizing on a ``correct'' $W^{1,p}$-space, as 
pointed out in \cite{Ball87SMSE,Ball02SOPE,BalMiz85ODVP,FoHrMi03LGPN}, 
following the original old observation in \cite{Lavr26QPCV}. This brings 
difficulties in a Galerkin approximation and which are copied also for evolution
problems, while time discretization would be even more problematic
due to necessarily very nonconvex stored energy.

We will consider an initial-value problem, prescribing the initial conditions
\begin{align}\label{Euler-thermodynam-IC}
\varrho|_{t=0}^{}=\varrho_0\,,
\ \ \ \ \vv|_{t=0}^{}=\vv_0\,,\ \ \ \ 
\FF|_{t=0}^{}=\FF_0\,,\ \ \text{ and }\ \ \theta|_{t=0}^{}=\theta_0\,.
\end{align}
Referring to the referential mass density $\rhoR$, the initial conditions
should satisfy $\varrho_0=\rhoR/\!\det\FF_0$.

\subsection{Weak solutions and the main existence result}
%           ~~~~~~~~~~~~~~~~~~~~~~~~~~~~~~~~~~~~~~~~~~~~
We will use the standard notation concerning the Lebesgue and the Sobolev
spaces, namely $L^p(\varOmega;\R^n)$ for Lebesgue measurable functions
$\varOmega\to\R^n$ whose Euclidean norm is integrable with $p$-power, and
$W^{k,p}(\varOmega;\R^n)$ for functions from $L^p(\varOmega;\R^n)$ whose
all derivatives up to the order $k$ have their Euclidean norm integrable with
$p$-power. We also write briefly $H^k=W^{k,2}$. The notation
$p^*$ will denote the exponent from the embedding
$W^{1,p}(\varOmega)\subset L^{p^*}(\varOmega)$, i.e.\ $p^*=dp/(d{-}p)$
for $p<d$ while $p^*\ge1$ arbitrary for $p=d$ or $p^*=+\infty$ for $p>d$.
Moreover, for a Banach space
$X$ and for $I=[0,T]$, we will use the notation $L^p(I;X)$ for the Bochner
space of Bochner measurable functions $I\to X$ whose norm is in $L^p(I)$
while $W^{1,p}(I;X)$ denotes for functions $I\to X$ whose distributional
derivative is in $L^p(I;X)$. Also, $C(\cdot)$ and $C^1(\cdot)$
will denote spaces of continuous and continuously differentiable functions.

Moreover, as usual, we will use $C$ for a generic constant which may vary
from estimate to estimate.

To devise a weak formulation of the initial-boundary-value problem
\eq{Euler-thermodynam-BC} and \eq{Euler-thermodynam-IC} for
the system \eq{Euler-thermodynam}, we use the by-part integration in time and
the Green formula for the inertial force in the form of the left-hand
side of \eq{inertial} multiplied by a smooth test function $\widetilde\vv$. 
The by-part integration in time and
the Green formula is used also for $\DT\Fe$ in the evolution rule
\eq{Euler-thermodynam2} tested by a smooth $\widetilde\SS$ with
$\widetilde\SS(T)=0$ together $\vv{\cdot}\nn=0$, which gives
\begin{align}\nonumber
\!\int_0^T\!\!\!\int_\varOmega\!
\DT\Fe{:}\widetilde\SS\,\d\xx\d t
&=\int_0^T\!\!\!\int_\varOmega\Big(\pdt\Fe+(\vv{\cdot}\Nabla)\Fe\Big){:}\widetilde\SS\,\d\xx\d t
=\int_0^T\!\!\!\int_\varGamma(\vv{\cdot}\nn)(\Fe{:}\widetilde\SS)\,\d S\d t
\\[-.1em]&-
\int_0^T\!\!\!\int_\varOmega\!\Fe{:}\pdt{\widetilde\SS}
+({\rm div}\,\vv)\Fe{:}\widetilde\SS
+\Fe{:}((\vv{\cdot}\Nabla)\widetilde\SS)
\,\d\xx\d t
-\!\int_\varOmega\!\Fe(0){:}\widetilde\SS(0)\,\d\xx\,.
\label{one-integral-id}\end{align}

\begin{definition}[Weak solutions to \eq{Euler-thermodynam}]\label{def}
For $p,q\in[1,\infty)$, a quadruple $(\varrho,\vv,\Fe,\theta)$ with
$\varrho\in L^\infty(I{\times}\varOmega)$,
$\vv\in L^q(I;W^{1,q}(\varOmega;\R^d))\,\cap\,L^p(I;W^{2,p}(\varOmega;\R^d))$,
$\FF\in L^\infty(I{\times}\varOmega;\R^{d\times d})$,
and $\theta\in L^1(I;W^{1,1}(\varOmega))$ will be called a weak
solution to the system \eq{Euler-thermodynam} with the boundary conditions
\eq{Euler-thermodynam-BC} and the initial condition \eq{Euler-thermodynam-IC}
if $\det\Fe>0$ a.e.,
$\psi_\Fe'(\Fe,\theta)\Fe^\top\!/\!\det\Fe\in L^1(I{\times}\varOmega;\R_{\rm sym}^{d\times d})$,
$\COUPLING'_{\Fe}(\Fe,\theta)\Fe^\top\!/\!\det\Fe\in
L^{q'}(I{\times}\varOmega;\R_{\rm sym}^{d\times d})$,
$\DIS(\Fe,\theta;\ee(\vv))\in L^{q'}(I{\times}\varOmega;\R_{\rm sym}^{d\times d})$,
 such that the integral identities 
\begin{subequations}\label{Euler-weak}\begin{align}\label{Euler0-weak}
&\int_0^T\!\!\!\!\int_\varOmega\varrho\pdt{\widetilde\varrho}
+\varrho\vv{\cdot}\nabla\widetilde\varrho\,\d\xx\d t+\!\int_\varOmega\!\varrho_0
\widetilde\varrho(0)\,\d\xx=0\,,
\intertext{for any $\widetilde\varrho$ smooth with $\widetilde\varrho(T)=0$
on $\varOmega$, and}
&\nonumber
\int_0^T\!\!\!\!\int_\varOmega\bigg(\Big(\frac{
\psi_\Fe'(\Fe,\theta)\Fe^\top\!\!\!}{\det\Fe}\,+\DIS(\Fe,\theta;\ee(\vv))
-\varrho\vv{\otimes}\vv\Big){:}\ee(\widetilde\vv)
+\nu|\nabla\ee(\vv)|^{p-2}\nabla\ee(\vv)\Vdots\Nabla\ee(\widetilde\vv)
\label{Euler1-weak}\\[-.4em]&\hspace{2.5em}
-\varrho\vv{\cdot}\pdt{\widetilde\vv}\bigg)\,\d\xx\d t
=\!\int_0^T\!\!\!\!\int_\varOmega\varrho\GRAVITY{\cdot}\widetilde\vv\,\d\xx\d t
+\int_0^T\!\!\!\!\int_\varGamma(\ff{-}\nu_\flat\vv){\cdot}\widetilde\vv\,\d S\d t
+\int_\varOmega\varrho_0\vv_0{\cdot}\widetilde\vv(0)\,\d\xx
\intertext{holds for any $\widetilde\vv$ and smooth with
$\widetilde\vv{\cdot}\nn={\bm0}$ and $\widetilde\vv(T)=0$, and}
&\int_0^T\!\!\!\int_\varOmega\bigg(\Fe{:}
\pdt{\widetilde\SS}+\Big(({\rm div}\,\vv)
\Fe{+}(\Nabla\vv)\Fe\Big){:}\widetilde\SS
+\Fe{:}((\vv{\cdot}\Nabla)\widetilde\SS)\bigg)\,\d\xx\d t
+\!\!\int_\varOmega\!\FF_0{:}\widetilde\SS(0)\,\d\xx=0
\label{Euler2-weak}
\intertext{holds for any $\widetilde\SS$ smooth with
$\widetilde\SS(T)={\bm0}$, and also the integral identity 
}
\nonumber
&\int_0^T\!\!\!\int_\varOmega\!\bigg(\!\OMEGA(\FF,\theta)\pdt{\widetilde\theta}
+\big(\OMEGA(\FF,\theta)\vv{-}\kappa(\FF,\theta)\nabla\theta\big){\cdot}\nabla\widetilde\theta
\\[-.4em]&\hspace{5em}\nonumber
+\Big(\DIS(\FF,\theta;\EE(\vv)){:}\EE(\vv)+\nu|\nabla\EE(\vv)|^p\!
+
\frac{\COUPLING'_{\Fe}(\Fe,\theta)\Fe^\top\!\!\!}{\det\Fe}\,{:}\ee(\vv)
\!\Big)\widetilde\theta\bigg)\d\xx\d t
\\[-.4em]&\hspace{11em}
+\!\int_0^T\!\!\!\int_\varGamma\Big(h(\theta){+}\frac{\nu_\flat}2|\vv|^2\Big)\widetilde\theta\,\d S\d t
+\!\int_\varOmega\!\OMEGA(\FF_0,\theta_0)\widetilde\theta(0)\,\d\xx=0
\label{Euler3-weak}
\end{align}
\end{subequations}
with $\OMEGA(\cdot,\cdot)$ from \eq{Euler-thermodynam3}
holds for any $\widetilde\theta$ smooth with $\widetilde\theta(T)=0$.
\end{definition}

Before stating the main analytical result, let us summarize the data
qualification which will be, to some extent, fitted to the
examples in Section~\ref{sec-examples}. 
For some $\delta>0$, and some $1<q<\infty$ and $d<p<\infty$, we assume:
\begin{subequations}\label{Euler-ass}\begin{align}
&\varOmega\ \text{ a smooth bounded domain of $\R^d$, }\ d=2,3
\\&\label{Euler-ass-phi}
\varphi\in C^1({\rm GL}^+(d)),\ \ \inf\varphi({\rm GL}^+(d))>0\,,
\\&\nonumber
\COUPLING\,{\in}\, C^1({\rm GL}^+(d){\times}\R^+),\ \forall(\FF,\theta,\widetilde\theta)\,{\in}\,{\rm GL}^+(d){\times}\R^+{\times}\R^+{:}\ 
\\[-.3em]&\hspace{3em}
\frac{\!\COUPLING_{\theta}'(\Fe,\mm,\theta){-}\COUPLING_{\theta}'(\Fe,\mm,\widetilde\theta)\!}
{\theta{-}\widetilde\theta}
\le\frac{-\delta}{\det\FF}
\ \ \text{ and }\ \ 
\Big|\frac{\COUPLING_{\Fe}'(\Fe,\theta)\Fe^\top\!\!\!}{\det\Fe}\;\Big|
\le C\Big(1\,{+}\,\frac{\varphi(\Fe){+}\theta}{\det\Fe}
\Big),\label{Euler-ass-adiab}
\\[-.0em]&\nonumber
\forall K{\subset}\,{\rm GL}^+(d)\text{ compact }\ \exists C_K<\infty
\ \forall(\FF,\theta)\in K{\times}\R^+:\ \ \ 
\\&\hspace{7em}
\OMEGA_\theta'(\FF,\theta)+|\OMEGA_{\FF\theta}''(\FF,\theta)|\le  C_K,\ \ \
|\OMEGA_\FF'(\FF,\theta)|\le  C_K(1{+}\theta)\,,
\label{Euler-ass-primitive-c}
\\[-.0em]&\nonumber\nu>0,\ \nu_\flat>0,\ \ \text{ and }\ 
\DIS(\FF,\theta;\cdot)\ \text{ monotone}\,,\ 
\DIS\in
C({\rm GL}^+(d){\times}\R^+{\times}\R_{\rm sym}^{d\times d};\R_{\rm sym}^{d\times d})\,,
\\[-.1em]&\nonumber
\forall(\FF,\theta,\EE,\widetilde\EE)\,{\in}\,{\rm GL}^+(d){\times}\R^+{\times}(\R_{\rm sym}^{d\times d})^2:
\ \ \DIS(\FF,\theta;\bm0)=\bm0\,,
\ \ \big|\DIS(\FF,\theta;\EE)\big|\le\frac{\!1{+}|\EE|^{q-1}\!\!\!\!}\delta\,,
\\[-.0em]&\hspace{17.em}
\big(\DIS(\FF,\theta;\EE){-}\DIS(\FF,\theta;\widetilde\EE)\big){:}(\EE{-}\EE)\ge\delta|\EE{-}\widetilde\EE|^q,
\label{Euler-ass-xi}
\\&\kappa\in C({\rm GL}^+(d){\times}\R^+)\ \text{ bounded},\ \ \
\inf\kappa({\rm GL}^+(d){\times}\R^+)>0\,,
\\&\nonumber
h:I{\times}\varGamma{\times}\R^+\to\R\ \text{ Carath\'eodory function},\ \ \
\theta\,h(t,\xx,\theta)\le C(1{+}\theta^2)\ \ \text{ and}
\\[-.1em]&\hspace{11em}
h(t,\xx,\theta)\le h_{\max}(t,\xx)\ \ \text{ for some }\ h_{\max}\in L^1(I{\times}\varGamma)\,,
\label{Euler-ass-h}
\\&
{\bm g}\in L^1(I;L^\infty(\varOmega;\R^d))\,,
\ \ \ff\in L^2(I{\times}\varGamma;\R^d),\ \ \ff{\cdot}\nn=0\,,
\label{Euler-ass-f-g}
\\&\Fe_0\in W^{1,r}(\varOmega;\R^{d\times d})\,,\ \ r>d\,,\ \ \ \text{ with }\ \ \
{\rm ess\hspace*{.05em}\inf}_{\varOmega}^{}\det\Fe_0>0\,,
\label{Euler-ass-Fe0}
\\&\rhoR\in L^\infty(\varOmega)\cap W^{1,r}(\varOmega)\,,\ \ r>d\,,\ \ \text{ with }\ \ 
{\rm ess\hspace*{.05em}\inf}_{\varOmega}^{}\rhoR>0\,,
\label{Euler-ass-rhoR}
\\&\theta_0\in L^1(\varOmega),\ \ \ \theta_0\ge0\ \text{ a.e.\ on }\ \varOmega\,,
\label{Euler-ass-theta0}
\end{align}\end{subequations}
where $\OMEGA=\OMEGA(\FF,\theta)$ in 
\eq{Euler-ass-primitive-c} is from \eq{Euler-thermodynam3} and
where ${\rm GL}^+(d)=\{F\in\R^{d\times d};\ \det F>0\}$ denotes the
orientation-preserving general linear group. One should note that
\eq{Euler-ass-phi} is used in \eq{Euler-est-of-rhs} below and,
although $\varphi$ occurs in \eq{Euler-thermodynam} only through its derivative,
actually requires $\varphi$ to be bounded from below. 
Independently, the blow-up
of energy in compression, i.e.\ $\varphi(\FF)\to\infty$ if $\det\FF\to0+$,
is allowed. The first condition in \eq{Euler-ass-adiab} is just a condition
on the heat capacity $c=c(\FF,\theta)=\OMEGA_\theta'(\Fe,\theta)$ and,
in particular, implies the coercivity
$\OMEGA(\Fe,\theta)\ge \delta\theta/\!\det\FF$ since $\OMEGA(\Fe,0)=0$. The
condition \eq{Euler-ass-h} is well fitted with the standard situation that
the boundary flux is $h(\theta)=f(\theta_{\rm ext})-f(\theta)$ with an
increasing function $f$ and with $\theta_{\rm ext}\ge0$ a prescribed external
temperature, so that one can choose $h_{\max}=f(\theta_{\rm ext})$
provided we prove that
$\theta\ge0$. Also the condition $\theta\,h(t,\xx,\theta)\le C(1{+}\theta^2)$
is well compatible with this ansatz provided $f(0)=0$ and $f(\theta_{\rm ext})\in
L^2(I{\times}\varGamma)$. 

Moreover, the expected {\it symmetry} of such Cauchy stress $\TT$
is granted by {\it frame indifference} of $\psi(\cdot,\theta)$. This means that 
\begin{align}
\forall (\FF,\theta)\in{\rm GL}^+(d){\times}\R,\ \ Q\in{\rm SO}(d):\ \ \ \ 
\psi(\FF,\theta)=\psi(Q\FF,\theta)\,,
\label{frame-indifference}\end{align}
where $Q\in{\rm SO}(d)=\{Q\in\R^{d\times d};\ Q^\top Q=QQ^\top=\bbI\}$ is the
special orthogonal group.
Let us note that the frame indifference \eq{frame-indifference} 
means that $\psi(\FF,\theta)=\widehat\psi(\FF^\top\!\FF,\theta)$
for some $\widehat\psi$, which further means that
$\TT=\FF\widehat\psi_\FF'(\FF^\top\!\FF,\theta)\FF^\top\!/\!\det\FF$. Such
$\TT$ is obviously symmetric.

\begin{theorem}[Existence and regularity of weak solutions]\label{prop-Euler}
Let $\min(p,q)>d$ and the assumptions \eq{Euler-ass} and
\eq{frame-indifference} hold. Then:\\
\Item{(i)}{there exist a weak solution $(\varrho,\vv,\FF,\theta)$ according
Definition~\ref{def} with a non-negative mass density
$\varrho\in L^\infty(I;W^{1,r}(\varOmega))$ such that
$\pdt{}\varrho\in L^s(I;L^{rs/(r+s)}(\varOmega))$ with $3\le s<p(pd{+}4p{-}2d)/(4p{-}2d)$,
and a non-negative temperature $\theta\in L^\infty(I;L^1(\varOmega))
\,\cap\,L^\EXP(I;W^{1,\EXP}(\varOmega))$ with $1\le\EXP<(d{+}2)/(d{+}1)$,
and further $\pdt{}{\Fe}\in L^{\min(p,q)}(I;L^{r}(\varOmega;\R^{d\times d}))$ and
$\Nabla\Fe\in L^\infty(I;L^r(\varOmega;\R^{d\times d\times d}))$.}
\Item{(ii)}{Moreover, this solution complies with energetics in the sense that
the energy dissipation balance \eq{thermodynamic-Euler-mech-engr}
as well as the total energy balance \eq{thermodynamic-Euler-engr}
integrated over the time interval
$[0,t]$ with the initial conditions \eq{Euler-thermodynam-IC} hold.}
\end{theorem}

It should be said that the uniqueness of the weak solution is left open and,
if it ever holds, would likely need more regularity than stated above.

\subsection{Formal estimates}
%           ~~~~~~~~~~~~~~~~
Formally, the assumptions \eq{Euler-ass} yield some a-priori bounds which can
be obtained from the total energy balance \eq{thermodynamic-Euler-engr} and the 
mechanical energy-dissipation balance \eq{thermodynamic-Euler-mech-engr} for
any sufficiently regular solution $(\varrho,\vv,\FF,\theta)$ with
$\theta\ge0$ a.e.\
in $I{\times}\varOmega$. Later in Section~\ref{sec-proof}, we will
prove the existence of such solutions, but unfortunately
we are not able to claim that every weak solution has $\theta$ non-negative.

First, we use the total energy balance \eq{thermodynamic-Euler-engr}
which does not see any adiabatic and dissipative-heat terms that are
problematic as far as estimation is concerned. Assuming that $\theta\ge0$, we
have also $\OMEGA(\FF,\theta)\ge0$ and thus we are ``only'' to estimate the
right-hand side in \eq{thermodynamic-Euler-engr}. One first issue is
estimation of the gravity force $\varrho\GRAVITY$  when tested by the
velocity $\vv$, which can be estimated by the H\"older/Young
inequality as
\begin{align}\nonumber
\int_\varOmega\varrho\GRAVITY{\cdot}\vv\,\d \bm{x}
&=\int_\varOmega\sqrt{\frac{\rhoR}{\det\FF}}\sqrt{\varrho}\vv{\cdot}\GRAVITY\,\d \bm{x}
\le\Big\|\sqrt{\frac{\rhoR}{\det\FF}}\Big\|_{L^{2}(\varOmega)}
\big\|\sqrt{\varrho}\vv\big\|_{L^2(\varOmega;\R^d)}^{}
\big\|\GRAVITY\big\|_{L^\infty(\varOmega;\R^d)}^{}
\\&\nonumber\le
\frac12\bigg(\Big\|\sqrt{\frac{\rhoR}{\det\FF}}\Big\|_{L^2(\varOmega)}^2
\!+\big\|\sqrt{\varrho}\vv\big\|_{L^2(\varOmega;\R^d)}^2\bigg)
\,\big\|\GRAVITY\big\|_{L^\infty(\varOmega;\R^d)}^{}
\\&\nonumber=\big\|\GRAVITY\big\|_{L^\infty(\varOmega;\R^d)}^{}\int_\varOmega
\frac{\rhoR}{2\det\FF\!}+\frac\varrho2|\vv|^2\,\d\xx
\\&
\!\!\stackrel{\eq{Euler-ass-phi}}{=}\!\big\|\GRAVITY\big\|_{L^\infty(\varOmega;\R^d)}^{}
\bigg(\frac{\max
\rhoR(\barOmega)\!\!}{2\inf\varphi({\rm GL}^+(d))}\!\int_\varOmega
\frac{\varphi(\FF)}{\det\FF\!}\,\d\xx
+\int_\varOmega\frac\varrho2|\vv|^2\,\d\xx\bigg)\,.
\label{Euler-est-of-rhs}\end{align}
The integral on the right-hand side of \eq{Euler-est-of-rhs} can then be
treated by the Gronwall lemma, for which 
one needs the qualification \eq{Euler-ass-f-g} of $\GRAVITY$.

The boundary terms in \eq{thermodynamic-Euler-engr} can be
estimated, at the current time instant $t\in I$, as
\begin{align}
\int_\varGamma\!\TRACTION{\cdot}\vv
+h(\theta)\,\d S\le
\frac1{\nu_\flat}\|\TRACTION\|_{L^2(\varGamma;\R^d)}^2+\frac{\nu_\flat}4\|\vv\|_{L^2(\varGamma;\R^d)}^2
+\|h_{\rm max}\|_{L^1(\varGamma)}^{}\,
\end{align}
with the term  $\nu_\flat\|\vv\|_{L^2(\varGamma;\R^d)}^2/4$ to be absorbed
in the left-hand side of \eq{thermodynamic-Euler-engr}; here we used
the modelling assumption that the part of the heat produced by the
boundary viscosity leaves the system, otherwise we would have to
confine ourselves to $\ff=0$. Then, integrating 
\eq{thermodynamic-Euler-engr} in time,
one can use the qualification of $\TRACTION$ in \eq{Euler-ass-f-g}.

As a result, we obtain the (formal) a-priori estimates
\begin{subequations}\label{Euler-est}\begin{align}
&\big\|\sqrt\varrho\vv\big\|_{L^\infty(I;L^2(\varOmega;\R^d))}\le C\,,\ \ \
\label{est-rv2}
\\&\label{est-phi}
\Big\|\frac{\varphi(\Fe)}{\det\Fe}\Big\|_{L^\infty(I;L^1(\varOmega))}^{}\le C
\intertext{and, when recalling that we consider now only
solutions with $\theta\ge0$ and realizing that
$\OMEGA(\FF,\theta)\ge\delta\theta/\!\det\FF$ in \eq{Euler-ass-adiab}, also}
&\Big\|\frac{\theta}{\det\FF}\Big\|_{L^\infty(I;L^1(\varOmega))}^{}\le C\,.
\label{est-w/detF}
\end{align}\end{subequations}

Now we come to \eq{thermodynamic-Euler-mech-engr}.
The issue is now the estimation of the adiabatic term in
\eq{thermodynamic-Euler-mech-engr}. Here we exploit the frame indifference
\eq{frame-indifference} so that $\COUPLING_{\Fe}'(\Fe,\theta)\Fe^\top$ is
symmetric and thus $\COUPLING_{\Fe}'(\Fe,\theta)\Fe^\top\!{:}\nabla \vv=
\COUPLING_{\Fe}'(\Fe,\theta)\Fe^\top\!{:}\ee(\vv)$. Then we estimate
\begin{align}\nonumber
&\!\!\int_\varOmega
\Big|\frac{\!\COUPLING_{\Fe}'(\Fe,\theta)\Fe^\top\!\!\!}{\det\Fe}{:}\ee(\vv)
\Big|\,\d\xx
\le\Big\|\frac{\!\COUPLING_{\Fe}'(\Fe,\theta)\Fe^\top\!\!\!}{\det\Fe}\
\Big\|_{L^1(\varOmega;\R^{d\times d})}\|\ee(\vv)\|_{L^\infty(\varOmega;\R^{d\times d})}
\\&\nonumber
\stackrel{\eq{Euler-ass-adiab}}\le\!
C\Big\|1{+}\frac{\upvarphi(\Fe)+\theta}{\det\Fe}
\Big\|_{L^1(\varOmega)}
\|\ee(\vv)\|_{L^\infty(\varOmega;\R^{d\times d})}
\!\!\!\!\stackrel{(\ref{Euler-est}b,c)}\le\!\!
C'\Big(\|\nabla\ee(\vv)\|_{L^p(\varOmega;\R^{d\times d\times d})}^{}\!
+\big\|\vv|_\varGamma\big\|_{L^2(\varGamma;\R^d)}\Big)
\\&\ \ \le
C_{p,\nu,\nu_\flat}\!\!+\frac\nu2\|\nabla\ee(\vv)\|_{L^p(\varOmega;\R^{d\times d\times d})}^p\!
+\frac{\nu_\flat}2\big\|\vv|_\varGamma\big\|_{L^2(\varGamma;\R^d)}^2\,.
\label{Euler-est-of-rhs+}\end{align}
The last two terms in \eq{Euler-est-of-rhs+} can be absorbed in the
left-hand side of \eq{thermodynamic-Euler-mech-engr}.
Thus we obtain also the estimates 
\begin{align}
\|\ee(\vv)\|_{L^q(I{\times}\varOmega;\R^{d\times d})}^{}\le C\ \ \text{ and }\ \
\|\nabla\ee(\vv)\|_{L^p(I{\times}\varOmega;\R^{d\times d\times d})}^{}\le C\,.\ \
\label{est-e(v)}\end{align}

When assuming $p>d$, the estimate \eq{est-e(v)} is essential by preventing
the evolution of singularities of the quantities transported by such a smooth
velocity field. Here, due to the qualification of $\FF_0$ and 
$\varrho_0=\rhoR/\!\det\FF_0$ in \eq{Euler-ass-Fe0}  and \eq{Euler-ass-rhoR}, 
it yields the estimates
\begin{subequations}\label{est+}
\begin{align}
\label{est+Fes}&\|\FF\|_{L^\infty(I;W^{1,r}(\varOmega;\R^{d\times d}))}\le C_r\,,
\ \ \ \:\Big\|\frac1{\det\FF}\Big\|_{L^\infty(I;W^{1,r}(\varOmega))}\le C_r\,,
\\&\label{est+rho}\|\varrho\|_{L^\infty(I;W^{1,r}(\varOmega))}^{}\le C_r\,,
\ \ \text{ and }\ \ \Big\|\frac1\varrho\Big\|_{L^\infty(I;W^{1,r}(\varOmega))}\!\le C_r
\ \ \ \text{ for any $1\le r<+\infty$};
\intertext{cf.\ Lemma~\ref{lem-transport} used as \eq{usage} and also the
regularity of the initial conditions that follows from the assumptions (\ref{Euler-ass}i-k):}
&\nonumber
\nabla\Big(\frac1{\det\FF_0}\Big)=-\,\frac{\det'(\FF_0){:}\nabla\FF_0}
{(\det\FF_0)^2}=-\,\frac{{\rm Cof}\FF_0{:}\nabla\FF_0}{(\det\FF_0)^2}\in L^r(\varOmega;\R^d)\,,
\\&\nonumber
\nabla\varrho_0=\frac{\nabla\rhoR}{\det\FF_0}-\rhoR\frac{{\rm Cof}\FF_0{:}\nabla\FF_0}{(\det\FF_0)^2}\in L^r(\varOmega;\R^d)\,,\ \ \text{ and }\ \ 
\\&\nonumber
\nabla\Big(\frac1{\varrho_0}\Big)=\nabla\Big(\frac{\det\FF_0}{\rhoR}\Big)=
\frac{{\rm Cof}\FF_0{:}\nabla\FF_0}{\rhoR}
-\frac{\det\FF_0\nabla\rhoR}{\varrho_\text{\sc r}^2}\in L^r(\varOmega;\R^d)\,.
\intertext{
From (\ref{Euler-est}a,d) and (\ref{est+}a,b), we then have also}
&\|\vv\|_{L^\infty(I;L^2(\varOmega;\R^d))}^{}\le
\|\sqrt\varrho\vv\|_{L^\infty(I;L^2(\varOmega;\R^d))}^{}\Big\|\frac1{\sqrt\varrho}\Big\|_{L^\infty(I\times\varOmega)}^{}\le C\ \ \ \text{ and}
\label{basic-est-of-v}
\\&\|\theta\|_{L^\infty(I;L^\ALPH(\varOmega))}^{}\le
\Big\|\frac{\theta}{\det\FF}\Big\|_{L^\infty(I;L^1(\varOmega))}^{}
\|\det\FF\|_{L^\infty(I\times\varOmega)}^{}\le C\,.
\label{basic-est-of-theta}\end{align}\end{subequations}

\subsection{Proof of Theorem~\ref{prop-Euler}}\label{sec-proof}
%           ~~~~~~~~~~~~~~~~~~~~~~~~~~~~~~~~

Let us outline the main technical difficulties more in detail:
The time discretization
(Rothe's method) standardly needs convexity of $\varphi$ (which
is not a realistic assumption in finite-strain mechanics) possibly
weakened if there is some viscosity in $\FF$ (which is not
directly considered here, however). The conformal space
discretization (i.e.\ the Faedo-Galerkin method) cannot
directly copy the energetics because the ``nonlinear'' test
of \eqref{Euler-thermodynam2} by $[\varphi/\!\det]'(\FF)$ needed in
\eq{Euler-large-thermo} is problematic in this approximation
as $[\varphi/\!\det]'(\FF)$ is not in the respective finite-dimensional
space in general and similarly also the test of \eq{Euler-thermodynam0}
by $|\vv|^2$ is problematic.

For the approximation method used in the proof below, we assume
the data $\psi$, $\DIS$, $\kappa$, and $h$ to be defined also for 
the negative temperature by extending them as
\begin{align}\nonumber
&\psi(\FF,\theta):=\varphi(\FF)
+\theta({\rm ln}(-\theta)-1)
\,,\ \ \ \ \DIS(\FF,\theta;\ee):=\DIS(\FF,-\theta;\ee)\,,
\\&\kappa(\FF,\theta):=\kappa(\FF,-\theta)\,,
\ \ \ \text{ and }\ \ \ h(t,\xx,\theta):=h(t,\xx,-\theta)
\ \ \ \text{ for }\ \ \theta<0
\label{extension-negative}\end{align}
with $\varphi$ from the split \eq{ansatz}. This definition
makes $\psi:{\rm GL}^+(d)\times\R\to\R$ continuous and implies
 $\COUPLING_\FF'(\FF,\theta)={\bm0}$ and 
$\OMEGA(\FF,\theta)=\theta/\!\det\FF$ 
so that $\OMEGA_\FF'(\FF,\theta)=
\theta\FF^{-\top}\!/\!\det\FF$ for $\theta$ negative. In
particular, both $\OMEGA(\FF,\cdot)$ and $\COUPLING_\FF'(\FF,\cdot)$ are
continuous. 

For clarity, we will divide the proof into nine steps.

\medskip\noindent{\it Step 1: a regularization}.
Referring to the formal estimates \eq{est+Fes}, we can choose $\LAM>0$ so
small that, for any possible sufficiently regular solution, it holds 
\begin{align}
&\det\FF>\LAM\ \ \ \ \text{ and }\ \ \ \ |\FF|<\frac1\LAM
\ \ \text{ a.e.\ on }\ I{\times}\varOmega\,.
\label{Euler-quasistatic-est-formal4}
\end{align}
We first regularize the stress $\TT$ in \eq{Euler-thermodynam1} by considering
a smooth cut-off $\pi_\LAM\in C^1(\R^{d\times d})$ defined as
\begin{align}&\pi_\LAM(\FF)
:=\begin{cases}
\qquad\qquad1&\hspace{-8em}
\text{for $\det \FF\ge\LAM$ and $|\FF|\le1/\LAM$,}
\\
\qquad\qquad0&\hspace{-8em}\text{for $\det \FF\le\LAM/2$ or $|\FF|\ge2/\LAM$,}
\\
\displaystyle{\Big(\frac{3}{\LAM^2}\big(2\det\FF-\LAM\big)^2
-\frac{2}{\LAM^3}\big(2\det\FF-\LAM\big)^3\Big)\,\times}\!\!&
\\[.2em]
\qquad\qquad\displaystyle{\times\,\big(3(\LAM|\FF|-1)^2
-2(\LAM|\FF|-1)^3\big)}\!\!&\text{otherwise}.
\end{cases}
\label{cut-off-general}
\end{align}
Here $|\cdot|$ stands for the Frobenius norm $|\FF|=(\sum_{i,j=1}^d
F_{ij}^2)^{1/2}$ for $\FF=[F_{ij}]$, which guarantees that $\pi_\LAM$
is frame indifferent.

Furthermore, we also regularize the singular nonlinearity $1/\!\det(\cdot)$
which is still employed in the right-hand-side force in the momentum
equation, although the mass-density continuity equation is
is kept nonregularized for the inertial term. To this aim, 
we introduce the short-hand notation 
\begin{subequations}\begin{align}\label{cut-off-det}
&{\det}_\LAM(\FF):=
\min\Big(\max\Big(\!\det\FF,\frac\LAM2\,\Big),\frac2\LAM\Big)\ \ \text{ and }
\\&
\TT_{\LAM,\EPS}(\FF,\theta):=
\Big([\pi_\LAM\varphi]'(\FF)+\pi_\LAM(\FF)\frac{\COUPLING_\FF'(\FF,\theta)}{1{+}\varepsilon|\theta|
}\Big)\frac{\FF^\top}{\det\FF}\,.
\end{align}\end{subequations}
Note that also $\pi_\LAM\varphi\in C^1(\R^{d\times d})$ if
$\varphi\in C^1({\rm GL}^+(d))$ and that 
$[\pi_\LAM\varphi]'$ together with the regularized Cauchy stress $\TT_{\LAM,\EPS}$
are bounded, continuous, and vanish if $\FF$ ``substantially'' violates
the constraints \eq{Euler-quasistatic-est-formal4},
specifically:
$$
\bigg(\det\FF\le\frac\LAM2\ \ \text{ or }\ \ |\FF|\ge\frac2\LAM\bigg)
\ \ \ \Rightarrow\ \ \ \TT_{\LAM,\EPS}(\FF,\theta)={\bm0}\,.
$$
Also, $\pi_\LAM(\FF)\COUPLING_\FF'(\FF,\theta)\FF^\top\!/\!\det_\LAM(\FF)=
\pi_\LAM(\FF)\COUPLING_\FF'(\FF,\theta)\FF^\top\!/\!\det\FF$ so that, due to
\eq{Euler-ass-adiab}, $\TT_{\LAM,\EPS}:\R^{d\times d}\times\R\to\R_{\rm sym}^{d\times d}$ is
bounded. Recalling the extension \eq{extension-negative}, let us 
note that it is defined also for negative temperatures.

Altogether, for the above chosen $\LAM$ and for any $\EPS>0$, we consider
the regularized system
\begin{subequations}\label{Euler-thermo-reg}
\begin{align}\label{Euler-thermo-reg0}&
\pdt\varrho=-{\rm div}(\varrho\vv)\,,
\\ &\pdt{}(\varrho\vv)+{\rm div}(\varrho\vv{\otimes}\vv)
  ={\rm div}\big(\TT_{\LAM,\EPS}(\FF,\theta)
  {+}\DIS(\Fe,\theta;\EE(\vv)){-}{\rm div}\,\mathcal{H}\big)
+\sqrt{\frac{\rhoR\varrho}{\!\det_\LAM(\FF)\!}\!}\ \GRAVITY
\nonumber\\&\hspace*{20em}
\ \ \ \text{ with  }\ \
\mathcal{H}=\nu|\nabla\EE(\vv)|^{p-2}\nabla\EE(\vv)\,,
\label{Euler-thermo-reg1}
\\\label{Euler-thermo-reg2}
&\pdt\FF=(\Nabla\vv)\FF-(\vv{\cdot}\nabla)\FF\,,
\\&\nonumber
\pdt\W={\rm div}\big(\kappa(\FF,\theta)\nabla\theta-\W
\vv\big)
+\frac{\DIS(\FF,\theta;\EE(\vv)){:}\EE(\vv){+}\nu|\nabla\EE(\vv)|^p}
{1{+}\EPS|\EE(\vv)|^q{+}\EPS|\nabla\EE(\vv)|^p}
\\&\hspace*{12.4em}
+\frac{\pi_\LAM(\Fe)\COUPLING'_{\Fe}(\Fe,\theta)\Fe^\top\!\!
}{(1{+}\EPS|\theta|)\det\Fe}{:}\ee(\vv)
\ \ \ \text{ with }\ \ \W=\OMEGA(\Fe,\theta)\,,
\label{Euler-thermo-reg3}\end{align}\end{subequations}
where $\OMEGA(\cdot,\cdot)$ is from \eq{Euler-thermodynam3}.
We complete this system with the correspondingly regularized boundary
conditions on $I{\times}\varGamma$:
\begin{subequations}\label{Euler-thermo-reg-BC-IC}
\begin{align}
&\vv{\cdot}\nn=0,\ \ \ \ 
\big[
(\TT_{\LAM,\EPS}(\FF,\theta){+}\DD(\FF,\theta;\ee(\vv)){-}{\rm div}\mathcal{H})\nn{-}\divS\big(\mathcal{H}
\nn\big)\big]_\text{\sc t}^{}\!+\nu_\flat\vv=\ff\,,\ \ 
\label{Euler-thermo-reg-BC-IC-1}
\\&
\Nabla\ee(\vv){:}(\nn{\otimes}\nn)={\bm0}\,,\ \ \ \text{ and }\ \ \ 
\kappa(\FF,\theta)\nabla\theta{\cdot}\nn-\frac{\nu_\flat|\vv|^2}{2{+}\EPS|\vv|^2}=h_\EPS(\theta):=\frac{h(\theta)}{1{+}\varepsilon|h(\theta)|}
\label{Euler-thermo-reg-BC-IC-2}
\intertext{and initial conditions}
&\varrho|_{t=0}^{}=\varrho_0\,,
\ \ \ \ \vv|_{t=0}^{}=\vv_0\,,\ \ \ \ 
\FF|_{t=0}^{}=\FF_0\,,\ \ \ \ \theta|_{t=0}^{}=\theta_{0,\varepsilon}:=
\frac{\theta_{0}}{1{+}\varepsilon\theta_{0}}\,.
\end{align}\end{subequations}

The corresponding weak formulation of
\eq{Euler-thermo-reg}--\eq{Euler-thermo-reg-BC-IC} a'la
Definition~\ref{def} is quite straightforward, and we will not 
explicitly write it, also because it will be obvious from its Galerkin
version \eq{Euler-weak-Galerkin} below. The philosophy of the
regularization \eq{Euler-thermo-reg} is that the estimation of the
mechanical part (\ref{Euler-thermo-reg}a-c) and the thermal part 
\eq{Euler-thermo-reg3} decouples since $\TT_{\LAM,\EPS}$ is bounded
and that heat sources in the heat equation are bounded for fixed $\EPS>0$.
Simultaneously, the heat equation has a non-negative solution and,
for $\EPS\to0$, the physical a~priori estimates
are the same as the formal estimates \eq{Euler-est}--\eq{est+} and,
when taking $\LAM>0$ small to comply with \eq{Euler-quasistatic-est-formal4},
the $\LAM$-regularization becomes eventually inactive, cf.\ Step~8 below.

\medskip\noindent{\it Step 2: a semi-discretization}.
For $\varepsilon>0$ fixed, we use a spatial semi-discretization, keeping the
transport equations
\eq{Euler-thermo-reg0} and \eq{Euler-thermo-reg2} continuous
(i.e.\ non-discretised) when exploiting Lemma~\ref{lem-transport}.
More specifically, we make a conformal Galerkin approximation of
\eq{Euler-thermo-reg1} by using a collection of 
nested finite-dimensional subspaces $\{V_k\}_{k\in\N}$ whose union is dense in
$W^{2,p}(\varOmega;\R^d)$ and a conformal Galerkin approximation of
\eq{Euler-thermo-reg3} by using a collection of nested
finite-dimensional subspaces $\{Z_k\}_{k\in\N}$ whose union is dense
in $H^1(\varOmega)$. Without loss of generality, we can assume
$\vv_0\in V_1$ and $\theta_{0,\varepsilon}\in Z_1$.

Let us denoted the approximate solution of the regularized system
\eq{Euler-thermo-reg} by $(\varrho_{\EPS k},\vv_{\EPS k},\FFepsk,\theta_{\EPS k}):
I\to W^{1,r}(\varOmega)\times V_k\times W^{1,r}(\varOmega;\R^{d\times d})\times Z_k$.
Specifically, such a quadruple should satisfy
\begin{align}\label{transport-equations}
\pdt{\varrho_{\EPS k}}=-{\rm div}(\varrho_{\EPS k}\vv_{\EPS k})\ \ \ \text{ and }
\ \ \ \pdt{\FFepsk}=(\nabla\vv_{\EPS k})\FFepsk
-(\vv_{\EPS k}{\cdot}\nabla)\FFepsk
\end{align}
in the weak form like
(\ref{Euler-weak}a,c) together with the following integral identities
\begin{subequations}\label{Euler-weak-Galerkin}\begin{align}
&\nonumber
\int_0^T\!\!\!\int_\varOmega
\bigg(\Big(\TT_{\LAM,\EPS}(\FFepsk,\theta_{\EPS k})
+\DIS(\FFepsk,\theta_{\EPS k};\ee(\vv_{\EPS k}))
-\varrho_{\EPS k}\vv_{\EPS k}{\otimes}\vv_{\EPS k}\Big){:}\ee(\widetilde\vv)
\\[-.1em]&\hspace{2.5em}\nonumber
+\nu|\nabla\ee(\vv_{\EPS k})|^{p-2}\nabla\ee(\vv_{\EPS k})\Vdots
\Nabla\ee(\widetilde\vv)
-\varrho_{\EPS k}\vv_{\EPS k}{\cdot}\pdt{\widetilde\vv}
\bigg)\,\d\xx\d t
+\int_0^T\!\!\!\!\int_\varGamma\nu_\flat\vv_{\EPS k}{\cdot}\widetilde\vv\,\d S\d t
\\[-.2em]&\hspace{5em}
=\!\int_0^T\!\!\!\int_\varOmega
\sqrt{\frac{\rhoR\varrho_{\EPS k}}{\det_\LAM(\FFepsk)}}\GRAVITY{\cdot}\widetilde\vv\,\d\xx\d t
+\!\int_0^T\!\!\!\int_\varGamma\ff{\cdot}\widetilde\vv\,\d S\d t
+\int_\varOmega\varrho_0\vv_0{\cdot}\widetilde\vv(0)\,\d\xx
\label{Euler1-weak-Galerkin}
\intertext{for any $\widetilde\vv\in L^\infty(I;V_k)$ with
$\widetilde\vv{\cdot}\nn=0$ on $I{\times}\varGamma$
and $\widetilde\vv(T)={\bm0}$, and}
\nonumber
&\!\int_0^T\!\!\!\int_\varOmega\bigg(\W_{\EPS k}\pdt{\widetilde\theta}
+\big(\W_{\EPS k}\vv_{\EPS k}
{-}\kappa(\FFepsk,\theta_{\EPS k})\nabla\theta_{\EPS k}\big)
{\cdot}\nabla\widetilde\theta
\\[-.2em]&\hspace{.0em}\nonumber
+\Big(\frac{\DIS(\FFepsk,\theta_{\EPS k};\EE(\vv_{\EPS k})){:}\EE(\vv_{\EPS k})+\nu|\nabla\EE(\vv_{\EPS k})|^p\!\!}{1{+}\EPS|\EE(\vv_{\EPS k})|^q{+}\EPS|\nabla\EE(\vv_{\EPS k})|^p}
+\frac{\pi_\LAM(\FFepsk)\COUPLING'_{\Fe}(\FFepsk,\theta_{\EPS k})\FFepsk^\top
}{(1{+}\EPS|\theta_{\EPS k}|)\det\FFepsk}
{:}\ee(\vv_{\EPS k})\Big)\widetilde\theta\,\bigg)\d\xx\d t
\\[-.1em]&\hspace{0em}
+\!\int_\varOmega\!\OMEGA(\FF_0,\theta_{0,\varepsilon})\widetilde\theta(0)\,\d\xx
+\!\int_0^T\!\!\!\int_\varGamma\!\Big(h_\EPS(\theta_{\EPS k}){+}
\frac{\nu_\flat|\vv_{\EPS k}|^2}{2{+}\EPS|\vv_{\EPS k}|^2}\Big)
\widetilde\theta\,\d S\d t=0
\ \ \ \text{ with }\ \ \W_{\EPS k}=\OMEGA(\FFepsk,\theta_{\EPS k})
\label{Euler3-weak-Galerkin}
\end{align}
\end{subequations}
holding for any $\widetilde\theta\in C^1(I;Z_k)$ with $\widetilde\theta(T)=0$.

Existence of this solution is based on the standard theory of systems of
ordinary differential equations first locally in time combined here
with the abstract $W^{1,r}(\varOmega)$-valued differential equations based
on Lemma~\ref{lem-transport} used with $n=1$ and $n=d{\times}d$ for the
scalar and the tensor transport equations 
\eq{transport-equations}, and then by successive prolongation on the whole
time interval based on the $L^\infty$-estimates below.

Actually, Lemma~\ref{lem-transport} with $\vv=\vv_{\EPS k}$
and with the fixed initial conditions $\FF_0$ and $ \varrho_0$, 
defines the nonlinear operators
$\mathfrak{F}:I\times L^p(I;W^{2,p}(\varOmega;\R^d))
\to W^{1,r}(\varOmega;\R^{\d\times d})$ and
$\mathfrak{R}:I\times L^p(I;W^{2,p}(\varOmega;\R^d))
\to W^{1,r}(\varOmega)$ by 
\begin{align}
\FFepsk(t)=\mathfrak{F}\big(t,\vv_{\EPS k}\big)\ \ \text{ and }\ \
\varrho_{\EPS k}(t)=\mathfrak{R}\big(t,\vv_{\EPS k}\big)\,.
\end{align}

\medskip\noindent{\it Step 3: first a~priori estimates}.
In the Galerkin approximation, it is legitimate to use
$\widetilde\vv=\vv_{\EPS k}$ for \eq{Euler1-weak-Galerkin} and
$\widetilde\theta=\theta_{\EPS k}$ for \eq{Euler3-weak-Galerkin}.
We take the benefit from having the transport equations
\eq{transport-equations} non-discretized and, thus, we can test them
by the nonlinearities $|\vv_{\EPS k}|^2/2$ and
$[[\pi_\LAM\varphi]'/\!\det](\FFepsk)=[(\varphi\pi_\LAM'{+}\pi_\LAM\varphi')/\!\det](\FFepsk)$,
respectively. In particular, we can use the calculus \eq{Euler-large-thermo}  
and  the calculus \eq{rate-of-kinetic}--\eq{calculus-convective-in-F}
also for the semi-Galerkin approximate solution.
Specifically,  from \eq{Euler1-weak-Galerkin} tested by $\vv_{\EPS k}$,
like \eq{thermodynamic-Euler-mech-engr} we obtain the identity
\begin{align}\nonumber
  \hspace*{0em}\frac{\d}{\d t}
  \int_\varOmega\!\frac{\varrho_{\EPS k}}2|\vv_{\EPS k}|^2&+
  \frac{\pi_\LAM(\FFepsk)\varphi(\FFepsk)\!}{\det\FFepsk}
  \,\d\xx
+\!\int_\varOmega\!\DIS(\FFepsk,\theta_{\EPS k};\EE(\vv_{\EPS k})){:}\EE(\vv_{\EPS k})
+\nu|\Nabla\EE(\vv_{\EPS k})|^p\,\d\xx
+\!\int_\varGamma\!\nu_\flat|\vv_{\EPS k}|^2\,\d S\ 
\\[-.1em]\hspace{.0em}&
=\int_\varOmega\frac{\rhoR\GRAVITY{\cdot}\vv_{\EPS k}}{\det_\LAM(\FFepsk)\!\!}
-\frac{\pi_\LAM(\FFepsk)\COUPLING_{\Fe}'(\FFepsk,\theta_{\EPS k})
\FFepsk^\top}{(1{+}\EPS|\theta_{\EPS k}|
)\det\FFepsk}{:}\ee(\vv_{\EPS k})\,\d\xx
+\!\int_\varGamma\!\TRACTION{\cdot}\vv_{\EPS k}\,\d S\,.
\label{thermodynamic-Euler-mech-disc}
\end{align}
The bulk-force term can be estimated similarly as \eq{Euler-est-of-rhs}:
\begin{align}\nonumber
\int_\varOmega\sqrt{\frac{\rhoR}{\det_\LAM(\FFepsk)}}\sqrt{\varrho_{\EPS k}}\vv_{\EPS k}{\cdot}\GRAVITY\,\d \bm{x}
&\le
\Big\|\sqrt{\frac{\rhoR}{\det_\LAM(\FFepsk)}}\Big\|_{L^{2}(\varOmega)}\big\|\sqrt{\varrho_{\EPS k}}\vv_{\EPS k}\big\|_{L^2(\varOmega;\R^d)}^{}\big\|\GRAVITY\big\|_{L^\infty(\varOmega;\R^d)}^{}
\\&\nonumber\le
\big\|\GRAVITY\big\|_{L^\infty(\varOmega;\R^d)}^{}\int_\varOmega
\frac{\rhoR}{2\det_\LAM(\FFepsk)\!}+\frac{\varrho_{\EPS k}}2|\vv_{\EPS k}|^2\,\d\xx
\\&
=\big\|\GRAVITY\big\|_{L^\infty(\varOmega;\R^d)}^{}
\bigg(\frac{{\rm meas}(\varOmega)\max_{\varOmega}\rhoR\!}{2\LAM}\,
+\!\int_\varOmega\frac{\varrho_{\EPS k}}2|\vv_{\EPS k}|^2\,\d\xx\bigg)\,
\label{Euler-est-bulk-force}\end{align}
while the last term in \eq{thermodynamic-Euler-mech-disc} can be estimated
as $\TRACTION{\cdot}\vv_{\EPS k}\le|\TRACTION|^2/\nu_\flat+\nu_\flat|\vv_{\EPS k}|^2/2$.

By the Gronwall inequality, we obtain the estimates
\begin{subequations}\label{Euler-quasistatic-est1}
\begin{align}\label{Euler-quasistatic-est1-1}
&\|\EE(\vv_{\EPS k})\|_{L^{\min(p,q)}(I;W^{1,p}(\varOmega;\R^{d\times d}))}^{}\le C\ \ \ \text{ and }\ \ \
\big\|\sqrt{\varrho_{\EPS k}}\vv_{\EPS k}\big\|_{L^\infty(I;L^2(\varOmega;\R^d))}\le C\,.
\intertext{The former estimate in \eq{Euler-quasistatic-est1-1} allows us to
use Lemma~\ref{lem-transport} to obtain the estimate}
&\label{Euler-est1-2}
\|\FFepsk\|_{L^\infty(I;W^{1,r}(\varOmega;\R^{d\times d}))}^{}\le C\ \ \text{ with }\ \
\Big\|\frac1{\det\FFepsk}\Big\|_{L^\infty(I;W^{1,r}(\varOmega))}^{}\!\le C
\ \ \text{ and }
\\&\|\varrho_{\EPS k}\|_{L^\infty(I;W^{1,r}(\varOmega))}^{}\le C
\ \ \ \text{ with }\ \ \ \Big\|\frac1{\varrho_{\EPS k}}\Big\|_{L^\infty(I;W^{1,r}(\varOmega))}^{}\!\le C
\intertext{
and, like \eq{basic-est-of-v}, we also have that}
&
\|\vv_{\EPS k}\|_{L^\infty(I;L^2(\varOmega;\R^d))}^{}\le C\,.
\end{align}\end{subequations}
For $\LAM$ and $\EPS$ fixed, it is important
that these estimates can be made independently of $\theta$ since $\TT_{\LAM,\EPS}$ is
a-priori bounded. It is also important that, due to the latter estimate in
\eq{Euler-est1-2}, the singularity in $\COUPLING(\cdot,\theta)$ is
not active and $\OMEGA(\FFepsk,\theta_{\EPS k})$ is well defined.

Further estimates can be obtained by testing the Galerkin approximation
of \eq{Euler3-weak-Galerkin} by $\widetilde\theta=\theta_{\EPS k}$.
For $\W_{\EPS k}=\OMEGA(\FFepsk,\theta_{\EPS k})$ and
using \eq{transport-equations}, we will use the calculus
\begin{align}\nonumber
\theta_{\EPS k}\DT\W_{\EPS k}&=\theta_{\EPS k}
\OMEGA_\Fe'(\FFepsk,\theta_{\EPS k}){:}\DT\FFepsk
+\theta_{\EPS k}\OMEGA_\theta'(\FFepsk,\theta_{\EPS k})\DT\theta_{\EPS k}
\\&=\Big(\theta_{\EPS k}\OMEGA_\Fe'(\FFepsk,\theta_{\EPS k})
-\widehat\OMEGA_\Fe'(\FFepsk,\theta_{\EPS k})\Big){:}
(\nabla\vv_{\EPS k})\FFepsk
+\DT{\overline{\widehat\OMEGA(\FFepsk,\theta_{\EPS k})}}\,,
\label{Euler-thermodynam3-test-}
\end{align}
where the dot-notation refers to the convective time derivative here
with respect to the velocity $\vv_{\EPS k}$ and where 
$\widehat\OMEGA(\FF,\theta)$ is a primitive function to
$\theta\mapsto\theta\OMEGA_\theta'(\Fe,\theta)$ depending smoothly
on $\Fe$, specifically
\begin{align}
\widehat\OMEGA(\Fe,\theta)=\int_0^1\!\!r\theta^2\OMEGA_\theta'(\Fe,r\theta)\,\d r\,.
\label{primitive}\end{align}
Integrating the last term in \eq{Euler-thermodynam3-test-}
over $\varOmega$ gives, by the Green formula,
\begin{align}\nonumber
\int_\varOmega\DT{\overline{\widehat\OMEGA(\FFepsk,\theta_{\EPS k})}}\,\d\xx=
\frac{\d}{\d t}\!\int_\varOmega\!\widehat\OMEGA(\FFepsk,\theta_{\EPS k})\,\d\xx
-\!\int_\varOmega\!\widehat\OMEGA(\FFepsk,\theta_{\EPS k}){\rm div}\vv_{\EPS k}\,\d\xx+\!\int_\varGamma\!\widehat\OMEGA(\FFepsk,\theta_{\EPS k})\vv_{\EPS k}{\cdot}\nn\,\d S\,.
\end{align}
Thus, the mentioned test by $\theta_{\EPS k}$ then gives
\begin{align}\nonumber
&\frac{\d}{\d t}\int_\varOmega\widehat\OMEGA(\FFepsk,\theta_{\EPS k})\,\d\xx
+\int_\varOmega\kappa(\FFepsk,\theta_{\EPS k})|\nabla\theta_{\EPS k}|^2\,\d\xx
\\[-.2em]&\hspace{.5em}\nonumber
=\int_\varOmega\!\bigg(
\theta_{\EPS k}\frac{\DIS(\FFepsk,\theta_{\EPS k};\EE(\vv_{\EPS k})){:}\EE(\vv_{\EPS k})+\nu|\nabla\EE(\vv_{\EPS k})|^p\!\!}{1{+}\EPS|\EE(\vv_{\EPS k})|^q{+}\EPS|\nabla\EE(\vv)|^p}
+\theta_{\EPS k}\frac{\pi_\LAM(\FFepsk)
\COUPLING'_{\Fe}(\FFepsk,\theta_{\EPS k})
\FFepsk^\top}{(1{+}\EPS|\theta_{\EPS k}|
)\det\FFepsk}{:}\ee(\vv_{\EPS k})
\\[.1em]&\nonumber\hspace{3.5em}
-\theta_{\EPS k}\OMEGA(\FFepsk,\theta_{\EPS k})
{\rm div}\,\vv_{\EPS k}-\Big(\theta_{\EPS k}
\OMEGA_\Fe'(\FFepsk,\theta_{\EPS k})
-\widehat\OMEGA_\Fe'(\FFepsk,\theta_{\EPS k})\Big){:}
(\nabla\vv_{\EPS k})\FFepsk
\\[-.1em]&\hspace{3.5em}
+\widehat\OMEGA(\FFepsk,\theta_{\EPS k}){\rm div}\,\vv_{\EPS k}
\bigg)\,\d\xx
+\int_\varGamma\Big(h_\EPS(\theta_{\EPS k})-\frac{\nu_\flat|\vv_{\EPS k}|^2}
{2{+}\EPS|\vv_{\EPS k}|^2}\Big)\theta_{\EPS k}\,\d S\,.
\label{Euler3-Galerkin-est}\end{align}

We integrate \eq{Euler3-Galerkin-est} in time over an interval
$[0,t]$ with $t\in I$. For the left-hand side, let us realize
that $\widehat\OMEGA(\Fe,\theta)\ge c_K^{}\theta^2$
due to \eq{Euler-ass-adiab} with $c_K^{}>0$ depending
on the fixed $\LAM>0$ used in \eq{Euler-quasistatic-est-formal4}.
Then the integrated right-hand side of \eq{Euler3-Galerkin-est} is to
be estimated from above, in particular relying on \eq{Euler-ass-adiab}
and on \eq{Euler-ass-primitive-c}. Let us discuss the difficult terms.
In view of \eq{primitive}, it holds that
$\wh\OMEGA_\Fe'(\Fe,\theta)=\int_0^1r\theta^2
\OMEGA_{\Fe\theta}''(\Fe,r\theta)\,\d r$.
Recalling \eq{Euler-ass-primitive-c}, we have 
$|\wh\OMEGA_\Fe'(\Fe,\theta)
-\theta\OMEGA_\Fe'(\Fe\theta)
-\theta\OMEGA_\Fe'(\Fe,\theta)|\le C(1{+}|\theta|^2)$.
It allows for estimation 
\begin{align}
&\bigg|\int_\varOmega\Big(\wh\OMEGA_\Fe'(\FFepsk,\theta_{\EPS k})-\theta_{\EPS k}
\OMEGA_\Fe'(\FFepsk,\theta_{\EPS k})\Big){:}(\nabla\vv_{\EPS k})\FFepsk\,\d\xx\bigg|
%\nonumber\\&\hspace{10em}
\le C\big(|\varOmega|+\|\theta_{\EPS k}\|_{L^{2}(\varOmega)}^{2}\big)
\big\|(\nabla\vv_{\EPS k})\FFepsk\big\|_{L^\infty(\varOmega;\R^{d\times d})}
\label{adiab-est-F}\,.
\end{align}
Using \eq{Euler-ass-primitive-c} together with $\OMEGA(\FF,0)=0$
so that $|\OMEGA(\FF,\theta|\le C_K|\theta|$, the convective terms
$\OMEGA(\FFepsk,\theta_{\EPS k})({\rm div}\vv_{\EPS k})\theta_{\EPS k}$
and $\wh\OMEGA(\FFepsk,\theta_{\EPS k}){\rm div}\,\vv_{\EPS k}$
in \eq{Euler3-Galerkin-est} can be estimated as
\begin{align}
\!\!\!\!\!\bigg|\int_\varOmega\!\!\Big(\wh\OMEGA(\FFepsk,\theta_{\EPS k})
{-}\theta_{\EPS k}\OMEGA(\FFepsk,\theta_{\EPS k})\Big){\rm div}\vv_{\EPS k}
\d\xx\bigg|\le\frac{\!C_K}2\big\|\theta_{\EPS k}\big\|_{L^2(\varOmega)}^2
\big\|{\rm div}\,\vv_{\EPS k}\big\|_{L^\infty(\varOmega)}.\!\!
\label{Euler3-Galerkin-est++}\end{align}
The terms $\|\theta_{\EPS k}\|_{L^{2}(\varOmega)}^{2}$ in \eq{adiab-est-F} and
in \eq{Euler3-Galerkin-est++} are to be treated by the Gronwall inequality. 
The boundary term in \eq{Euler3-Galerkin-est} can be estimated by \eq{Euler-ass-h},
taking also into the account the extension \eq{extension-negative}, as 
\begin{align}
\!\int_\varGamma\!\Big(h_\EPS(\theta_{\EPS k})-\frac{\nu_\flat|\vv_{\EPS k}|^2}
{2{+}\EPS|\vv_{\EPS k}|^2}\Big)\theta_{\EPS k}\,\d S
&\le C_{\EPS,\nu_\flat,a}\!\!+a\|\theta_{\EPS k}\|_{L^2(\varGamma)}^2
%\nonumber\\&
\le C_{\EPS,\nu_\flat,a}\!\!+a N^2\big(\|\theta_{\EPS k}\|_{L^2(\varOmega)}^2
+\|\nabla\theta_{\EPS k}\|_{L^2(\varOmega;\R^d)}^2\big)\,,
\label{boundary-heat-est}\end{align}
where $C_{\EPS,\nu_\flat,\delta}$ depends also on $C$ from \eq{Euler-ass-h} and
where $N$ is the norm of the trace operator $H^1(\varOmega)\to L^2(\varGamma)$.
For $a>0$ in \eq{boundary-heat-est} sufficiently small, the last term
can be absorbed in the left-hand side of \eq{Euler3-Galerkin-est}.

By the Gronwall inequality, exploiting again the bound 
(\ref{Euler-quasistatic-est1}a,b), we obtain the estimate
\begin{subequations}\label{Euler-quasistatic-est2}
\begin{align}\label{Euler-quasistatic-est2-theta}
&\|\theta_{\EPS k}\|_{L^\infty(I;L^{\ONEALPH}(\varOmega))\,\cap\,L^2(I;H^1(\varOmega))}^{}\le C\,.
\intertext{In addition, realizing that $|\OMEGA_\FF'(\FFepsk,\theta_{\EPS k})|$
is bounded in $L^2(I;L^{2^*}(\varOmega))$ due to \eq{Euler-ass-primitive-c}
 and that $|\nabla\FFepsk|$ is bounded in
 $L^\infty(I;L^r(\varOmega))$, from the calculus
$\nabla\W_{\EPS k}=\OMEGA_\theta'(\FFepsk,\theta_{\EPS k})\nabla\theta_{\EPS k}
$ $+$ $\OMEGA_\FF'(\FFepsk,\theta_{\EPS k})\nabla\FFepsk$ $\in$ $ L^2(I{\times}\varOmega;\R^d)$,
we have also the bound}
&\|\W_{\EPS k}\|_{L^\infty(I;L^{\ONEALPH}(\varOmega))\,\cap\,L^2(I;H^1(\varOmega))}^{}\le C.
\label{Euler-weak-sln-w}
\end{align}\end{subequations}

\medskip\noindent{\it Step 4: Limit passage for $k\to\infty$}.
Using the Banach selection principle, we can extract some subsequence of
$\{(\varrho_{\EPS k},\vv_{\EPS k},\FFepsk,\W_{\EPS k})\}_{k\in\N}$
and its limit
$(\varrho_\EPS,\vv_\EPS,\FFeps,\W_\EPS):I\to W^{1,r}(\varOmega)\times
L^2(\varOmega;\R^d)\times W^{1,r}(\varOmega;\R^{d\times d})\times L^2(\varOmega)$
such that
\begin{subequations}\label{Euler-weak-sln}
\begin{align}
&\!\!\varrho_{\EPS k}\to\varrho_\EPS&&\text{weakly* in $\
L^\infty(I;W^{1,r}(\varOmega))\,\cap\,W^{1,\min(p,q)}(I;L^r(\varOmega))$}\,,
\\\label{Euler-weak-sln-v}
&\!\!\vv_{\EPS k}\to\vv_\EPS&&\text{weakly* in $\
L^\infty(I;L^2(\varOmega;\R^d))\cap
L^{\min(p,q)}(I;W^{2,p}(\varOmega;\R^d))$,}\!\!&&
\\
&\!\!\FFepsk\to\FFeps
\!\!\!&&\text{weakly* in $\ L^\infty(I;W^{1,r}(\varOmega;\R^{d\times d}))\,\cap\,
W^{1,\min(p,q)}(I;L^r(\varOmega;\R^{d\times d}))$,}\!\!
\\
&\!\!\W_{\EPS k}\to\W_\EPS&&
\text{weakly* in $\ L^\infty(I;L^{\ONEALPH}(\varOmega))\,\cap\,L^2(I;H^1(\varOmega))$}\,.
\label{Euler-weak-sln-w++}
\end{align}\end{subequations}
By the Aubin-Lions lemma here relying on the assumption $r>d$, we also have that
\begin{align}\label{rho-conv}
&\varrho_{\EPS k}\to\varrho_\EPS\ \text{ strongly in }\ C(I{\times}\barOmega) 
\ \ \text{ and }\ \ \FFepsk\to\FFeps\text{ strongly in
$C(I{\times}\barOmega;\R^{d\times d})$.}
\end{align}
This already allows for the limit passage in the evolution equations
\eq{transport-equations}, cf.\ \eq{v-mapsto-F} below.

Further, by comparison in the equation \eq{Euler-thermo-reg3} with the boundary condition
\eq{Euler-thermo-reg-BC-IC-2} in its Galerkin approximation,
we obtain a bound on $\pdt{}\W_{\EPS k}$ in seminorms $|\cdot|_l$ on $L^2(I;H^1(\varOmega)^*)$
arising from this Galerkin approximation:
$$
\big|f\big|_l^{}:=\sup\limits_{\stackrel{{\scriptstyle{\widetilde\theta(t)\in Z_l\ \text{for }t\in I}}}{{\scriptstyle{\|\widetilde\theta\|_{L^2(I;H^1(\varOmega))}^{}\le1}}}}\int_0^T\!\!\!\int_\varOmega f\widetilde\theta\,\d\xx\d t
$$
More specifically, for any $k\ge l$, we can estimate
\begin{align}\nonumber&
\Big|\pdt{\W_{\EPS k}}\Big|_l^{}=\!\!\!
\sup\limits_{\stackrel{{\scriptstyle{\widetilde\theta(t)\in Z_l\ \text{for }t\in I}}}{{\scriptstyle{\|\widetilde\theta\|_{L^2(I;H^1(\varOmega))}^{}\le1}}}}
\int_0^T\!\!\!\int_\varOmega\!\bigg(\Big(\frac{\DIS(\FFepsk,\theta_{\EPS k};\EE(\vv_{\EPS k})){:}\EE(\vv_{\EPS k})+\nu|\nabla\EE(\vv_{\EPS k})|^p\!\!}{1{+}\EPS|\EE(\vv_{\EPS k})|^q{+}\EPS|\nabla\EE(\vv)|^p}
+\frac{\pi_\LAM(\FFepsk)\COUPLING'_{\Fe}(\FFepsk,\theta_{\EPS k})\FFepsk^\top
}{(1{+}\EPS|\theta_{\EPS k}|)\det\FFepsk}{:}\ee(\vv_{\EPS k})
\Big)\widetilde\theta
\\[-1.8em]&\hspace{13em}
{-}\kappa(\FFepsk,\theta_{\EPS k})\nabla\theta_{\EPS k}\big)
{\cdot}\nabla\widetilde\theta\,\bigg)\d\xx\d t
%\nonumber\\[-.1em]&\hspace{17em}
+\!\int_0^T\!\!\!\int_\varGamma\!
\Big(h_\EPS(\theta_{\EPS k}){+}
\frac{\nu_\flat|\vv_{\EPS k}|^2}{2{+}\EPS|\vv_{\EPS k}|^2}\Big)
\widetilde\theta\,\d S\d t\le C
\label{Euler3-weak-Galerkin+}\end{align}
with some $C$ depending on the estimates (\ref{Euler-quasistatic-est1}a,b)
and \eq{Euler-quasistatic-est2-theta} but independent on $l\in N$. 
Thus, by \eq{Euler-weak-sln-w++} and by a generalized Aubin-Lions theorem
\cite[Ch.8]{Roub13NPDE}, we obtain
\begin{subequations}\label{Euler-weak+}
\begin{align}
&\label{w-conv}
\W_{\EPS k}\to \W_\EPS\hspace*{-0em}&&
\hspace*{-1em}\text{strongly in $L^s(I{\times}\varOmega)$ for $1\le s<2+4/d$}.
\intertext{Since $\OMEGA(\FFepsk,\cdot)$ is increasing, we can write 
$\theta_{\EPS k}=[\OMEGA(\FFepsk,\cdot)]^{-1}(\W_{\EPS k})$. 
Thanks to the continuity of
$(\FF,\W)\mapsto[\OMEGA(\FF,\cdot)]^{-1}(\W):\R^{d\times d}\times\R\to\R$
and the at-most linear growth in $\W$ uniformly with respect to $\FF$
from any compact $K\subset{\rm GL}^+(d)$, cf.\ \eq{Euler-ass-adiab}, we have
also}
&\label{z-conv}
\theta_{\EPS k}\to \theta_\EPS=[\OMEGA(\FFeps,\cdot)]^{-1}(\W_\EPS)
\hspace*{-0em}&&\hspace*{-1em}\text{strongly in $L^s(I{\times}\varOmega)$
for $1\le s<2+4/d$};
\intertext{note that we do not have any direct information about
$\pdt{}\theta_{\EPS k}$
so that we could not use the Aubin-Lions arguments straight for
$\{ \theta_{\EPS k}\}_{k\in\N}$. Thus, by the continuity of the corresponding
Nemytski\u{\i} (or here simply superposition) mappings, also the 
conservative part of the regularized Cauchy stress
as well as the heat part of the internal energy, namely}
&
\TT_{\LAM,\EPS}(\FFepsk,\theta_{\EPS k})\to \TT_{\LAM,\EPS}(\FFeps,\theta_\EPS)
\hspace*{-0em}&&\hspace*{0em}\text{strongly in
$L^c(I{\times}\varOmega;\R_{\rm sym}^{d\times d}),\ \ 1\le c<\infty$,}
\label{Euler-T-strong-conv}
\\&
\frac{\pi_\LAM(\FFepsk)\COUPLING'_{\Fe}(\FFepsk,\theta_{\EPS k})\FFepsk^\top
}{(1{+}\EPS|\theta_{\EPS k}|)\det\FFepsk}
%\nonumber\\&\hspace{2em}
\to \frac{\pi_\LAM(\FFeps)\COUPLING'_{\Fe}(\FFeps,\theta_\EPS)\FFeps^\top
}{(1{+}\EPS|\theta_\EPS|)\det\FFeps}
\hspace*{0em}&&\hspace*{0em}\text{strongly in
$L^c(I{\times}\varOmega;\R_{\rm sym}^{d\times d}),\ \ 1\le c<\infty$,}
\\&\OMEGA(\FFepsk,\theta_{\EPS k})\to
\OMEGA(\FFeps,\theta_\EPS)\hspace*{-0em}&&\hspace*{0em}\text{strongly in $L^c(I{\times}\varOmega),\ \ 1\le c<2+4/d$.}\label{Euler-weak+w}
\end{align}\end{subequations}
It is important that $\nabla(\varrho_{\EPS k}\vv_{\EPS k})=
\nabla\varrho_{\EPS k}{\otimes}\vv_{\EPS k}+\varrho_{\EPS k}\nabla\vv_{\EPS k}$
is bounded in $L^{\infty}(I;L^r(\varOmega;\R^{d\times d}))$ due to the already
obtained bounds (\ref{Euler-quasistatic-est1}a,c,d).
Therefore, $\varrho_{\EPS k}\vv_{\EPS k}$ converges weakly* in
$L^{\infty}(I;W^{1,r}(\varOmega;\R^d))$. The limit of
$\varrho_{\EPS k}\vv_{\EPS k}$ can be identified as $\varrho_{k}\vv_\EPS$ because
we already showed that $\varrho_{\EPS k}$ converges strongly in \eq{rho-conv}
and $\vv_{\EPS k}$ converges weakly due to \eqref{Euler-weak-sln-v}. 

By comparison, we also obtain
$\pdt{}(\varrho_{\EPS k}\vv_{\EPS k})$. Let us note, indeed, that 
\eq{inertial} still holds for the semi-discretized system since
the continuity equation has not been discretized. Specifically,
for any $\widetilde\vv\in L^2(I;W^{2,p}(\varOmega;\R^d))$ with
$\widetilde\vv(t)\in V_k$ for a.a.\ $t\in I$, we have 
\begin{align}\nonumber
\int_0^T\!\!\!\int_\varOmega
&\pdt{}(\varrho_{\EPS k}\vv_{\EPS k})\,\widetilde\vv\,\d\xx\d t
=\!\int_0^T\!\!\!\int_\varGamma\!\!\big(\ff{+}\nu_\flat\vv_{\EPS k}\big){\cdot}\widetilde\vv\,\d S\d t
+\int_0^T\!\!\!\int_\varOmega\bigg(
\sqrt{\frac{\rhoR\varrho_{\EPS k}}{\det_\LAM(\FFepsk)}}
\GRAVITY{\cdot}\widetilde\vv
\\&\hspace{0em}\nonumber
+\big(\varrho_{\EPS k}\vv_{\EPS k}{\otimes}\vv_{\EPS k}
-\TT_{\LAM,\EPS}(\FFepsk,\theta_{\EPS k})
-\DIS(\FFepsk,\theta_{\EPS k};\EE(\vv_{\EPS k})\big){:}\ee(\widetilde\vv)
\\&\hspace{0em}
-\nu|\nabla\EE(\vv_{\EPS k})|^{p-2}\nabla\EE(\vv_{\EPS k})\Vdots\nabla\EE(\widetilde\vv)\bigg)\,\d\xx\d t\le C\|\widetilde\vv\|_{L^q(I;W^{1,q}(\varOmega;\R^d))\,\cap\,L^p(I;W^{2,p}(\varOmega;\R^d))}^{}
\,,\label{est-of-DT-rho.v}
\end{align}
which yields a bound for $\pdt{}(\varrho_{\EPS k}\vv_{\EPS k})$ in a
seminorm on $L^{q'}(I;W^{1,q}(\varOmega;\R^d)^*)+L^{p'}(I;W^{2,p}(\varOmega;\R^d)^*)$
induced by the Galerkin
discretization by $V_k$ (and by any $V_l$ with $l\le k$, too). It
is important that $C$ in \eq{est-of-DT-rho.v} does not depend on $k$.
Using again a generalization of the Aubin-Lions compact-embedding theorem,
cf.\ \cite[Lemma 7.7]{Roub13NPDE}, we then obtain 
\begin{subequations}\label{rho-v-conv}\begin{align}\label{rho.v-conv}
&\varrho_{\EPS k}\,\vv_{\EPS k}\to\varrho_{k}\vv_\EPS
&&\hspace*{-1em}\text{strongly in }L^{c}(I{\times}\varOmega;\R^d)\ \ \text{ for any $1\le c<4$}\,.
\intertext{Since obviously
$\vv_{\EPS k}=(\varrho_{\EPS k}\vv_{\EPS k})(1/\varrho_{\EPS k})$,
thanks to \eq{rho-conv} and \eq{rho.v-conv}, we also have that}
&\vv_{\EPS k}\to\vv_\EPS
&&\hspace*{-1em}\text{strongly in }L^c(I{\times}\varOmega;\R^d)\ \ \text{ with
any $1\le c<4$}\,.
\end{align}\end{subequations}

For the limit passage in the momentum equation, one uses the monotonicity
of the dissipative stress, i.e., the monotonicity of the quasilinear operator
$$
\vv\mapsto{\rm div}\Big({\rm div}(|\Nabla\EE(\cdot)|^{p-2}\Nabla\EE(\vv))
-\DIS(\FF,\theta;\EE(\vv))\Big)
$$
as well as of the time-derivative operator. One could use the already obtained
weak convergences and the so-called Minty trick but, later, we will need a
strong convergence of $\ee(\vv_{\EPS k})$ to pass to the limit in the heat
equation. Thus, we first prove this strong convergence, which then allows for
the limit passage in the momentum equation directly.
We will use the weak convergence of the inertial force
\begin{align}\nonumber
&
\int_0^T\!\!\!\int_\varOmega
\Big(\pdt{}(\varrho_{\EPS k}\vv_{\EPS k})+{\rm div}(\varrho_{\EPS k}\vv_{\EPS k}{\otimes}\vv_{\EPS k})
\Big){\cdot}\widetilde\vv\,\d \xx\d t
\\[-.4em]&\ \ \ \ \nonumber=\int_\varOmega\varrho_{\EPS k}(T)\vv_{\EPS k}(T){\cdot}\widetilde\vv(T)
-\varrho_0\vv_0{\cdot}\widetilde\vv(0)\,\d \xx
-\!\int_0^T\!\!\!\int_\varOmega\varrho_{\EPS k}\vv_{\EPS k}{\cdot}\pdt{\widetilde\vv}
+(\varrho_{\EPS k}\vv_{\EPS k}{\otimes}\vv_{\EPS k}){:}\nabla\widetilde\vv\,\d \xx\d t
\\&\ \ \ \ \nonumber
{\buildrel{k\to\infty}\over{\longrightarrow}}\!
\int_\varOmega\varrho_\EPS(T)\vv_\EPS(T){\cdot}\widetilde\vv(T)
-\varrho_0\vv_0{\cdot}\widetilde\vv(0)\,\d \xx
-\!\int_0^T\!\!\!\int_\varOmega\varrho_\EPS\vv_\EPS{\cdot}\pdt{\widetilde\vv}
+(\varrho_\EPS\vv_\EPS{\otimes}\vv_\EPS){:}\nabla\widetilde\vv\,\d \xx\d t
\\&\ \ \ \ =
\int_0^T\!\!\!\int_\varOmega
\Big(\pdt{}(\varrho_\EPS\vv_\EPS)+{\rm div}(\varrho_\EPS\vv_\EPS{\otimes}\vv_\EPS)\Big){\cdot}\widetilde\vv\,\d \xx\d t\,.
\label{conv-of-inirtia.v}\end{align}
Further, relying on the calculus \eq{calculus-convective-in-F}, we will use
the identity
\begin{align}\nonumber
\!\!\int_\varOmega\!\frac{\varrho_{\EPS k}(T)}2\big|\vv_{\EPS k}(T){-}\vv_\EPS(T)\big|^2\d \xx
&=\int_0^T\!\!\!\int_\varOmega
\Big(\pdt{}(\varrho_{\EPS k}\vv_{\EPS k})
+{\rm div}(\varrho_{\EPS k}\vv_{\EPS k}{\otimes}\vv_{\EPS k})
\Big){\cdot}\vv_{\EPS k}\,\d \xx\d t
\\[-.1em]&
+\!\int_\varOmega\!\frac{\varrho_0}2|\vv_0|^2\!-\varrho_{\EPS k}(T)\vv_{\EPS k}(T){\cdot}\vv_\EPS(T)
+\frac{\!\varrho_{\EPS k}(T)}2|\vv_\EPS(T)|^2\d \xx.
\label{Euler-one-substitution}\end{align}
We further used that the $\varrho_{\EPS k}(T) $ is
also bounded in $W^{1,r}(\varOmega)$ 
and $\vv_{\EPS k}(T)$ is bounded in $ L^2(\varOmega;\R^d)$, together
with some information about the time derivative 
$\pdt{}(\varrho_{\EPS k}\vv_{\EPS k})$, cf.\ \eq{est-of-DT-rho.v}, so that
we can identify the weak limit of $\varrho_{\EPS k}(T)\vv_{\EPS k}(T)$.
Specifically, we have that
\begin{align}
&&&\varrho_{\EPS k}(T)\vv_{\EPS k}(T)\to\varrho_\EPS(T)\vv_\EPS(T)&&
\text{weakly in $\ L^2(\varOmega;\R^d)$.}&&&&
\label{Euler-weak-rho.v(T)}\end{align}

We now use the Galerkin approximation of the regularized momentum equation
\eq{Euler1-weak-Galerkin} tested by $\widetilde\vv=\vv_{\EPS k}-\widetilde\vv_k$
with $\widetilde\vv_k:I\to V_k$ an approximation of $\vv_\EPS$ in the sense that
$\widetilde\vv_k\to\vv_\EPS$ strongly in $L^\infty(I;L^2(\varOmega;\R^d))$,
$\EE(\widetilde\vv_k)\to\EE(\vv_\EPS)$ strongly in
$L^q(I{\times}\varOmega;\R^{d\times d})$, and
$\Nabla\EE(\widetilde\vv_k)\to\Nabla\EE(\vv_\EPS)$ for $k\to\infty$
strongly in $L^p(I{\times}\varOmega;\R^{d\times d\times d})$ for $k\to\infty$.
Using also the first inequality in \eq{Euler-quasistatic-est-formal4}
and \eq{Euler-one-substitution}, we can estimate
\begin{align}\nonumber
&\frac{\delta}2\big\|\vv_{\EPS k}(T){-}\vv_\EPS(T)\big\|_{L^2(\varOmega;\R^{ d})}^2
+\delta\|\EE(\vv_{\EPS k}{-}\vv_\EPS)\|_{L^q(I{\times}\varOmega;\R^{d\times d})}^q
+\nu c_{p}\|\nabla\EE(\vv_{\EPS k}{-}\vv_\EPS)\|_{L^p(I{\times}\varOmega;\R^{d\times d\times d})}^p
\\&\nonumber
\le\int_\varOmega\frac{\varrho_{\EPS k}(T)}2\big|\vv_{\EPS k}(T){-}\vv_\EPS(T)\big|^2\,\d\xx
+\int_0^T\!\!\!\int_\varGamma\nu_\flat|\vv_{\EPS k}{-}\vv_\EPS|^2\,\d S\d t
\\[-.3em]&\nonumber\hspace*{4em}
+\int_0^T\!\!\!\int_\varOmega\!\bigg(
\big(\DIS(\FFepsk,\theta_{\EPS k};\EE(\vv_{\EPS k})){-}\DIS(\FFepsk,\theta_{\EPS k};\EE(\vv_\EPS))\big)
{:}\EE(\vv_{\EPS k}{-}\vv_\EPS)
  \\[-.6em]&\hspace*{8em}\nonumber
 +\nu\big(|\nabla\EE(\vv_{\EPS k})|^{p-2}\nabla\EE(\vv_{\EPS k})
-|\nabla\EE(\vv_\EPS)|^{p-2}\nabla\EE(\vv_\EPS)\big)\Vdots
  \nabla\EE(\vv_{\EPS k}{-}\vv_\EPS)\bigg)\,\d\xx\d t
 \\[-.5em]&=\nonumber
 \int_0^T\!\!\!\int_\varOmega\bigg(\sqrt{\frac{\rhoR\varrho_{\EPS k}\!\!}{\det_{\LAM}(\FFepsk)\!\!}}\ \GRAVITY{\cdot}(\vv_{\EPS k}{-}\widetilde\vv_k)
  -\TT_{\LAM,\EPS}(\FFepsk,\theta_{\EPS k}){:}\ee(\vv_{\EPS k}{-}\widetilde\vv_k)
 \\[.1em]&\nonumber\hspace{4em}
 -\DIS(\FF_{\EPS k},\theta_{\EPS k};\EE(\widetilde\vv_k))
{:}\EE(\vv_{\EPS k}{-}\widetilde\vv_k)
-\nu\big(|\nabla\EE(\widetilde\vv_k)|^{p-2}\nabla\EE(\widetilde\vv_k)\big)\Vdots
 \nabla\EE(\vv_{\EPS k}{-}\widetilde\vv_k)
 \\[-.1em]&\nonumber\hspace*{4em}+
\Big(\pdt{}(\varrho_{\EPS k}\vv_{\EPS k})
+{\rm div}(\varrho_{\EPS k}\vv_{\EPS k}{\otimes}\vv_{\EPS k})\Big)
{\cdot}\widetilde\vv_k\bigg)\,\d\xx\d t
+\int_0^T\!\!\!\int_\varGamma\big(\TRACTION{+}\nu_\flat\widetilde\vv_k\big)
{\cdot}(\vv_{\EPS k}{-}\widetilde\vv_k)\,\d S\d t
\\[-.3em]&\hspace*{4em}
+\int_\varOmega\frac{\varrho_0}2|\vv_0|^2-\varrho_{\EPS k}(T)\vv_{\EPS k}(T){\cdot}\widetilde\vv_k(T)
+\frac{\varrho_{\EPS k}(T)}2|\widetilde\vv_k(T)|^2\,\d\xx+\mathscr{O}_k
\ {\buildrel{k\to\infty}\over{\longrightarrow}}\ 0\,
\label{strong-hyper+}\end{align}
with $\delta>0$ from \eq{Euler-ass-xi} and with some $c_{p}>0$ related to the inequality
$c_{p}|G-\widetilde G|^p\le(|G|^{p-2}G-|\widetilde G|^{p-2}\widetilde G)\Vdots(G-\widetilde G)$
holding for $p\ge2$. The term $\mathscr{O}_k$ in \eq{strong-hyper+} is
\begin{align}\nonumber
\mathscr{O}_k&=\int_\varOmega\frac{\varrho_{\EPS k}(T)\!}2\,\vv_{\EPS k}(T)
{\cdot}\big(\widetilde\vv_k(T){-}\vv_\EPS(T)\big)\,\d\xx
+\int_0^T\!\!\!\int_\varGamma\nu_\flat\vv_{\EPS k}{\cdot}(\widetilde\vv_k{-}\vv_\EPS)\,\d S\d t
\\[-.3em]&\nonumber\hspace*{1em}
+\int_0^T\!\!\!\int_\varOmega\!
\DIS(\FFepsk,\theta_{\EPS k};\EE(\vv_{\EPS k})){:}\EE(\widetilde\vv_k{-}\vv_\EPS)
 +\nu|\nabla\EE(\vv_{\EPS k})|^{p-2}\nabla\EE(\vv_{\EPS k})\Vdots
  \nabla\EE(\widetilde\vv_k{-}\vv_\EPS)\,\d\xx\d t
\end{align}
and it converges to zero due to the strong approximation properties of the
approximation $\widetilde\vv_k$ of $\vv_\EPS$. Here we used
\eq{conv-of-inirtia.v}--\eq{Euler-one-substitution} and also the strong
convergence \eq{rho-conv}, \eq{z-conv}, and \eq{rho-v-conv}.
Knowing already \eq{Euler-T-strong-conv} and that
$\ee(\vv_{\EPS k}{-}\widetilde\vv_k)\to0$ weakly in
$L^p(I{\times}\varOmega;\R^{d\times d})$,
we have that $\int_0^T\!\!\!\int_\varOmega\TT_{\LAM,\EPS}(\FFepsk,\theta_{\EPS k})
{:}\ee(\vv_{\EPS k}{-}\widetilde\vv_k)\,\d\xx\d t\to0$. 
Thus, we obtain the desired strong convergence
\begin{subequations}\label{strong-conv}\begin{align}
&\ee(\vv_{\EPS k})\to\ee(\vv_{\EPS})&&\text{strongly in 
$\,L^q(I{\times}\varOmega;\R_{\rm sym}^{d\times d})\,$ and}\\
&\nabla\ee(\vv_{\EPS k})\to\nabla\ee(\vv_{\EPS})\hspace{-2em}&&\text{strongly in 
$\,L^p(I{\times}\varOmega;\R^{d\times d\times d})$\,,}
\intertext{and also of $\vv_{\EPS k}(T)\to\vv_{\EPS}(T)$ in $L^2(\varOmega;\R^d)$.
In fact, executing this procedure for
a current time instants $t$ instead of $T$, we obtain}
&\vv_{\EPS k}(t)\to\vv_\EPS(t)&&\text{strongly in $\,L^2(\varOmega;\R^d)\,$ for any $t\in I$.}&&
\intertext{It also implies, by continuity of the trace operator
$L^2(I;H^1(\varOmega))\to L^2(I{\times}\varGamma)$, that}
&\vv_{\EPS k}\big|_{I\times\varGamma}\to\vv_\EPS\big|_{I\times\varGamma}&&
\text{strongly in $\,L^2(I{\times}\varGamma;\R^d)$}\,.
\end{align}\end{subequations}
Having \eq{strong-conv} at disposal, the limit passage in the Galerkin-approximation 
of \eq{Euler1-weak-Galerkin} to the weak solution of \eq{Euler-thermo-reg1}
is then easy.

Strong convergence (\ref{strong-conv}a,b,d) 
allows for passing to the limit also in the dissipative heat sources
and the other terms in the Galerkin approximation 
of the heat equation \eq{Euler3-weak-Galerkin} are even easier.

For further purposes, let us mention that the energy dissipation
balance \eq{thermodynamic-Euler-mech-disc} is inherited in the limit, i.e.
\begin{align}\nonumber
  &\hspace*{0em}\frac{\d}{\d t}
  \int_\varOmega\!\frac{\varrho_\EPS}2|\vv_\EPS|^2+
  \frac{\pi_\LAM(\FFeps)\varphi(\FFeps)\!}{\det\FFeps}
  \,\d\xx
+\!\int_\varOmega\!\DIS(\FFeps,\theta_\EPS;\EE(\vv_\EPS)){:}\EE(\vv_\EPS)
+\nu|\Nabla\EE(\vv_\EPS)|^p\,\d\xx
\\[-.1em]&\hspace{1.5em}
+\!\int_\varGamma\nu_\flat|\vv_\EPS|^2\,\d S
=\int_\varOmega\frac{\rhoR\GRAVITY{\cdot}\vv_\EPS}{\det_\LAM(\FFeps)\!\!}
-\frac{\pi_\LAM(\FFeps)\COUPLING_{\Fe}'(\FFeps,\theta_\EPS)\FFeps^\top}
{(1{+}|\theta_\EPS|)\det\FFeps}{:}\ee(\vv_\EPS)\,\d\xx
+\!\int_\varGamma\TRACTION{\cdot}\vv_\EPS\,\d S\,.
\label{thermodynamic-Euler-mech+}
\end{align}

\medskip\noindent{\it Step 5 -- non-negativity of temperature}:
We can now perform various nonlinear tests of the regularized but
non-discretized heat equation. The first test can be by the negative part
of temperature $\theta_\EPS^-:=\min(0,\theta_\EPS)$. Let us recall the
extension \eq{extension-negative}, which in particular gives
$\OMEGA(\FF_{\EPS},\theta^-_{\EPS})=\theta^-_{\EPS}$ and 
$\OMEGA_\FF'(\FF_{\EPS},\theta^-_{\EPS})=0$. Note also that
$\theta_\EPS^-\in L^2(I;H^1(\varOmega))$,
so that it is indeed a legal test for \eq{Euler-thermo-reg3}.
Here we rely on the data qualification $\nu,\nu_\flat\ge0$,
$\kappa=\kappa(\FF,\theta)\ge0$, $\DIS(\FF,\theta;\EE){:}\EE\ge0$,
$\theta_0\ge0$, and $h(\theta)\ge0$ for $\theta\le0$, cf.\
(\ref{Euler-ass}f--h,l).
Realizing that $\nabla\theta^-\!=0$ wherever $\theta>0$ so that
$\nabla\theta{\cdot}\nabla\theta^-=|\nabla\theta^-|^2$ and that
$\COUPLING'_{\Fe}(\Fe,\theta)\theta^-=\COUPLING'_{\Fe}(\Fe,\theta^-)\theta^-=0$
and $h(\theta)\theta^-=h(\theta^-)\theta^-=0$, this test gives
\begin{align}\nonumber
&\frac12\frac{\d}{\d t}\|\theta_\EPS^-\|_{L^2(\varOmega)}^2\le
\int_\varOmega\theta_\EPS^-\pdt{\W_\EPS}+
\kappa(\FFeps,\theta_\EPS)\nabla\theta_\EPS{\cdot}\nabla\theta_\EPS^-\,\d\xx=
\int_\varGamma\Big(h(\theta_\EPS){+}\frac{\nu_\flat|\vv_\EPS|^2}{2{+}\EPS|\vv_\EPS|^2}\Big)\theta_\EPS^-\,\d S
\\[-.1em]&\nonumber\
+\!\int_\varOmega\!\W_\EPS\vv_\EPS{\cdot}\nabla\theta_\EPS^-\!
+\Big(\!\DIS(\FFeps,\theta_\EPS;\EE(\vv_\EPS)){:}\EE(\vv_\EPS)+\nu|\nabla\EE(\vv_\EPS)|^p\!
+\frac{\pi_\LAM(\FFeps)\COUPLING'_{\Fe}(\FFeps,\theta_\EPS)\FFeps^\top\!\!}{(1{+}\EPS|\theta_\EPS|)\det\FFeps}{:}
\ee(\vv_\EPS)\Big)\theta_\EPS^-\,\d\xx
\\[-.1em]\nonumber
&\qquad\le\int_\varOmega\W\vv_\EPS{\cdot}\nabla\theta_\EPS^-\,\d\xx
=\int_\varOmega\theta_\EPS^-\vv_\EPS{\cdot}\nabla\theta_\EPS^-\,\d\xx=
-\int_\varOmega\theta_\EPS^-\vv_\EPS{\cdot}\nabla\theta_\EPS^-+|\nabla\theta_\EPS^-|^2{\rm div}\,\vv_\EPS\,\d\xx
\\[-.1em]
&\hspace{12em}
=-\frac12\int_\varOmega|\nabla\theta_\EPS^-|^2{\rm div}\,\vv_\EPS\,\d\xx
\le
\|\theta_\EPS^-\|_{L^2(\varOmega)}^2\|{\rm div}\,\vv_\EPS\|_{L^\infty(\varOmega)}^{}\,.
\label{Euler-thermo-test-nonnegative}
\end{align}
Recalling the assumption $\theta_0\ge0$ so that $\theta_{0,\EPS}^-=0$
and exploiting the information
$\vv_\EPS\in L^{\min(p,q)}(I;W^{1,p}(\varOmega;\R^d))$ with $p>d$ inherited
from \eq{Euler-quasistatic-est1-1}, by the Gronwall inequality we obtain
$\|\theta_\EPS^-\|_{L^\infty(I;L^2(\varOmega))}=0$, so that $\theta_\EPS\ge0$
a.e.\ on $I{\times}\varOmega$. 

Having proved the non-negativity of temperature, we can now execute
the strategy based of the $L^1$-theory for the heat equation
which led to the estimates \eq{est-e(v)}--\eq{est+}, i.e.\ here
\begin{subequations}\label{est-eps}\begin{align}
&\|\vv_\EPS\|_{L^\infty(I;L^2(\varOmega;\R^d))}^{}\le C,\ \ \
\|\ee(\vv_\EPS)\|_{L^q(I{\times}\varOmega;\R^{d\times d})}^{}\le C,\ \ \
\|\nabla\ee(\vv_\EPS)\|_{L^p(I{\times}\varOmega;\R^{d\times d\times d})}^{}\le C,\ \
\label{est-e(v)-eps}
\\&\label{est+Fes-eps}\|\FFeps\|_{L^\infty(I;W^{1,r}(\varOmega;\R^{d\times d}))}\le C_r\,,
\ \ \ \:\Big\|\frac1{\det\FFeps}\Big\|_{L^\infty(I;W^{1,r}(\varOmega))}\le C_r\,,
\\&\label{est+rho-eps}\|\varrho_\EPS\|_{L^\infty(I;W^{1,r}(\varOmega))}^{}\le C_r\,,
\ \ \ \ \Big\|\frac1{\varrho_\EPS}\Big\|_{L^\infty(I;W^{1,r}(\varOmega))}\!\le C_r
\ \ \ \text{ for any $1\le r<+\infty$},
\\&\|\W_\EPS\|_{L^\infty(I;L^1(\varOmega))}^{}\le C\,,\ \ \text{ and }\ \ 
\|\theta_\EPS\|_{L^\infty(I;L^\ALPH(\varOmega))}^{}\le C\,.
\label{basic-est-of-theta-eps}
\intertext{By interpolation exploiting the Gagliardo-Nirenberg inequality
between $L^2(\varOmega)$ and $W^{2,p}(\varOmega)$, we have
$\|\cdot\|_{L^\infty(\varOmega)}^{}\le C\|\cdot\|_{L^2(\varOmega)}^{r}
\|\cdot\|_{W^{2,p}(\varOmega)}^{1-r}$ with $0<r<pd/(pd+4p-2d)$.
Using also Korn's inequality, from \eq{est-e(v)-eps}
we thus obtain the estimate}
&\label{est+v}
\|\vv_\EPS\|_{L^s(I;L^\infty(\varOmega;\R^d))}^{}\le C_s\ \ \ \text{ with }\
1\le s<p\,\frac{pd{+}4p{-}2d}{4p-2d}\,.
\intertext{By comparison from
$\pdt{}\varrho_\EPS=({\rm div}\vv_\EPS)\varrho_\EPS-\vv_\EPS{\cdot}\nabla\varrho_\EPS$
and from $\pdt{}\FFeps=(\nabla\vv_{\EPS})\FF_{\EPS}
-(\vv_{\EPS}{\cdot}\nabla)\FF_{\EPS}$, we also have}
&\label{est+dF/dt}
\Big\|\pdt{\varrho_\EPS}\Big\|_{L^{\min(p,q)}(I;L^r(\varOmega))}^{}\!\le C
\ \ \ \text{ and }\ \ \ 
\Big\|\pdt{\FFeps}\Big\|_{L^{\min(p,q)}(I;L^{r}(\varOmega;\R^{d\times d}))}^{}\!\le C\,.
\end{align}\end{subequations}
The estimates \eq{basic-est-of-theta-eps} are naturally weaker than
\eq{Euler-quasistatic-est2} but, importantly,
are uniform with respect to $\EPS>0$, in contrast to
\eq{Euler-quasistatic-est2} which is not uniform in this sense.
The total energy balance \eq{thermodynamic-Euler-engr} holds for
$\EPS$-solution only as an inequality because the heat sources do not
exactly cancel; more in detail, while the regularized adiabatic heat again
cancels, the dissipative heat terms are regularized (and smaller) in
\eq{Euler-thermo-reg3} and in \eq{Euler-thermo-reg-BC-IC-2} but the
corresponding viscous stress in \eq{Euler-thermo-reg1} and force in
\eq{Euler-thermo-reg-BC-IC-1} are not regularized. This inequality still
allows to execute the above-mentioned estimation.

Let us also note that the extension \eq{extension-negative} becomes now
inactive and we can work with the original data defined for
non-negative $\theta$ only.

\medskip\noindent{\it Step 6 -- further a-priori estimates}:
We are to prove an estimate of $\nabla\theta_\EPS$ based on the test of the heat
equation \eq{Euler-thermo-reg3} by $\chi_\zeta(\theta_\EPS)$ with an increasing
nonlinear function $\chi_\zeta:[0,+\infty)\to[0,1]$ defined as
\begin{align}\label{test-chi}
\chi_\zeta(\theta):=1-\frac1{(1{+}\theta)^\zeta}\,,\ \ \ \ \zeta>0\,,
\end{align}
simplifying the original idea of L.\,Boccardo and T.\,Gallou\"et 
\cite{BDGO97NPDE,BocGal89NEPE} in the spirit of \cite{FeiMal06NSET},
expanding the estimation strategy in \cite[Sect.\,8.2]{KruRou19MMCM}.
Importantly, here we have $\chi_\zeta(\theta_\EPS(t,\cdot))\in H^1(\varOmega)$,
hence it is a legal test function, because 
$0\le\theta_\varepsilon(t,\cdot)\in H^1(\varOmega)$ has already been proved
and because $\chi_\zeta$ is Lipschitz continuous on $[0,+\infty)$. 

We consider $1\le \EXP<2$ and estimate the $L^\EXP$-norm  of $\nabla\theta_\varepsilon$
by H\"older's inequality as 
\begin{align}\nonumber
&\!\!\int_0^T\!\!\!\int_\varOmega|\nabla\theta_\varepsilon|^\EXP\,\d\xx\d t
\le C_1\bigg(\underbrace{\int_0^T
\big\|1{+}\theta_\varepsilon(t,\cdot)\big\|^{(1+\zeta)\EXP/(\TWO-\EXP)}
_{L^{(1+\zeta)\EXP/(\TWO-\EXP)}(\varOmega)}\,\d t}_{\displaystyle\ \ \ =:I_{\EXP,\zeta}^{(1)}(\theta_\varepsilon)}\bigg)^{1-\EXP/\TWO}
\bigg(\underbrace{\int_0^T\!\!\!\int_\varOmega
\chi_\zeta'(\theta_\EPS)|\nabla\theta_\EPS|^\TWO}
_{\displaystyle\ \ \ =:I_{\zeta}^{(2)}(\theta_\varepsilon)}\bigg)^{\EXP/\TWO}.
\\[-2em]\label{8-**-+}
\end{align}
with $\chi_\zeta$ from \eq{test-chi} so that
$\chi_\zeta'(\theta)=\zeta/(1{+}\theta)^{1+\zeta}$ 
and with a constant $C_1$ dependent
on $\zeta$, $\EXP$, and $T$. Then we interpolate the Lebesgue space
$L^{(1+\zeta)\EXP/(\TWO-\EXP)}(\varOmega)$ between  $W^{1,\EXP}(\varOmega)$ and 
$L^1(\varOmega)$ to exploit the already obtained
$L^\infty(I;L^1(\varOmega))$-estimate  in \eq{basic-est-of-theta-eps}.
More specifically, by the Gagliardo-Nirenberg inequality, we obtain
\begin{align}
\big\|1{+}\theta_\varepsilon(t,\cdot)\big\|_{L^{\EXP/\sigma}(\varOmega)}^{\mu/\sigma}
\le C_2\Big(1+\big\|\nabla\theta_\varepsilon(t,\cdot)\big\|_{L^\EXP(\varOmega;\R^d)}\Big)^\EXP
\ \ \ \text{ with }\ \sigma=\frac{\TWO{-}\EXP}{1{+}\zeta}
\label{8-cond}
\end{align}
with $C_2$ depending on $\sigma$, $C_1$, and $C$ from \eq{basic-est-of-theta-eps},
so that $I_{\EXP,\zeta}^{(1)}(\theta_\varepsilon)\le C_3(1+\int_0^T\!\int_\varOmega\big|\nabla \theta_\varepsilon\big|^\EXP\,\d\xx\d t)$
with $C_3$ depending on $C_2$. Combining it with \eq{8-**-+}, we obtain
\begin{align}\|\nabla\theta_\varepsilon\|_{L^\EXP(I\times\Omega;\R^d)}^\EXP=C_1C_3\big(1+\|\nabla\theta_\varepsilon\|_{L^\EXP(I\times\Omega)}^\EXP\big)^{1-\EXP/2}I_{\EXP,\zeta}^{(2)}(\theta_\varepsilon)^{\EXP/2_{_{_{}}}}_{^{^{^{}}}}\,.
\label{8-***}
\end{align}
Furthermore, we estimate $I_{\zeta}^{(2)}(\theta_\varepsilon)$ in \eq{8-**-+}.
Let us denote by ${\mathcal X}_\zeta$ a primitive function to
$\theta\mapsto\chi_\zeta(\theta)\OMEGA_\theta'(\Fe,\theta)$ depending
smoothly on $\Fe$, specifically
\begin{align}
{\mathcal X}_\zeta(\Fe,\theta)
=\int_0^1\!\!\theta\chi_\zeta(r\theta)\OMEGA_\theta'(\Fe,r\theta)\,\d r\,.
\label{primitive+}\end{align}

Like \eq{Euler-thermodynam3-test-} but using the partial (not convective) time
derivative, we have now the calculus
\begin{align}\nonumber
&\!\!\!\int_\varOmega\!\chi_\zeta(\theta)\pdt\W\,\d\xx
=\!\int_\varOmega\!\chi_\zeta(\theta)\OMEGA_\theta'(\Fe,\theta)\pdt{\theta}
+\chi_\zeta(\theta)\OMEGA_\Fe'(\Fe,\theta){:}\pdt{\Fe\!}\,\,\d\xx
=\frac{\d}{\d t}\int_\varOmega{\mathcal X}_\zeta(\Fe,\theta)\,\d\xx
\\&\qquad\qquad
-\int_\varOmega\big[{\mathscr X}_\zeta\big]_\Fe'(\FF,\theta){:}\pdt{\Fe}\,\d\xx
\ \ \text{ where }\
{\mathscr X}_\zeta(\FF,\theta):={\mathcal X}_\zeta(\Fe,\theta)
-\chi_\zeta(\theta)\OMEGA(\Fe,\theta)\,.
\label{Euler-thermodynam3-test--}\end{align}
In view of \eq{primitive+}, it holds $[{\mathscr X}_\zeta]_\Fe'(\Fe,\theta)
=\int_0^1\theta\chi_\zeta(r\theta)\OMEGA_{\Fe\theta}''(\Fe,r\theta)\,\d r
-\chi_\zeta(\theta)\OMEGA_{\Fe}'(\Fe,\theta)$.
Altogether, testing \eq{Euler-thermo-reg3} with
\eq{Euler-thermo-reg-BC-IC-2} by $\chi_\zeta(\theta_\EPS)$ gives
\begin{align}\nonumber
&\frac{\d}{\d t}\int_\varOmega\!{\mathcal X}_\zeta(\FFeps,\mm_\EPS,\theta_\EPS)\,\d\xx
+\int_\varOmega
\chi_\zeta'(\theta_\EPS)\kappa(\FFeps,\theta_\EPS)|\nabla\theta_\EPS|^2\,\d\xx
=\!\int_\varGamma\!\Big(h_\EPS(\theta_\EPS){+}\frac{\nu_\flat|\vv_\EPS|^2}{2{+}\EPS|\vv_\EPS|^2}\Big)\chi_\zeta(\theta_\EPS)\,\d S
\\&\nonumber\qquad 
+\!\int_\varOmega\!\bigg(
\frac{\DIS(\FFeps,\theta_\EPS;\EE(\vv_\EPS))\,{:}\,\EE(\vv_\EPS)+
\nu|\nabla\EE(\vv_\EPS)|^p}{1{+}\EPS|\EE(\vv_\EPS)|^q{+}\EPS|\nabla\EE(\vv_\EPS)|^p\!}
\,\chi_\zeta(\theta_\EPS)
+\OMEGA(\FFeps,\theta_\EPS)\chi_\zeta'(\theta_\EPS)\vv_\EPS{\cdot}\nabla\theta_\EPS
\\&\qquad\qquad\qquad\quad
+\big[{\mathscr X}_\zeta\big]_\Fe'(\FFeps,\theta_\EPS){:}\pdt{\FFeps}
+\chi_\zeta(\theta_\EPS)\frac{\pi_\LAM(\FFeps)\COUPLING'_{\Fe}(\FFeps,\theta_\EPS)\FFeps^\top
{:}\ee(\vv_\EPS)}{(1{+}\EPS\theta_\EPS)\det\FFeps}
\bigg)\,\d\xx\,.
\label{Euler-thermodynam3-test+++}\end{align}
We realize that $\chi_\zeta'(\theta)=\zeta/(1{+}\theta)^{1+\zeta}$ as used
already in \eq{8-**-+} and that ${\mathcal X}_\zeta(\FFeps,\theta_\EPS)\ge
c_K\theta_\EPS$ with some $c_K$ for $\theta_\EPS\ge0$ due to
\eq{Euler-ass-adiab}; again $K$ is a compact subset of ${\rm GL}^+(d)$
related here with the already proved estimates \eq{est+Fes-eps}.
The convective term in \eq{Euler-thermodynam3-test+++} can be estimated,
 for any $\delta>0$, as
\begin{align}\nonumber
\int_\varOmega\W_\EPS\chi_\zeta'(\theta_\EPS)\vv_\EPS{\cdot}\nabla\theta_\EPS\,\d\xx
&\le\frac1\delta\int_\varOmega\chi_\zeta'(\theta_\EPS)|\vv_\EPS|^2\W_\EPS^2\,\d\xx
+\delta\int_\varOmega\chi_\zeta'(\theta_\EPS)|\nabla\theta_\EPS|^\TWO\,\d\xx
\\&=\frac1\delta\int_\varOmega\chi_\zeta'(\theta_\EPS)|\vv_\EPS|^2
\W_\EPS^2\,\d\xx+\delta I_{\zeta}^{(2)}(\theta_\EPS)\,.
\label{est-of-convective}\end{align}
Denoting by $0<\kappa_0=\inf_{\FF,\theta}\kappa(\FF,\theta)$ and using
\eq{Euler-thermodynam3-test+++} integrated over $I=[0,T]$,
we further estimate:
\begin{align}\nonumber
&
I_{\zeta}^{(2)}(\theta_\EPS)=
\frac1\zeta\int_0^T\!\!\!\int_\varOmega
\chi_\zeta'(\theta_\EPS)|\nabla\theta_\EPS|^\TWO\,\d\xx\d t
\le\frac1{\kappa_0\zeta}\int_0^T\!\!\!\int_\varOmega\!\kappa(\FFeps,\theta_\EPS)
\nabla\theta_\EPS{\cdot}\nabla\chi_\zeta(\theta_\EPS)\,\d\xx\d t 
\\&\nonumber
\le\frac1{\kappa_0\zeta}\bigg(
\int_0^T\!\!\!\int_\varOmega\!\kappa(\FFeps,\theta_\EPS)
\nabla\theta_\EPS{\cdot}\nabla\chi_\zeta(\theta_\EPS)
\,\d\xx\d t+\int_\varOmega\!{\mathcal X}_\zeta(\FFeps(T),\theta_\EPS(T))\,\d\xx\bigg)
\\&\nonumber
=\frac1{\kappa_0\zeta}\bigg(\int_\varOmega\!{\mathcal X}_\zeta(\Fe_0,\theta_{0,\EPS})\,\d\xx
+\!\int_0^T\!\!\!\int_\varOmega\!
\bigg(\xi_\EPS\big(\FFeps,\theta_\EPS;\EE(\vv_\EPS),\nabla^2\vv_\EPS,\rr_\EPS\big)
\chi_\zeta(\theta_\EPS)
\\&\nonumber\ \ 
+\big[{\mathscr X}_\zeta\big]_\Fe'(\FFeps,\theta_\EPS){:}\pdt{\FFeps}
+\OMEGA(\FFeps,\theta_\EPS)\chi_\zeta'(\theta_\EPS)\vv_\EPS{\cdot}\nabla\theta_\EPS
\\&\nonumber\ \ 
+\frac{\pi_\LAM(\FFeps)\COUPLING'_{\Fe}(\FFeps,\theta_\EPS)\FFeps^\top
{:}\ee(\vv_\EPS)}{(1{+}\EPS\theta_\EPS
)\det\FFeps}\chi_\zeta(\theta_\EPS)
\bigg)\,\d\xx
+\!\int_0^T\!\!\!\int_\varGamma\!\Big(h_\EPS(\theta_\EPS)
+\frac{\nu_\flat|\vv_\EPS|^2}{2{+}\EPS|\vv_\EPS|^2}\Big)
\chi_\zeta(\theta_\EPS)\,\d S\d t\bigg)
\\[-.2em]&\nonumber
\!\!\stackrel{\eq{Euler-thermodynam3-test+++}}{\le}\!\!
\frac1{\kappa_0\zeta}\bigg(
\big\|{\mathcal X}_\zeta(\Fe_0,\theta_{0,\EPS})\big\|_{L^1(\varOmega)}\!
+\big\|\DIS(\FFeps,\theta_\EPS;\EE(\vv_\EPS)){:}\EE(\vv_\EPS)+
\nu|\nabla\EE(\vv_\EPS)|^p\big\|_{L^1(I\times\varOmega)}\!
\\&\quad\nonumber
+\int_0^T\!\!\big\|\big[{\mathscr X}_\zeta\big]_\Fe'(\FFeps,\theta_\EPS)
\big\|_{L^{r'}(\varOmega;\R^{d\times d})}^{r'}\!
+\Big\|\pdt{\FFeps}\Big\|_{L^r(\varOmega;\R^{d\times d})}^r\d t
\\&\qquad\qquad\nonumber
+\Big\|\frac{\pi_\LAM(\FFeps)\COUPLING'_{\Fe}(\FFeps,\theta_\EPS)\FFeps^\top}
{\det\FFeps}{:}\ee(\vv_\EPS)\Big\|_{L^1(I\times\varOmega)}\!\!
+\big\|h_{\max}\big\|_{L^1(I\times\varGamma)}\!
\\[-.2em]&
\qquad\qquad\qquad
+\frac1\delta\|\vv_\EPS\|_{L^2(I;L^\infty(\Omega;\R^d))}^{\TWOprime}
\|\chi_\zeta'(\theta_\EPS)\OMEGA^{\TWOprime}(\FFeps,\theta_\EPS)
\|_{L^\infty(I;L^1(\varOmega))}^{}\!\bigg)
+\frac\delta{\kappa_0}I_{\zeta}^{(2)}(\theta_\EPS);
\label{+++}
\end{align}
noteworthy, when choosing $\delta<\kappa_0$, we can absorb the last
term in the left-hand side. Due to the assumption \eq{Euler-ass-adiab},
we can estimate the adiabatic rates
$\pi_\LAM(\FFeps)\COUPLING'_{\Fe}(\FFeps,\theta_\EPS)\FFeps^\top{:}\ee(\vv_\EPS)/\!\det\FFeps$ 
in \eq{+++}, cf.\ \eq{Euler-est-of-rhs+}. We also use the estimate
\eq{est-eps} and the assumption \eq{Euler-ass-h} relying on the already
proved non-negativity of temperature. By the qualification \eq{Euler-ass-primitive-c},
we have $|[{\mathscr X}_\zeta]_\Fe'(\FF,\theta)|\le C(1{+}\theta)$.
This allows for estimation 
\begin{align}\nonumber
&\big\|\big[{\mathscr X}_\zeta\big]_\Fe'(\FFeps,\theta_\EPS)
\big\|_{L^{r'}(\varOmega;\R^{d\times d})}^{r'}
\le C^{r'}\big\|1+\theta_\EPS\big\|_{L^{r'}(\varOmega)}^{r'}
\\&\hspace{9em}\le C_4+C_4\|\theta_\EPS\|_{L^1(\varOmega)}^{r'(1-{\EXP^*}'/r)}\big(\|\theta_\EPS\|_{L^1(\varOmega)}^{}
\!{+}\|\nabla\theta_\EPS\|_{L^\EXP(\varOmega;\R^d)}\big)^{{\EXP^*}'/(r-1)}\!\,,
\label{est-adiabatic}\end{align}
where we use the Gagliardo-Nirenberg inequality to interpolate
$L^{r'}(\varOmega)$ between $L^1(\varOmega)$ and $W^{1,\EXP}(\varOmega)$.
The penultimate term in \eq{+++} is a-priori bounded independently of
$\EPS$ for $\zeta>0$ fixed because, as
$\OMEGA(\FF,\theta)=\mathscr{O}(\theta)$ due to \eq{Euler-ass-primitive-c}
and due to $\chi_\zeta'(\theta)=\mathscr{O}(1/\theta)$ uniformly for $\zeta>0$,
so that we have
$[\chi_\zeta'(\cdot)\OMEGA^\TWO(\FF,\cdot)](\theta)=\mathscr{O}(\theta)$.
Thus the estimate \eq{basic-est-of-theta-eps} guarantees
$\chi_\zeta'(\theta_\EPS)\OMEGA^{\TWOprime}(\FFeps,\theta_\EPS)$ bounded
in $L^\infty(I;L^1(\Omega))$ while $|\vv_\EPS|^{\TWOprime}$
is surely bounded in $L^1(I;L^\infty(\Omega))$, cf.\ \eq{est+v}.

In view of \eq{est-adiabatic}, one can summarize \eq{+++} as
$I_{\zeta}^{(2)}(\theta_\EPS)\le C(1+\|\nabla\theta_\EPS\|_{L^\EXP(\varOmega;\R^d)}\big)^{{\EXP^*}'/(r-1)}$. Combining it with \eq{8-***}, one obtain the inequality as
\begin{align}\|\nabla\theta_\varepsilon\|_{L^\EXP(I\times\Omega;\R^d)}^{}
\le C\big(1+\|\nabla\theta_\varepsilon\|_{L^\EXP(I\times\Omega;\R^d)}^{1-\EXP/2+{\EXP^*}'/(2r-2)}\big)\,.
\label{final-est-nabla-theta}
\end{align}
Reminding the choice 
$\sigma:=(\TWO{-}\EXP)/(1{+}\zeta)$ from \eq{8-cond} with $\zeta>0$
arbitrarily small, one gets after some algebra, the condition $\EXP<(d{+}2)/(d{+}1)$.
Obviously, for $r$ big enough (in particular if $r>d$ as assumed),
the exponent  on the right-hand side is lower than 1,
which gives a bound for $\nabla\theta_\varepsilon$ in $L^\EXP(I\times\Omega;\R^d)$.
Altogether, we proved
\begin{subequations}\label{est-W-eps}\begin{align}
&\|\theta_\EPS\|_{L^\infty(I;L^1(\varOmega))\,\cap\,L^\EXP(I;W^{1,\EXP}(\varOmega))}^{}\le C_\EXP
\ \ \text{ with }\ 1\le\EXP<\frac{d{+}2}{d{+}1}\,.
\intertext{Exploiting the calculus
$\nabla\W_{\EPS}=\OMEGA_\theta'(\FF_{\EPS},\theta_{\EPS})\nabla\theta_{\EPS}+
\OMEGA_\FF'(\FF_{\EPS},\theta_{\EPS})\nabla\FF_{\EPS}$
with $\nabla\FF_{\EPS}$ bounded in $L^\infty(I;L^r(\varOmega;\R^{d\times d\times d}))$
and relying on the assumption \eq{Euler-ass-primitive-c}, we have also the
bound on $\nabla\W_{\EPS}$ in $L^\EXP(I;L^{\EXP^*d/(\EXP^*+d)}(\varOmega;\R^d))$, so
that}
&\|\W_\EPS\|_{L^\infty(I;L^1(\varOmega))\,\cap\,L^\EXP(I;W^{1,\EXP^*d/(\EXP^*+d)}(\varOmega))}^{}\le C_\EXP\,.
\end{align}\end{subequations}

\medskip\noindent{\it Step 7: Limit passage for $\EPS\to0$}.
We use the Banach selection principle as in Step~4, now also taking 
\eq{est-eps} and \eq{est-W-eps} into account instead of
the estimates \eq{Euler-quasistatic-est1} and \eq{Euler-quasistatic-est2}.
For some subsequence and some $(\varrho,\vv,\FF,\theta)$, we now have
\begin{subequations}\label{Euler-weak++}
\begin{align}\nonumber
&\varrho_\EPS\to\varrho&&\hspace*{-11em}
\text{weakly* in $\ L^\infty(I;W^{1,r}(\varOmega))\,\cap\,
 W^{1,\min(p,q)}(I;L^2(\varOmega))$}
\\[-.2em]\label{Euler-weak-rho}
&&&\hspace*{-2em}\text{and strongly in  $C(I{\times}\barOmega)$}\,,\\
&\vv_\EPS\to\vv&&\hspace*{-11em}\text{weakly* in $\
L^\infty(I;L^2(\varOmega;\R^d))\cap
L^{ q}(I;W^{2,p}(\varOmega;\R^d))$,}\!\!&&
\label{Euler-weak-v}\\\nonumber
&\FFeps\to\FF\!\!\!&&\hspace*{-11em}\text{weakly* in $\
L^\infty(I;W^{1,r}(\varOmega;\R^{d\times d}))\,\cap\,
W^{1,\min(p,q)}(I;L^\infty(\varOmega;\R^{d\times d}))$}\!\!
\\[-.2em]&&&\hspace*{-2em}\text{and strongly in
$C(I{\times}\barOmega;\R^{d\times d})$, }
\label{Euler-weak-F}
\\%\nonumber
&\theta_\EPS\to\theta\!\!\!&&\hspace*{-11em}
\text{weakly* in $\ L^\EXP(I;W^{1,\EXP}(\varOmega)),\ 1\le\EXP<(d{+2})/(d{+}1)$.}
\intertext{Like \eq{w-conv}, by the Aubin-Lions theorem, we now have}
&\W_\EPS\to\W\!\!\!&&\hspace*{-11em}\text{strongly in $L^c(I{\times}\varOmega),\ \
1\le c<1{+}2/d$,}
\intertext{and then, using again continuity of
$(\FF,\W)\mapsto[\OMEGA(\FF,\cdot)]^{-1}(\W)$ as in \eq{z-conv}, we also have}
&\theta_\EPS\to\theta=[\OMEGA(\FF,\cdot)]^{-1}(\W)\!\!\!&&
\hspace*{-11em}\text{strongly in $\ L^c(I{\times}\varOmega),\ \ 1\le c<1{+}2/d$.}
\intertext{By the continuity of $\pi_\LAM\varphi_\FF'$, $\COUPLING_\FF'$,
$\det(\cdot)$, and $\kappa$, we have also}
&\label{m-strongly+}
\kappa(\FFeps,\theta_\EPS)\to\kappa(\FF,\theta)
&&\hspace*{-11em}\text{strongly in $L^c(I{\times}\varOmega)$ for any $1\le c<\infty$, and}
\\
&\TT_{\LAM,\EPS}(\FFeps,\theta_\EPS)\to\TT_\LAM=
\frac{[\pi_\LAM\varphi]'(\FF)\!}{\det\FF}\FF^\top\!\!
%\nonumber\\&\hspace*{7em}
+\frac{\!\pi_\LAM(\FF)\COUPLING_\FF'(\FF,\theta)\!}{\det\FF}\FF^\top
\hspace*{.5em}&&\hspace*{-1em}\text{strongly in
$L^{c}(I{\times}\varOmega;\R_{\rm sym}^{d\times d}),\ 1\le c<1{+}2/d$\,.}
\label{Euler-weak-stress}\end{align}\end{subequations}

The momentum equation \eq{Euler-thermo-reg1} (still regularized by
$\varepsilon$) is to be treated like in Step~4. Here we exploit 
the information about $\pdt{}(\varrho_\EPS\vv_\EPS)$ in 
$L^{q'}(I;W^{1,q}(\varOmega;\R^d)^*)+L^{p'}(I;W^{2,p}(\varOmega;\R^d)^*)$
obtained like in \eq{est-of-DT-rho.v}; here we used also \eq{est+v}.
By the Aubin-Lions compact-embedding theorem, we then obtain 
\begin{align}\label{rho-v-conv+}
&\varrho_{\EPS k}\,\vv_{\EPS k}\to\varrho_\EPS\vv_\EPS
&&\hspace*{-1em}\text{strongly in }L^s(I{\times}\varOmega;\R^d)\ \ \text{ with
$s$ from \eq{est+v}}\,.
\end{align}
In fact, the argumentation \eq{strong-hyper+} is to be slightly modified by
using $(\varrho_\EPS,\vv_\EPS,\TT_{\LAM,\EPS}(\FFeps,\vv_\EPS),\theta_\EPS)$ in place
of $(\varrho_{\EPS k},\vv_{\EPS k},\TT_{\LAM,\EPS}(\FFepsk,\vv_{\EPS k}),\theta_{\EPS k})$
and with $\widetilde\vv_k$ replaced by $\vv$. Specifically,
the convergence \eq{Euler-weak-stress} is also weak* in
$L^\infty(I;L^1(\varOmega;\R_{\rm sym}^{d\times d}))$ and $\ee(\vv_\EPS)\to\ee(\vv)$
strongly in $L^q(I;L^\infty(\varOmega;\R_{\rm sym}^{d\times d}))$ so that
$\int_0^T\!\int_\varOmega\TT_{\LAM,\EPS}(\FFeps,\theta_\EPS){:}\ee(\vv_\EPS{-}\vv)\,\d\xx\d t\to0$. Also,
$\int_0^T\!\int_\varOmega\pdt{}(\varrho_{\EPS k}\vv_{\EPS k}){\cdot}\widetilde\vv_k\,\d\xx\d t$
is to be replaced by $\langle\pdt{}(\varrho_\EPS\vv_\EPS),\vv\rangle$
with $\langle\cdot,\cdot\rangle$ denoting here the duality between
$L^{q'}(I;W^{1,q}(\varOmega;\R^d)^*)+L^{p'}(I;W^{2,p}(\varOmega;\R^d)^*)$
and $L^q(I;W^{1,q}(\varOmega;\R^d))\cap L^p(I;W^{2,p}(\varOmega;\R^d))$.

Limit passage in the heat equation \eq{Euler-thermo-reg3} is then simple. 
Altogether, we proved that $(\varrho,\vv,\FF,\theta)$ solves in the weak
sense the problem \eq{Euler-thermo-reg}--\eq{Euler-thermo-reg-BC-IC}
with $\EPS=0$ and with $\TT_\LAM$ from \eq{Euler-weak-stress} in place of
$\TT_{\LAM,\EPS}(\FFeps,\theta_\EPS)$.

\medskip\noindent{\it Step 8: the original problem}.
Let us note that  the limit $\Fe$ lives in 
$L^\infty(I;W^{1,r}(\varOmega;\R^{d\times d}))\,\cap\,
W^{1,\min(p,q)}(I;L^{r}(\varOmega;\R^{d\times d}))$, cf.\ (\ref{est-eps}b,f),
and this space  is embedded  into 
$C(I{\times}\barOmega;\R^{d\times d})$ if $r>d$. Therefore $\Fe$ and its
determinant evolve continuously in time, being valued respectively in
$C(\barOmega;\R^{d\times d})$ and $C(\barOmega)$.
Let us recall that the initial condition $\FF_0$ complies with the bounds
\eq{Euler-quasistatic-est-formal4} and we used this $\FF_0$
also for the $\LAM$-regularized system.
Therefore $\Fe$ satisfies these bounds not only at $t=0$ but also at least
for small times. Yet, in view of the choice \eq{Euler-quasistatic-est-formal4}
of $\LAM$, this means that the $\LAM$-regularization is nonactive
and  $(\varrho,\vv,\Fe,\theta)$ solves, at least for a small time, the
original nonregularized problem
\eq{Euler-thermo-reg}--\eq{Euler-thermo-reg-BC-IC}
for which the a~priori $L^\infty$-bounds \eq{est+} hold.
By the continuation argument, we may see that the $\LAM$-regularization 
remains therefore inactive within the whole evolution of
$(\varrho,\vv,\Fe,\theta)$ on the whole time interval $I$.

\medskip\noindent{\it Step 9: energy balances}.
It is now important that the tests and then all the subsequent calculations
leading to the energy balances \eq{thermodynamic-Euler-mech-engr} and
\eq{thermodynamic-Euler-engr} integrated over a current
time interval $[0,t]$ are really legitimate.

More in detail, in the calculus \eq{Euler-large-thermo}, we rely on that 
$[\varphi(\FF)/\!\det\FF]'\in L^\infty(I{\times}\varOmega;\R^{d\times d})$
is surely in duality with
$\pdt{}\FF\in L^{\min(p,q)}(I;L^{r}(\varOmega;\R^{d\times d}))$ and
$(\vv{\cdot}\nabla)\FF\in L^s(I;L^{r}(\varOmega;\R^{d\times d}))$
with $s$ from \eq{est+v}. Moreover, $\pdt{}(\varrho\vv)\in
L^{q'}(I;W^{1,q}(\varOmega;\R^d)^*)+L^{p'}(I;W^{2,p}(\varOmega;\R^d)^*)$
is in duality with
$\vv\in L^q(I;W^{1,q}(\varOmega;\R^d))\cap L^p(I;W^{2,p}(\varOmega;\R^d))$,
as used in \eq{calculus-convective-in-F}. Further, the calculus
\eq{rate-of-kinetic} relies on that $\pdt{}\varrho$ and ${\rm div}(\varrho\vv)=
\vv{\cdot}\nabla\varrho+\varrho\,{\rm div}\,\vv$ live in $L^s(I;L^{rs/(r+s)}(\varOmega))$
and thus are surely in duality with $|\vv|^2\in
L^{s/2}(I;L^\infty(\varOmega))$ with $3\le s<p(pd{+}4p{-}2d)/(4p{-}2d)$, cf.\
\eq{est+v}. Eventually, since
$\Nabla\ee(\vv)\in L^{p}(I{\times}\varOmega;\R^{d\times d\times d})$,
we have ${\rm div}^2(\nu|\Nabla\ee(\vv)|^{p-2}\Nabla\ee(\vv))\in
L^{p'}(I;W^{2,p}(\varOmega;\R^d)^*)$ in
duality with $\vv$. Also ${\rm div}\DIS(\FF,\theta;\ee(\vv))\in
L^{q'}(I;W^{1,q}(\varOmega;\R^d)^*)$ is in duality with $\vv$ due to the
growth condition \eq{Euler-ass-xi}. Altogether, the calculations
\eq{Euler-large-thermo}--\eq{calculus-convective-in-F} are legitimate.

\begin{remark}[{\sl ``Physical'' versus ``mathematical'' estimates}]
\label{rem-qualif}\upshape
The above results hold for arbitrarily large time horizon $T$. Anyhow,
mathematical arguments vitally relied on the regularity of the initial 
conditions, which intuitively is an information gradually forgetting when 
$T\to\infty$. This is reflected by the fact that the
$W^{1,r}(\varOmega)$-regularity of $\varrho$ and $\FF$ blows up when
$T\to\infty$. Anyhow, some physically relevant estimates for the autonomous
thermodynamic system are indeed uniform in time, specifically those which
come from the energy balances \eq{thermodynamic-Euler-mech-engr} and
\eq{thermodynamic-Euler-engr}, i.e.\ \eq{Euler-est} and \eq{est-e(v)}.
In contrast, the estimates \eq{est+} depend on the assumed regularity of the
initial conditions $\varrho_0$ and $\FF_0$ through the hyper-viscosity $\nu$,
which is surely an analytically important ``mathematical quality'' of the
solutions but does not have a direct and sustainable 
physical relevance. The same qualitative difference is relevant for
hyper-viscosity, which is sometimes considered as a controversial modelling
aspect, and therefore the asymptotics for $\nu\to0$ is relevant even if the
limit for $\nu=0$ is analytically open. More specifically, the former estimates
in \eq{Euler-est} and \eq{est-e(v)} are uniform with respect
to $\nu\to0$, in contrast to the latter
estimates in \eq{Euler-est} and \eq{est+}. A good message is that also the
total-energy equality \eq{thermodynamic-Euler-engr} is uniform with respect
to $\nu$, being not explicitly dependent on $\nu$, in
contrast to the energy-dissipation balance \eq{thermodynamic-Euler-mech-engr}.
Integrated in time, this balance involves the dissipated-energy term
$\int_0^T\!\int_\varOmega\nu|\nabla\ee(\vv)|^p\,\d x\d t$ and it rises an
interesting question whether this hyperviscosity dissipated energy (which
is a-priori solely bounded) is indeed small if $\nu>0$ is small. This is a
very nontrivial problem and the positive answer would need more regularity
of solutions than in Theorem~\ref{prop-Euler}(ii); in a linear small-strain
model see \cite[Remark~2]{Roub23GTLV}.
\end{remark}

\section{Some examples}\label{sec-examples}
%        ~~~~~~~~~~~~~

The assumptions \eq{Euler-ass} are not easy to satisfy. It is thus
worth illustrating them with some examples motivated by several specific
phenomena.

\def\EXPANSION{\alpha}

\begin{example}[{\sl Neo-Hookean thermally expanding materials}]
\label{exa-neo-Hook}\upshape
For the elastic bulk modulus $K_\text{\sc e}^{}$ and the elastic
shear modulus $G_\text{\sc e}^{}$ and a volumetric thermal expansion function 
$\EXPANSION:[0,+\infty)\to(0,+\infty)$, and example of the free energy is
\begin{align}\label{neo-Hookean}
\psi(\Fe,\theta)=\!\!\!\!\lineunder{\frac12K_\text{\sc e}^{}\big(\det\Fe-1\big)^2\!+
\frac12G_\text{\sc e}^{}\Big(\frac{{\rm tr}(\Fe\Fe^\top)}{(\det\Fe)^{2/d}}-d\Big)}
{$=\varphi(\FF)$}
\!\!\!\!\!+\!\!\!\!
\lineunder{ c_{0}\theta(1{-}{\rm ln}\theta)-K_\text{\sc e}^{}\EXPANSION(\theta)\det\Fe_{_{_{_{_{}}}}}\!\!}
{$=\COUPLING(\Fe,\theta)^{^{^{^{}}}}$}
\end{align}
for $\det\FF>0$ otherwise $\psi(\FF,\theta)=+\infty$. Here
$c_{0}>0$ is the (referential) heat capacity. 
To comply with the ansatz \eq{ansatz}, it should hold $\EXPANSION(0)=0$.
The stored energy $\varphi$ in \eq{neo-Hookean}, considered per
referential volume, is the standard isothermal neo-Hookean model.
The minimum of $\psi(\cdot,\theta)$ with respect to $\det\FF$ is at
$\det\FF=1+\EXPANSION(\theta)$, which shows the role of the function
$\EXPANSION$ as the {\it volumetric thermal expansion}. The ``coupling stress''
is
\begin{align}
\frac{\COUPLING'_{\Fe}(\Fe,\theta)\Fe^\top\!\!\!}{\det\Fe}
=-K_\text{\sc e}^{}\frac{\EXPANSION(\theta){\rm Cof}\,\FF}{\det\Fe}\Fe^\top
=- K_\text{\sc e}^{}\EXPANSION(\theta)\bbI\,,
\label{neo-Hookean++++}\end{align}
where we again used the formulas $(\det F)'={\rm Cof}F$ and
$F^{-1}={\rm Cof}F^\top\!/\!\det F$, so that the power of the adiabatic effects
caused by thermal expansion in \eq{thermodynamic-Euler-mech-engr} is
\begin{align}
\frac{\COUPLING_{\Fe}'(\Fe,\theta)\Fe^\top\!\!\!}{\det\Fe}{:}\ee(\vv)
=-K_\text{\sc e}^{}\EXPANSION(\theta){\rm div}\,\vv\,.
\end{align}
The thermal part of the actual internal energy $\W=\OMEGA(\Fe,\theta)$ is
\begin{align}
\OMEGA(\Fe,\theta)=\frac{ c_{0}\theta}{\det\FF}
-K_\text{\sc e}^{}(\EXPANSION(\theta)-\theta\EXPANSION'(\theta))
\label{neo-Hookean++}\end{align}
and the (actual) heat capacity $c(\Fe,\theta)=\OMEGA_\theta'(\Fe,\theta)$ is
\begin{align}
c(\Fe,\theta)=\frac{ c_{0}}{\det\FF}
+\theta K_\text{\sc e}^{}\EXPANSION''(\theta)\,.
\label{neo-Hookean+++}\end{align}
Always, $\EXPANSION$ should be convex to ensure positivity of the heat
capacity, cf.\ the first condition in \eq{Euler-ass-adiab}. Yet,
\eq{neo-Hookean++++} fulfills the second condition in \eq{Euler-ass-adiab}
only if $\EXPANSION$ is bounded on $\R^+$, which could be fulfilled if
$\EXPANSION$ is decreasing, but it is not much physical in real materials.
More realistically, for an increasing bounded $\EXPANSION$, the positivity
of $c(\Fe,\theta)$ can be ``locally'' achieved if $\det\FF$ is bounded
from above, i.e.\ not too much stretch of the body. This last attribute can
be expected only for ``reasonable'' loading regimes and will intuitively be
violated for extremely big loading regimes causing drastic volumetric
stretching.
\end{example}

\begin{example}[{\sl Thermal expansion: an alternative ansatz}]\label{exa-neo-Hook-alt}\upshape
Another standard concept is a multiplicative decomposition of the deformation
gradient $\FF=\FF_\text{\sc e}\FF_\text{\sc t}$ with the elastic strain
$\FF_\text{\sc e}$ and the thermal-expansion strain
$\FF_\text{\sc t}=\FF_\text{\sc t}(\theta)=(1{+}\EXPANSION(\theta))\bbI$.
Here, $\EXPANSION(\cdot)$ means a {\it length thermal expansion}. The
neo-Hookean stored energy should then be a function of the elastic strain
$\FF_\text{\sc e}=\FF/(1{+}\EXPANSION(\theta))$. Realizing that
$\FF_\text{\sc t}$ influences only the spherical but not deviatoric part of
$\FF_\text{\sc e}$ and that
$\det\FF_\text{\sc t}(\theta)=(1{+}\EXPANSION(\theta))^d$, we arrive at
\begin{align}\nonumber
\!\psi(\Fe,\theta)&=\frac12K_\text{\sc e}^{}\Big(\frac{\det\Fe}{(1{+}\EXPANSION(\theta))^d}-1\Big)^2\!+\frac12G_\text{\sc e}^{}\Big(\frac{{\rm tr}(\Fe\Fe^\top)}{(\det\Fe)^{2/d}}-d\Big)+ c_{0}\theta(1{-}{\rm ln}\theta)
\\&
=\!\!\!\!\!\lineunder{\frac12K_\text{\sc e}^{}\big(\det\Fe{-}1\big)^2\!+
\frac12G_\text{\sc e}^{}\Big(\frac{{\rm tr}(\Fe\Fe^\top)}{(\det\Fe)^{2/d}}-d\Big)}
{$=\varphi(\FF)$}
\!\!\!\!\!
\nonumber\\[-.5em]&\hspace*{4em}
+\!\!\!\!\!
\lineunder{ c_{0}\theta(1{-}{\rm ln}\theta)-K_\text{\sc e}^{}
\frac{(1{+}\EXPANSION(\theta))^{2d}\!-1}{2(1{+}\EXPANSION(\theta))^{2d}}\det\FF^2
+K_\text{\sc e}^{}\frac{(1{+}\EXPANSION(\theta))^{d}\!-1}{2(1{+}\EXPANSION(\theta))^{d}}\det\FF
\!}{$=\COUPLING(\Fe,\theta)^{^{^{^{}}}}$}\!\!\!\!.
\label{neo-Hookean-alt}\end{align}
One analytical advantage is that the coupling stress 
\begin{align}\nonumber
\frac{\COUPLING_{\Fe}'(\Fe,\theta)\Fe^\top\!\!\!}{\det\Fe}=
K_\text{\sc e}^{}\Big(\frac{(1{+}\EXPANSION(\theta))^{2d}\!-1\!}{(1{+}\EXPANSION(\theta))^{2d}}\det\Fe
+\frac{(1{+}\EXPANSION(\theta))^{d}\!-1}{2(1{+}\EXPANSION(\theta))^{d}}\Big)\bbI
\end{align}
complies with the growth condition \eq{Euler-ass-adiab}
 since $\varphi$ has a sufficiently fast-growing response on
stretching, specifically
 here $\det\FF/\varphi(\FF)=\mathscr{O}(\det\FF)$ for $\det\FF\to+\infty$.
 On the other hand, like in Example~\ref{exa-neo-Hook}, the heat capacity
$-\theta\COUPLING_{\theta\theta}''(\FF,\theta)$ can be positive
for ``reasonable'' loading regimes which do not lead to too drastic
volumetric stretch.
\end{example}

\begin{example}[{\sl Phase transitions}]\label{rem-SMA}\upshape 
An interesting example of the free energy $\psi$ occurs in (somewhat
simplified) modelling of austenite-martensite transition in so-called
{\it shape-memory allows} (SMA):
\begin{align}\nonumber
\psi(\FF,\theta)=(1{-}\lambda(\theta))\varphi_{_{\rm A}}(\FF)
+\lambda(\theta)\varphi_{_{\rm M}}(\FF)+c_{0}\theta(1{-}{\rm ln}\theta)
=\!\!\!\!\!\lineunder{\varphi_{_{\rm A}}(\FF)}
{$=\varphi(\FF)^{^{^{}}}$}\!\!\!\!\!+\!\!\!\!\!
\lineunder{\lambda(\theta)\varphi_{_{\rm MA}}(\FF)
+ c_{0}\theta(1{-}{\rm ln}\theta)\!}
{$=\COUPLING(\Fe,\theta)^{^{^{}}}$}
\end{align}
with $\varphi_{_{\rm MA}}:=\varphi_{_{\rm M}}{-}\varphi_{_{\rm A}}$; cf.\
\cite[Example 2.5]{MieRou20TKVR}. Here $\lambda:\R\to[0,1]$ denotes the
volume fraction of the austenite versus martensite which is supposed to
depend only on temperature. In this case, the coupling stress
$\COUPLING_{\Fe}'(\Fe,\theta)\Fe^\top\!/\!\det\Fe$ equals to
 $\lambda(\theta)\varphi_{_{\rm MA}}'(\FF)\FF^\top\!/\!\det\FF$
but need not be bounded unless $\FF$ ranges compact sets in ${\rm GL}^+(d)$
 and $\lambda(\cdot)$ is bounded. The heat capacity then reads as
\begin{align*}c(\FF,\theta)=
\frac{ c_{0}+\theta\lambda''(\theta)\varphi_{_{\rm MA}}(\FF)}{\det\FF}\,.
\end{align*}
Again, like in Examples \ref{exa-neo-Hook} and \ref{exa-neo-Hook-alt}, the 
positivity of the heat capacity is a subtle issue since $\lambda(\cdot)$ cannot
be convex and a sufficiently big $c_{0}$ 
to dominate $\theta\lambda''(\theta)\varphi_{_{\rm MA}}(\FF)$ is needed.
\end{example}

\begin{example}[{\sl Volumetric phase transitions}]\upshape\label{exa-PT}
The phase transition in SMA in Remark~\ref{rem-SMA} is primarily isochoric.
In contrast, other phase transitions can be volumetric, which occurs
specifically in rocks which ``compactify'' under compression by big pressures;
e.g.\ in the Earth's silicate mantle there are important phase transitions at
pressure 14\,GPa and 24\,GPa in the depths 440 and 660~km, respectively. The
simplest neo-Hookean ansatz is
\begin{align}\label{psi-volumetric-PT}
\psi(\FF,\theta)=v(\det\FF)+
\frac12G_\text{\sc e}^{}\Big(\frac{{\rm tr}(\Fe\Fe^\top)}{(\det\Fe)^{2/d}}-d\Big)+\phi(\theta)\,.
\end{align}
This gives the pressure depending on $\FF$ as  
$p=\pl\varphi'(\FF)/\pl\det\FF=v'(\det\FF)$. This function $v$ has two
constant parts - two plateaus here with the values 14\,GPa and 24\,GPa.
Then $\varrho$ is a function of pressure $\varrho=\rho(p)$, namely
$\varrho=\rhoR/\!\det\FF=\varrho_0/[v']^{-1}(p)$ and this 
$\rho(\cdot)$ is discontinuous with two jumps at 14\,GPa and 24\,GPa. Making
$v$ dependent also on temperature, one could then model the complete {\it state equation}
of the type $\varrho=\rho(p,\theta)=\varrho_0/[[v(\cdot,\theta)]']^{-1}(p)$
with phase transitions at specific pressures which also depend on temperature,
as used in material science and (geo)physics. Such $v(\det\FF,\theta)$ would
then contribute to the heat capacity by a term
$-\theta v_{\theta\theta}''(\det\FF,\theta)=-\theta v_{\theta\theta}''(\rhoR/\rho(p,\theta),\theta)$
depending discontinuously also on pressure. Additional adiabatic
effects due to volume
change may make this phase transition either exothermic or endothermic.
 Since $\COUPLING_\FF'(\FF,\theta)=\bm0$, the ansatz \eq{psi-volumetric-PT}
easily complies with \eq{Euler-ass}. 
\end{example}

\section{Appendix: transport by Lipschitz velocity fields}\label{sec-app3}
%        ~~~~~~~~~~~~~~~~~~~~~~~~~~~~~~~~~~~~~~~~~~~~~~~~

\def\ZZ{{\bm z}}
\def\bilin{b}

For $n\in\N$ and the vector field $\ZZ$ valued in $\R^n$, it is useful to
state results for an initial-value problem with a general linear transport
and evolution equation:
\begin{align}\label{Z-dynamics}
\DT\ZZ=\bilin(\nabla\vv,\ZZ)\ \ \ \text{ with }\ \ \ \ZZ|_{t=0}^{}=\ZZ_0\,,
\end{align}
where $b:\R^{d\times d}\times\R^n\to\R^n$ is a bilinear mapping. 
We will primarily be interested in weak solutions to \eq{Z-dynamics}
defined by using the integral identity like \eq{one-integral-id}, i.e.\ here
\begin{align}
\int_0^T\!\!\!\int_\varOmega\ZZ{\cdot}\pdt{\widetilde\ZZ}
+\big(({\rm div}\,\vv)\ZZ+\bilin(\nabla\vv,\ZZ)\big){\cdot}\widetilde\ZZ
+\ZZ{\cdot}((\vv{\cdot}\Nabla)\widetilde\ZZ)
\,\d\xx\d t=
-\!\int_\varOmega\!\ZZ_0{\cdot}\widetilde\ZZ(0)\,\d\xx
\label{one-integral-id+}
\end{align}
to be valid for any $\widetilde\ZZ$ smooth with $\widetilde\ZZ(T)={\bm0}$.
The following assertion has been proved essentially in \cite{Roub22QHLS}
or, in a special case in \cite{RouSte22VESS}. We thus present a bit more
general situation here only with sketched proof. 

\begin{lemma}[Evolution-and-transport equation \eq{Z-dynamics}]
\label{lem-transport}
Let $n\in\N$, $p>d$, $r>2$, and 
$\bilin:\R^{d\times d}\times\R^{{n}}\to\R^{{n}}$ is continuous and bilinear.
Then, for any $\vv\in L^1(I;W^{2,p}(\varOmega;\R^d))$ with
$\vv{\cdot}\nn=0$ and any $\ZZ_0\in W^{1,r}(\varOmega;\R^{{n}})$,
there exists exactly one weak solution $\ZZ\in
C_{\rm w}(I;W^{1,r}(\varOmega;\R^{{n}}))\cap W^{1,1}(I;L^r(\varOmega;\R^{{n}}))$
to \eqref{Z-dynamics} and the estimate
\begin{align}\label{F-evol-est}
\|\ZZ\|_{L^\infty(I;W^{1,r}(\varOmega;\R^{{n}}))\,\cap\,
W^{1,1}(I;L^r(\varOmega;\R^{{n}}))}^{}\le
\mathfrak{C}\Big(\|\Nabla\vv\|_{L^1(I;W^{1,p}(\varOmega;\R^{d\times d}))}^{}\,,\,
\|\ZZ_0\|_{W^{1,r}(\varOmega;\R^{{n}})}^{}\Big)
\end{align}
holds with some $\mathfrak{C}\in C(\R^2)$. The equation
\eqref{Z-dynamics} actually holds a.e.\ on $I\times\varOmega$
and $\ZZ\in C(I\times\barOmega;\R^{{n}})$. Moreover, the mapping
\begin{align}\label{v-mapsto-F}
\vv\mapsto\ZZ:L^1(I;W^{2,p}(\varOmega;\R^d))\to
L^\infty(I;W^{1,r}(\varOmega;\R^{{n}}))
\end{align}
is (weak,weak*)-continuous.
\end{lemma}

\begin{proof}
First, for analytical reasons, let us mollify
$\vv$ in time in order to obtain $\vv\in L^2(I;W^{2,p}(\varOmega;\R^d))$.
This will facilitate the (anyhow not uniform) estimate \eq{test-FF+} below
and eventually this mollification can be forgotten.

Let us make a parabolic regularization of \eq{Z-dynamics} by considering
\begin{align}\label{F-evol-reg}
\DT\ZZ=\bilin(\Nabla\vv,\ZZ)+\varepsilon{\rm div}(|\nabla\ZZ|^{r-2}\nabla\ZZ)\,,
\end{align}
with an additional boundary condition $(\Nabla\ZZ)\nn=\bm0$.
Then we make a Faedo-Galerkin approximation of \eq{F-evol-reg} by using a
collection of nested
finite-dimensional subspaces $\{V_k\}_{k\in\N}$ whose union is dense in
$W^{1,p}(\varOmega;\R^{{n}})$. Without loss of generality, we can assume
that $\ZZ_0\in V_1$. Existence of this solution, let us denote
it by $\ZZ_k$, is based on the theory of
systems of ordinary differential equations first locally in time, and then
by successive prolongation on the whole time interval based on 
the $L^\infty$-estimates below.

Testing the Galerkin approximation of \eq{F-evol-reg} by $\ZZ_k$,
we can estimate 
\begin{align}\nonumber
\frac{\d}{\d t}&\int_\varOmega\frac12|\ZZ_k|^2\,\d\xx
+\varepsilon\!\int_\varOmega|\Nabla\ZZ_k|^r\,\d\xx
=\int_\varOmega\Big(
\bilin(\Nabla\vv,\ZZ_k)-(\vv{\cdot}\nabla)\ZZ_k\Big){\cdot}\ZZ_k
\,\d\xx
\\[-.4em]&\ \ 
=\int_\varOmega\!\bilin(\Nabla\vv,\ZZ_k){\cdot}\ZZ_k
+\frac{{\rm div}\,\vv}2|\ZZ_k|^2\,\d\xx
\le\Big(B+\frac12\Big)\|\nabla\vv\|_{L^\infty(\varOmega;\R^{d\times d})}^{}
\|\ZZ_k\|_{L^2(\varOmega;\R^{{n}})}^2
\label{test-FF}\end{align}
with $B:=\sup_{V\in\R^{d\times d},z\in\R^n{},|V|=1,|z|=1}^{}\bilin(V,z)$;
here we used also the calculus 
\begin{align}\nonumber
  \int_\varOmega(\vv{\cdot}\nabla)\ZZ_k{\cdot}\ZZ_k\,\d\xx
  &=\!\int_\varGamma|\ZZ_k|^2(\vv{\cdot}\nn)\,\d S
 %\\[-.4em]\nonumber&\quad
-\!\int_\varOmega\!\ZZ_k{\cdot}(\vv{\cdot}\nabla)\ZZ_k+({\rm div}\,\vv)|\ZZ_k|^2\,\d\xx
=-\frac12\int_\varOmega({\rm div}\,\vv)|\ZZ_k|^2\,\d\xx
\end{align}
together with the boundary condition $\vv{\cdot}\nn=0$.
Note that, beside regularity of $\Nabla\vv$, we needed an integrability
$\vv\in L^1(I{\times}\varOmega;\R^d)$ to ensure legitimacy of the integrals in
\eq{test-FF}. By the Gronwall inequality exploiting the first left-hand-side
term which does not contain the factor $\varepsilon$, we will obtain the
estimate
\begin{align}
\label{Euler-quasistatic-est1-2}
\|\ZZ_k\|_{L^\infty(I;L^2(\varOmega;\R^{{n}}))}\le C\ \ \text{ with }\ \ \|\Nabla\ZZ_k\|_{L^r(I{\times}\varOmega;\R^{d\times {n}})}\le C\varepsilon^{-1/r}\,.
\end{align}

On the Galerkin-discretization level, another legitimate test
 of \eq{F-evol-reg} is by $\pdt{}\ZZ_k$. This allows to estimate 
\begin{align}\nonumber
&\int_\varOmega\bigg|\pdt{\ZZ_k}\bigg|^2\,\d\xx
+\frac\varepsilon r\frac{\d}{\d t}\int_\varOmega|\Nabla\ZZ_k|^r\,\d\xx
=\int_\varOmega\Big(\bilin(\Nabla\vv,\ZZ_k)-(\vv{\cdot}\nabla)\ZZ_k
\Big){\cdot}\pdt{\ZZ_k}\,\d\xx
\\[.2em]&\nonumber\hspace{9em}
\le B^2\|\nabla\vv\|_{L^\infty(\varOmega;\R^{d\times d})}^2\|\ZZ_k\|_{L^2(\varOmega;\R^{{n}})}^2
\\[-.6em]&\hspace{10em}
+C_r\|\vv\|_{L^\infty(\varOmega;\R^d)}^2
\Big(1+\|\Nabla\ZZ_k\|_{L^r(\varOmega;\R^{{n}})}^r\Big)
+\frac12\bigg\|\pdt{\ZZ_k}\bigg\|_{L^2(\varOmega;\R^{{n}})}^2
\label{test-FF+}
\end{align}
with some $C_r\in\R$; here we used that $r>2$ is assumed. Using the
assumed regularity of $\vv$, in particular also
$\vv\in L^2(I;L^\infty(\varOmega;\R^d))$ for which we needed
to mollify temporarily the original $\vv\in L^1(I;W^{2,p}(\varOmega;\R^d))$,
and the already obtained estimate 
\eq{Euler-quasistatic-est1-2} and the Gronwall inequality, we obtain
\begin{align}
\label{Euler-quasistatic-est2-z}
&
\Big\|\pdt{\ZZ_k}\Big\|_{L^2(I{\times}\varOmega;\R^{{n}})}\le 
\varepsilon C{\rm e}^{1/(r\varepsilon)}\ \ \text{ and }\ \ 
\|\Nabla\ZZ_k\|_{L^\infty(I;L^r(\varOmega;\R^{d\times {n}}))}\le C
{\rm e}^{1/(r\varepsilon)}\,.
\end{align}
Note that here the Gronwall inequality uses not the first but the second
left-hand-side term which contains the factor $\varepsilon$ so that both
estimates in \eqref{Euler-quasistatic-est2-z} are $\varepsilon$-dependent.

Considering $\varepsilon>0$ fixed, these estimates allow for the limit passage
with $k\to\infty$ by standard arguments for quasilinear parabolic equations.
The limit is a weak solution to the initial-boundary value problem for
\eq{F-evol-reg}, let us denote it by $\ZZ_\varepsilon\in
L^\infty(I;L^2(\varOmega;\R^{{n}}))\cap L^r(I;W^{1,r}(\varOmega;\R^{{n}}))$.
Actually, as it is determined uniquely, even the whole sequence
$\{\ZZ_k\}_{k\in\N}^{}$ converges to it.

Since now
\begin{align}
\bigg\|\pdt{\ZZ_\varepsilon}+(\vv{\cdot}\nabla)\ZZ_\varepsilon-
\bilin(\nabla\vv,\ZZ_\varepsilon)
\bigg\|_{L^2(I{\times}\varOmega;\R^{{n}})}\!\le\varepsilon C{\rm e}^{1/(r\varepsilon)}\,,\end{align}
by comparison we also obtain
\begin{align}\label{Euler-quasistatic-est3-1}
&\|{\rm div}(|\Nabla\ZZ_\varepsilon|^{r-2}\Nabla\ZZ_\varepsilon^{})\|_{L^2(I{\times}\varOmega;\R^{{n}})}^{}\le
C{\rm e}^{1/(r\varepsilon)}\,.
\end{align}
Although this estimate blows up when $\varepsilon\to0$, we have now at
least the information
that ${\rm div}(|\Nabla\ZZ_\varepsilon|^{r-2}\Nabla\ZZ_\varepsilon)\in
L^2(I{\times}\varOmega;\R^{{n}})$ and we have
\eq{F-evol-reg} continuous. Therefore, we can legitimately
use ${\rm div}(|\Nabla\ZZ_\varepsilon|^{r-2}\Nabla\ZZ_\varepsilon)$ as a test.
Let us denoted the Sobolev exponent to $r$ by $r^*$, i.e.\ $r^*=dr/(d{-}r)$ if
$r<d$ while $r^*=+\infty$ if $r>d$ or arbitrary in $[1,+\infty)$ if $r=d$.
Since $p>d$, we have $p^{-1}+(r^*)^{-1}+(r')^{-1}\le1$, and thus by the H\"older
and Young inequalities, we can estimate
\begin{align}\nonumber
\frac{\d}{\d t}&
\int_\varOmega\frac1r|\nabla\ZZ_\varepsilon|^r\,\d\xx\le\frac{\d}{\d t}
\int_\varOmega\frac1r|\nabla\ZZ_\varepsilon|^r\,\d\xx
+
\varepsilon\int_\varOmega|{\rm div}(|\Nabla\ZZ_\varepsilon|^{r-2}\Nabla\ZZ_\varepsilon)|^2\,\d\xx
\\[-.0em]\nonumber&
=\int_\varOmega\nabla\big((\vv{\cdot}\nabla)\ZZ_\varepsilon-
\bilin(\nabla\vv,\ZZ_\varepsilon)
  \big)\Vdots\big(|\nabla\ZZ_\varepsilon|^{r-2}\nabla\ZZ_\varepsilon\big)\,\d\xx
\\[-.1em]\nonumber&=
\int_\varOmega\Big(|\nabla\ZZ_\varepsilon|^{r-2}(\nabla\ZZ_\varepsilon{\otimes}\nabla\ZZ_\varepsilon){:}\ee(\vv)-\frac1r|\nabla\ZZ_\varepsilon|^r{\rm div}\,\vv
\\[-.4em]&\hspace*{9em}\nonumber
-\big(\bilin_{\nabla\vv}'(\ZZ_\varepsilon)\nabla^2\vv_\varepsilon
{+}\bilin_{\ZZ}'(\nabla\vv)\nabla\ZZ_\varepsilon
\big)\Vdots\big(|\nabla\ZZ_\varepsilon|^{r-2}\nabla\ZZ_\varepsilon\big)\Big)\,\d\xx
\\&\nonumber
\le C_r\|\nabla\vv\|_{L^\infty(\varOmega;\R^{d\times d})}^{}
\|\nabla\ZZ_\varepsilon\|_{L^r(\varOmega;\R^{d\times {n}})}^r\!
+C_r\|\nabla^2\vv\|_{L^p(\varOmega;\R^{d\times d\times d})}^{}
\|\ZZ_\varepsilon\|_{L^{r^*}(\varOmega;\R^{{n}})}^{}
\|\nabla\ZZ_\varepsilon\|_{L^{r}(\varOmega;\R^{d\times {n}})}^{r-1}
\\[.1em]&\le\nonumber
C_r\|\nabla\vv\|_{L^\infty(\varOmega;\R^{d\times d})}^{}
\|\nabla\ZZ_\varepsilon\|_{L^r(\varOmega;\R^{d\times {n}})}^r\!+
C_rN\|\nabla^2\vv\|_{L^p(\varOmega;\R^{d\times d\times d})}^{}
\|\nabla\ZZ_\varepsilon\|_{L^{r}(\varOmega;\R^{d\times {n}})}^r
\\[.1em]&\hspace*{9em}+C_rN\|\nabla^2\vv\|_{L^p(\varOmega;\R^{d\times d\times d})}^{}
\|\ZZ_\varepsilon\|_{L^2(\varOmega;\R^{{n}})}
\big(1{+}\|\nabla\ZZ_\varepsilon\|_{L^{r}(\varOmega;\R^{d\times {n}})}^r\big)\,,
\label{test-Delta-r}
\end{align}
where we used $p>d$ also for the embedding of $\Nabla\vv$
into $L^\infty(\varOmega;\R^{d\times d})$ and where we further used
the calculus 
\begin{align}\nonumber
&\int_\varOmega\nabla\big((\vv{\cdot}\nabla)\ZZ_\varepsilon
  \big){\cdot}|\nabla\ZZ_\varepsilon|^{r-2}\nabla\ZZ_\varepsilon\,\d\xx
 \\[-.5em]&\hspace{2em}\nonumber=\int_\varOmega|\nabla\ZZ_\varepsilon|^{r-2}(\nabla\ZZ_\varepsilon{\otimes}\nabla\ZZ_\varepsilon){:}\ee(\vv)  
+(\vv{\cdot}\nabla)\nabla\ZZ_\varepsilon\Vdots|\nabla\ZZ_\varepsilon|^{r-2}\nabla\ZZ_\varepsilon\,\d\xx
\\[-.1em]&\hspace{2em}\nonumber
=\int_\varGamma|\nabla\ZZ_\varepsilon|^r\vv{\cdot}\nn\,d S
+\int_\varOmega\Big(|\nabla\ZZ_\varepsilon|^{r-2}(\nabla\ZZ_\varepsilon{\otimes}\nabla\ZZ_\varepsilon){:}\ee(\vv)
\\[-.8em]&\hspace{14em}\nonumber
-({\rm div}\,\vv)|\nabla\ZZ_\varepsilon|^r-(r{-}1)|\nabla\ZZ_\varepsilon|^{r-2}\nabla\ZZ_\varepsilon\Vdots
(\vv{\cdot}\nabla)\nabla\ZZ_\varepsilon\Big)\,\d\xx
\\[-.4em]&\hspace{2em}\nonumber
=\int_\varGamma\frac{|\nabla\ZZ_\varepsilon|^r\!\!}r\ \vv{\cdot}\nn\,d S
+\int_\varOmega|\nabla\ZZ_\varepsilon|^{r-2}(\nabla\ZZ_\varepsilon{\otimes}\nabla\ZZ_\varepsilon){:}\ee(\vv)-({\rm div}\,\vv)\frac{|\nabla\ZZ_\varepsilon|^r\!\!}r\ \d\xx\,.
\end{align}
Again, the boundary integral vanishes in \eq{test-Delta-r} if
$\vv{\cdot}\nn=0$. For the last inequality in \eq{test-Delta-r}, we will use
$\|\ZZ_\varepsilon\|_{L^{r^*}(\varOmega;\R^{{n}})}^{}\le
N(\|\ZZ_\varepsilon\|_{L^2(\varOmega;\R^{{n}})}^{}+\|\Nabla\ZZ_\varepsilon\|_{L^{r}(\varOmega;\R^{{n}})}^{})$
where $N$ is the norm of the embedding $W^{1,r}(\varOmega)\subset
L^{r*}(\varOmega)$ if $W^{1,r}(\varOmega)$ is endowed with the norm
$\|\cdot\|_{L^2(\varOmega)}+\|\nabla\cdot\|_{L^r(\varOmega;\R^d)}$.

Thus one can apply the Gronwall inequality to \eq{test-Delta-r}.
The estimates \eq{Euler-quasistatic-est1-2} and 
\eq{Euler-quasistatic-est2-z} can thus be strengthened. Specifically, using
the former estimate in 
\eqref{Euler-quasistatic-est1-2}
and having assumed $\ZZ_0\in W^{1,r}(\varOmega;\R^{{n}})$, one obtains the estimates
\begin{subequations}\label{Euler-quasistatic-est3}
\begin{align}\label{Euler-quasistatic-est3-1+}
&\|\Nabla\ZZ_\varepsilon\|_{L^\infty(I;L^r(\varOmega;\R^{d\times {n}}))}\le C\ \ \text{ and }
\\&\label{Euler-quasistatic-est3-2}
\big\|{\rm div}(|\Nabla\ZZ_\varepsilon|^{r-2}\Nabla\ZZ_\varepsilon^{})
\big\|_{L^2(I{\times}\varOmega;\R^{{n}})}\!\le C\varepsilon^{-1/2}\,.
\end{align}\end{subequations}

The limit passage for $\varepsilon\to0$ in linear terms is then easy and, 
due to \eq{Euler-quasistatic-est3-2} the quasilinear regularizing term
in \eq{F-evol-reg} vanishes as $\mathscr{O}(\varepsilon^{1/2})$ for
$\varepsilon\to0$. Alternatively, when tested it by $\widetilde\ZZ$
and using \eq{Euler-quasistatic-est3-1+}, it vanishes even faster as
$$
\bigg|\int_0^T\!\!\!\int_\varOmega\varepsilon|\Nabla\ZZ_\varepsilon|^{r-2}\Nabla\ZZ_\varepsilon\Vdots
\Nabla\widetilde\ZZ\,\d\xx\d t\bigg|\le
\varepsilon\|\Nabla\ZZ_\varepsilon\|_{L^r(I\times\varOmega;\R^{d\times {n}})}^{r-1}
\|\Nabla\widetilde\ZZ\|_{L^r(I\times\varOmega;\R^{d\times {n}})}^{}=\mathscr{O}(\varepsilon)\,.
$$
In any case, the limit for $\varepsilon\to0$ solves the original initial-boundary value
problem \eq{F-evol-reg}. As this equation is linear, this solution is unique.

The former estimate in \eq{Euler-quasistatic-est2-z} on $\pdt{}\ZZ_\varepsilon$ is
not inherited by the limit, but we can obtain by comparison
$\pdt{}\ZZ=\bilin(\Nabla\vv,\ZZ)-(\vv{\cdot}\Nabla)\ZZ$ at least the estimate 
\begin{align}
\Big\|\pdt{\ZZ}\Big\|_{L^1(I;L^r(\varOmega;\R^{{n}}))}\le C\,.
\end{align}
In particular, the equation in \eq{F-evol-reg} holds a.e.\ on
$I\times\varOmega$. By the embedding $L^\infty(I;W^{1,r}(\varOmega))
\cap W^{1,1}(I;L^r(\varOmega))\subset C(I{\times}\barOmega)$,
we have also $\ZZ\in C(I{\times}\barOmega;\R^{{n}})$.

The (weak,weak*)-continuity of the mapping $\vv\mapsto\ZZ$ as \eq{v-mapsto-F} is
simple is we realize also the bound of $\ZZ$ in
$W^{1,1}(I;L^r(\varOmega;\R^{{n}}))$ and use the Aubin-Lions theorem, which then
gives strong convergence of $\ZZ$'s in
$L^{1/\zeta}(I;L^{r^*-\zeta}(\varOmega;\R^{{n}}))$ for any $0<\zeta\le1$.
Then, if varying $\vv$ in the weak topology of
$L^1(I;W^{2,p}(\varOmega;\R^{{n}}))$, we can pass to the limit in
the bi-linear nonlinearities $(\vv{\cdot}\nabla)\ZZ$ and
$\bilin(\nabla\vv,\ZZ)$ in \eq{Z-dynamics} in its weak formulation
\eq{one-integral-id+}.
%$\hfill\Box$
\end{proof}

In Section~\ref{sec-anal}, we used \eq{Z-dynamics} with $n=d{\times}d$ for
\begin{subequations}\label{usage}\begin{align}\label{usage1}
&\text{$\ \bilin(\nabla\vv,\ZZ)=(\nabla\vv)\ZZ\ \ $ with $\ \ \ZZ=\FF\,$,}
\intertext{which yielded \eq{Euler-thermodynam2} and, in the weak formulation,
\eq{Euler2-weak}. Simultaneously, we used it in the scalar case $n=1$ for}
&\label{usage2}\text{$\ \bilin(\nabla\vv,\ZZ)=-\zz\,{\rm div}\,\vv\ \ $ with $\ \ \ZZ=$}
\begin{cases}\text{$\ \ \ \varrho$,\ \ \ or}&\\
\text{$1/\!\det\FF\,,$}\end{cases}
\intertext{which was applied to \eq{cont-eq+} or \eq{DT-det-1}, respectively.
In the case of the continuity equation \eq{cont-eq+}, the weak formulation
\eq{one-integral-id+} yields the weak formulation \eq{Euler0-weak}; note that
the term $({\rm div}\,\vv)\ZZ+\bilin(\nabla\vv,\ZZ)$ in \eq{one-integral-id+}
vanishes in this last case. Similarly, this scalar case was applied to
\eq{cont-eq-inverse} for}
&\label{usage4}\text{$\ \bilin(\nabla\vv,\ZZ)=\zz{\rm div}\,\vv\ \ $ with
$\ \ \ZZ=1/\varrho\,$.}
\intertext{Another tensorial application with $n=d{\times}d$ could be on
distortion $\FF^{-1}$, relying on:}
\nonumber
&\text{$\ \bilin(\nabla\vv,\ZZ)=-\ZZ(\nabla\vv)\ \ $ with $\ \ \ZZ=\FF^{-1}\ $}
\end{align}\end{subequations}
together with the qualification \eq{Euler-ass-Fe0}, which guarantees
$\nabla\FF_0^{-1}=\nabla({\rm Cof}\FF_0/\!\det\FF_0)$ $\in$
$L^r(\varOmega;\R^{d\times d\times d})$.
A trivial vectorial application with  $n=d$ and $\bilin\equiv0$ is
for the return mapping $\ZZ=\bm{\xi}$, cf.\ \eq{transport-xi}, which can
then lead to a formulation in terms of the distortion $\FF^{-1}=\nabla\bm{\xi}$,
cf.\ in particular \cite{GodPes10TCNM,PaPeKl20HCM}. For another nontrivial
affine and even time-dependent $\bilin$, which would allow here for 
a combination with Maxwellian rheology leading to a more general
 Jeffreys' visco-elastic rheology, we refer to
\cite{Roub22QHLS}. 

% end of small

\end{sloppypar}
\end{document}